%% file: SCBADMM-final.tex
\documentclass[11pt]{article}
\usepackage{listings}

\usepackage{pdflscape}
\usepackage{epsfig}
\usepackage{lscape}
\usepackage{longtable}
\usepackage{amsmath,amssymb}
\usepackage{graphicx}
\usepackage{color}
\usepackage{placeins}
\usepackage{url}
\newcommand{\blue}[1]{\begin{color}{blue}#1\end{color}}

\newcommand{\red}[1]{\begin{color}{red}#1\end{color}}
\newcommand{\green}[1]{\begin{color}{green}#1\end{color}}

\def\sqr#1#2{{\vcenter{\vbox{\hrule height.#2pt
        \hbox{\vrule width.#2pt height#1pt \kern#2pt
        \vrule width.#2pt}
        \hrule height.#2pt}}}}
\newcommand{\qed}{\hfill $ \mathchoice\sqr74\sqr74\sqr{2.1}3\sqr{1.5}3 $}
\def\approxleq{ \kern3pt \mbox{\raisebox{.6ex}{$<$}} \kern-8pt
  \mbox{\raisebox{-.6ex}{$\sim$}} \kern5pt}
\def\mc{\multicolumn}
\def\norm#1{\|#1 \|}
\def\inprod#1#2{\langle#1,\,#2 \rangle}

\def\cA{{\cal A}} \def\cB{{\cal B}}  \def\cK{{\cal K}}
 \def\cD{{\cal D}} \def\cP{{\cal P}} \def\cT{{\cal T}}
\def\cH{{\cal H}} \def\cU{{\cal U}}   \def\cG{{\cal G}} 
\def\cM{{\cal M}}  \def\cW{{\cal W}}
\def\cS{{\cal S}} \def\cE{{\cal E}}  \def\cF{{\cal F}}  \def\cV{{\cal V}}
\def\cI{{\cal I}} \def\cX{{\cal X}} \def\cY{{\cal Y}} \def\cZ{{\cal Z}}

\def\cQ{{\cal Q}}             
\def\R{{\sl I\kern-4.3pt R}}  \def\Rn{\R^n}      

\def\Rn2{\R^{n(n+1)/2}}
\def\Sn{{\cal S}^n}

\def\skron{{{\kern2.3pt\bigcirc\kern-8.5pt\ast\kern5pt }}}
\def\sskron{{{\kern2.3pt\bigcirc\kern-11.5pt\ast\kern7pt }}}
\def\cirdot{{{\kern2.3pt\bigcirc\kern-6.8pt\cdot\kern5pt}}}
 
\def\C{\kern4pt\hbox{\rm\vrule height 6.9pt width 0.7pt\kern-3.5pt C}}

\def\sig{\sigma}

\def\nn{\nonumber}

\def\disp{\displaystyle}

\def\ubar{\bar{u}}
\def\vbar{\bar{v}}  \def\vprime{v^\prime}
\def\xbar{\bar{x}} \def\ybar{\bar{y}} \def\zbar{\bar{z}}
\def\cbar{\bar{c}}
\def\deltabar{\bar{\delta}} \def\gammabar{\bar{\gamma}}
\def\alphabar{\bar{\alpha}} \def\betabar{{\beta}}

\def\yprime{y^\prime}
\def\Xprime{X^\prime}

\def\cX{{\cal X}} 
\def\cU{{\cal U}} 
\def\cR{{\cal R}}  \def\cL{{\cal L}}

\newlength{\len}
\newtheorem{theorem}{Theorem}[section]
\newtheorem{prop}{Proposition}[section]

\newtheorem{remark}{Remark}[section]
\newtheorem{assumption}{Assumption}[section]

\def\SCBADMM{\mbox{\sc Scb-spadmm}}
\def\ADMM{\mbox{\sc Admm}}
\def\ADMMGB{\mbox{\sc Admmgb}}

\begin{document}

\title{\bf A Schur Complement Based Semi-Proximal ADMM for Convex Quadratic Conic Programming
and Extensions}
%\author{Xudong Li, Defeng Sun and Kim-Chuan Toh\\ National University of Singapore}
\author{Xudong Li\thanks{Department of Mathematics, National University of Singapore, 10 Lower Kent Ridge Road, Singapore ({\tt lixudong@nus.edu.sg}).}, \;
Defeng Sun\thanks{Department  of  Mathematics  and  Risk  Management  Institute, National University of Singapore, 10 Lower Kent Ridge Road, Singapore ({\tt matsundf@nus.edu.sg}). } \  and \
Kim-Chuan Toh\thanks{Department of Mathematics, National University of Singapore, 10 Lower Kent Ridge Road, Singapore
({\tt mattohkc@nus.edu.sg}).  }
}

\date{September 9, 2014}

\maketitle
\begin{abstract} This paper is devoted to the design of an efficient and convergent  {semi-proximal} alternating direction method of multipliers (ADMM) for finding a solution of low to medium accuracy to convex quadratic conic programming and related problems. For this class of problems, the convergent two block semi-proximal ADMM can be employed to solve their primal form  in a straightforward way.  However, it is known  that it is more efficient to apply the directly extended multi-block semi-proximal ADMM,  though its convergence is not guaranteed,  to  the dual form  of these problems. Naturally, one may ask the following question: can one construct a convergent multi-block semi-proximal ADMM that is more efficient than the directly extended semi-proximal ADMM? Indeed, for linear conic programming with 4-block constraints this has been shown to be  achievable in a recent paper by Sun, Toh and Yang [arXiv preprint
arXiv:1404.5378, (2014)]. Inspired by the
aforementioned work and with the convex quadratic conic programming  in mind,
 we propose a Schur complement based convergent semi-proximal ADMM  for solving convex programming problems, with a coupling linear equality constraint,  whose objective function is the sum of two  proper closed  convex functions plus an arbitrary number of convex  quadratic or linear functions.
 Our convergent
semi-proximal ADMM is particularly suitable for solving convex quadratic semidefinite programming (QSDP) with constraints consisting of
linear equalities, a positive semidefinite cone and a simple convex polyhedral set.
 The  efficiency of our proposed algorithm is  demonstrated by    numerical experiments on various examples including QSDP.

 \vskip 15 true pt

\end{abstract}

\section{Introduction}
In this paper, we aim to design an efficient yet simple first order convergent method
for solving convex quadratic conic programming.
An important special case  is the following {convex quadratic semidefinite programming} (QSDP)
 \begin{eqnarray}
  \begin{array}{ll}
    \min &  \frac{1}{2} \inprod{X}{\cQ X} + \inprod{C}{X}  \\[5pt]
   \mbox{s.t.}
       &\cA_E X   =  b_E,
       \quad \cA_I X    \geq b_I, \quad
       X \in \Sn_+\cap \cK  \, ,
\end{array}
 \label{eq-qsdp}
\end{eqnarray}
where  $\cS_+^n$ is the cone of $n\times n$ symmetric and positive semi-definite matrices in the space of $n\times n$ symmetric matrices $\cS^n$ endowed with the standard trace
inner product $\inprod{\cdot}{\cdot}$ and the Frobenius norm $\norm{\cdot}$, $\cQ$ is a self-adjoint positive semidefinite linear operator from  $\Sn$ to $\Sn$,  $\cA_E:\Sn \rightarrow \Re^{m_E}$ and $\cA_I:\Sn\rightarrow \Re^{m_I}$ are two linear maps,  $C\in \Sn$,  $b_E\in \Re^{m_E}$ and $b_I\in \Re^{m_I}$ are given data,  $\cK$ is a nonempty simple closed convex set,
e.g., $\cK =\{W\in\Sn:\; L\leq W\leq U\}$ with $L,U\in \Sn$ being given matrices.
By introducing a slack variable $W\in \Sn$, we can equivalently recast (\ref{eq-qsdp}) as
  \begin{eqnarray}
  \begin{array}{lll}
    \min &  \frac{1}{2} \inprod{X}{\cQ X} + \inprod{C}{X} +\delta_{\cK}(W)\\[5pt]
   \mbox{s.t.}
       &\cA_E X     =  b_E, \quad \cA_I X    \geq b_I, \quad
       X = W, \quad X \in \Sn_+\, ,
\end{array}
 \label{eq-qsdp1}
\end{eqnarray}
where $\delta_\cK(\cdot)$ is the indicator function of $\cK$, i.e., $\delta_\cK(X) = 0$ if $X \in \cK$ and $\delta_\cK(X) = \infty$ if $X \notin \cK$.
The dual of problem  (\ref{eq-qsdp1}) is given by
 \begin{eqnarray}
  \begin{array}{rllll}
    \max &   -\delta_{\cK}^*(-Z) + \inprod{b_I}{y_I} - \frac{1}{2}\inprod{X}{\cQ X} + \inprod{b_E}{y_E}    \\[5pt]
   \mbox{s.t.} & Z + \cA_I^*y_I - \cQ X + S + \cA_E^* y_E  = C, \quad
                y_I\geq 0, \quad S \in \Sn_+\, ,
   \end{array}
   \label{eq-d-qsdp}
\end{eqnarray}
where for any $Z\in \Sn$, $ \delta_{\cK}^*(-Z) $ is given by
\begin{eqnarray}
  \delta_{\cK}^*(-Z) = -\inf_{W \in \cK} \inprod{Z}{W}
       = \sup_{W \in \cK} \inprod{-Z}{W}. \label{gz}
       %&=& \sum_{i=1}^n\sum_{j=1}^n \min\{l_{ij}z_{ij},u_{ij}z_{ij}\} \nn.
\end{eqnarray}
%Recently, convex quadratic semidefinite programming has been applied in the matrix completion with fixed basis coefficients \cite{miao2013matrix,wu2014high-dim}.
%Actually, in order to estimate the correlation matrix with missing observations, Wu
%in \cite{wu2014high-dim} has applied the following least squares SDP
%\begin{eqnarray}
%   \begin{array}{lll}
%    \min & \frac{1}{2} \norm{\cB(x) - d}^2 + \inprod{c}{x}  \\[5pt]
%   \mbox{s.t.}
%       &\cA_E x     =  b_E, \quad \cA_I x    \geq b_I, \quad
%        x \in \Sn_+
%\end{array}
% \label{eq-lssdp-wu}
%\end{eqnarray}
%with its dual given by
% \begin{eqnarray}
%  \begin{array}{rllll}
%    \max & \inprod{b_I}{y_I}  +\inprod{d}{v} - \frac{1}{2}\inprod{v}{v} + \inprod{b_E}{y_E}    \\[5pt]
%   \mbox{s.t.} & \cA_I^*y_I  + s + \cB^*v + \cA_E^* y_E  = c, \quad
%                y_I\geq 0, \quad s \in \Sn_+.
%   \end{array}
%   \label{eq-dd-lssdp-wu}
%\end{eqnarray}

 It is evident that the dual problem \eqref{eq-d-qsdp} is in the form of  the following convex optimization model:
\begin{eqnarray}
  \begin{array}{rlll}
    \min & f(u)+\sum_{i=1}^p \theta_i(y_i) + g(v) + \sum_{j=1}^q \varphi_j(z_j)    \\[5pt]
   \mbox{s.t.} & \cF^*u + \sum_{i=1}^p \cA_i^*y_i  + \cG^* v + \sum_{j=1}^q \cB_j^* z_j  = c, \\[5pt]
   \end{array}
   \label{eq-dd-M1}
\end{eqnarray}
where $p$ and $q$ are given nonnegative integers, $f:\cU \rightarrow (-\infty,+\infty]$, $g:\cV \rightarrow (-\infty,+\infty],$ $\theta_i:\cY_i \rightarrow (-\infty,+\infty],\,  i=1,\ldots,p,$  and  $\varphi_j:\cZ_j \rightarrow (-\infty,+\infty],\, j=1,\ldots,q$  are closed proper convex functions, $\cF:\cX\rightarrow \cU,$ $\cG: \cX\rightarrow \cV$, $\cA_i:\cX\rightarrow \cY_i, \,  i=1,\ldots,p$  and $\cB_j: \cX\rightarrow \cZ_j, \, j=1,\ldots,q$  are linear maps, $\cU,\cV, \cY_1, \ldots, \cY_p, \cZ_1, \ldots, \cZ_q $ and $\cX$ are all real finite dimensional Euclidean spaces each equipped with an inner product $\inprod{\cdot}{\cdot}$ and its induced norm $\norm{\cdot}.$

In this paper, we make the following blanket  assumption.
\begin{assumption}\label{main-assumption}
  %While convex $f$ and $g$ are potentially non-quadratic convex functions,
  For $i=1,\ldots,p$ and $j = 1,\ldots,q$, each $\theta_i(\cdot)$ and $\varphi_j(\cdot)$ are convex quadratic functions.
\end{assumption}

Note that,  in general,
  problem \eqref{eq-d-qsdp} does not satisfy Assumption \ref{main-assumption} unless  $y_I$  is vacuous from the model or $\cK\equiv \Sn$. However,  one can always reformulate  problem \eqref{eq-d-qsdp}  equivalently as
  \begin{eqnarray}
  \begin{array}{rllll}
    \min & (\delta^*_{\cK}(-Z) + \delta_{\Re^{m_I}_+}(u)) -\inprod{b_I}{y_I} + \frac{1}{2}\inprod{X}{\cQ X}  +\delta_{\cS^n_+}(S) - \inprod{b_E}{y_E}    \\[5pt]
   \mbox{s.t.} & Z + \cA_I^*y_I - \cQ X  + S  + \cA_E^* y_E  = C,  \\[5pt]
               & \cD^* u  - \cD^* y_I     = 0,
   \end{array}
   \label{eq-d-qsdp-r}
\end{eqnarray}
where $\cD:\Re^{m_I} \to \Re^{m_I}$  is any given nonsingular linear operator and  $\delta_{\Re^{m_I}_+}(\cdot)$ is the indicator function over $\Re^{m_I}_+$.   Now, one can see that problem (\ref{eq-d-qsdp-r}) satisfies
Assumption \ref{main-assumption}.

%
%\blue{ This can be regarded as an extension of our model \eqref{eq-dd-M1}.} Apart from problem \eqref{eq-d-qsdp} and \eqref{eq-rpca-g},
%our convex model \eqref{eq-dd-M1} actually covers a wide range of interesting examples
%and we will mention some of them in our examples section.
There are many other important cases that take  the form of model (\ref{eq-dd-M1}) satisfying Assumption  \ref{main-assumption}. One prominent example comes from the matrix completion with fixed basis coefficients \cite{miao2014arank,miao2013matrix,wu2014high-dim}. Indeed the nuclear semi-norm penalized least squares model in \cite{miao2013matrix} can be written
as
\begin{eqnarray}
   \begin{array}{lll}
    \min\limits_{X \in\Re^{m \times n}} &  \frac{1}{2} \norm{\cA_F X - d}^2 + \rho(\norm{X}_*- \inprod{C}{X}) \\[5pt]
   \mbox{s.t.}
       &\cA_E X     =  b_E, \quad X\in\cK := \{X\mid\norm{\cR_{\Omega}X}_{\infty} \leq \alpha\},
\end{array}
 \label{eq-lsn-miao}
\end{eqnarray}
where $\norm{X}_*$ is the nuclear norm of $X$ defined as the sum of all its singular values,  $\norm{\cdot}_{\infty}$ is the elementwise $l_{\infty}$ norm defined by $\norm{X}_{\infty} := \max_{i=1,\ldots,m}\{\max_{j=1,\ldots,n} |X_{ij}|\}$, $\cA_F : \Re^{m \times n} \to \Re^{n_F}$ and $\cA_E : \Re^{m \times n} \to \Re^{n_E}$ are two linear maps, $\rho$ and $\alpha$ are two given positive parameters,  $d\in \Re^{n_F}$, $C\in \Re^{m \times n}$ and $b_E \in \Re^{n_E}$ are given data, $\Omega\subseteq \{1,\ldots,m\}\times\{1,\ldots,n\}$ is the set of the indices relative to which the basis coefficients are not fixed, $\cR_{\Omega}:\Re^{m \times n} \to \Re^{|\Omega|}$ is the linear map such that $\cR_{\Omega} X := (X_{ij})_{ij\in \Omega}.$
   Note that when there are no fixed basis coefficients (i.e., $\Omega = \{1,\ldots,m\}\times\{1,\ldots,n\}$ and $\cA_E$ are vacuous),  the above problem reduces to the model considered
  by Negahban and Wainwright in \cite{negahban2012restricted} and Klopp in \cite{klopp2014noisy}.
By introducing slack variables $\eta$, $R$ and $W$, we can reformulate problem \eqref{eq-lsn-miao} as
\begin{eqnarray}
   \begin{array}{lll}
    \min &  \frac{1}{2} \norm{\eta}^2   + \rho\big(\norm{R}_* -  \inprod{C}{X}\big) +\delta_{\cK}(W) \\[5pt]
   \mbox{s.t.}
       &\cA_F X - d = \eta, \quad \cA_E x =  b_E, \quad X = R, \quad X = W.
\end{array}
 \label{eq-lsn-miao-r}
\end{eqnarray}
The dual of problem  (\ref{eq-lsn-miao-r}) takes the form of
 \begin{eqnarray}
  \begin{array}{rllll}
    \max &  -\delta_{\cK}^*(-Z) - \frac{1}{2}\norm{\xi}^2 + \inprod{d}{\xi} + \inprod{b_E}{y_E}    \\[5pt]
   \mbox{s.t.} &  Z  + \cA_F^* \xi + S + \cA_E^* y_E  = -\rho C, \quad \norm{S}_2\le \rho,
   \end{array}
   \label{eq-dlsn-miao-r}
\end{eqnarray}
where $\norm{S}_2$ is the operator norm of $S$, which is defined to be its largest singular value.
%where for $z\in \Re^{m \times n}$, $ \delta_{\Omega}^*(-z) $ is given by
%\begin{eqnarray}
%  \delta_{\Omega_\alpha}^*(-z) = -\inf_{w_2 \in \Omega_\alpha} \inprod{z}{w_2}
%       = \sup_{w_2 \in \Omega_\alpha} \inprod{-z}{w_2}. \label{miao-gz}
%       %&=& \sum_{i=1}^n\sum_{j=1}^n \min\{l_{ij}z_{ij},u_{ij}z_{ij}\} \nn.
%\end{eqnarray}

Another compelling example is the so called {robust PCA} (principle component analysis) considered in \cite{wright2009robust}:
\begin{eqnarray}
  \begin{array}{ll}
    \min & \norm{A}_* + \lambda_1\norm{E}_1 + \disp\frac{\lambda_2}{2}\norm{Z}^2_{F} \\[5pt]
   \mbox{s.t.}
       & A + E + Z = W, \quad A,E,Z\in \Re^{m\times n}\, ,
       %\quad \cA_I x    \geq b_I, \quad
%       x \in \Sn_+\cap \Omega  \, ,
\end{array}
 \label{eq-rpca}
\end{eqnarray}
where $W \in \Re^{m\times n}$ is the observed data matrix, $\norm{\cdot}_1$ is the elementwise $l_1$ norm given by $\norm{E}_1 := \sum_{i=1}^m\sum_{j=1}^n |E_{ij}|$, $\norm{\cdot}_F$ is the Frobenius norm, $\lambda_1$ and $\lambda_2$ are two
positive parameters.  There are many different variants to the robust PCA model. For example,
 one may consider the following model where  the observed data matrix $W$ is incomplete:
\begin{eqnarray}
  \begin{array}{ll}
    \min & \norm{A}_* + \lambda_1\norm{E}_1 + \disp\frac{\lambda_2}{2}\norm{\cP_{\Omega}(Z)}^2_{F} \\[5pt]
   \mbox{s.t.}
       & \cP_{\Omega}(A + E + Z) = \cP_{\Omega}(W),\quad A,E,Z\in \Re^{m\times n}\, ,
       %\cP_{\Omega} (A + E)+Z= \cP_{\Omega} (W),
       %\quad \cA_I x    \geq b_I, \quad
%       x \in \Sn_+\cap \Omega  \, ,
\end{array}
 \label{eq-rpca-g}
\end{eqnarray}
%\fbox{\blue{The equality constraint above should involve $\cP_\Omega$?? And the objective should
%be $\cP_{\Omega^c}(Z)$?} }
 i.e. one assumes that only a subset $\Omega \subseteq \{1,\ldots,m\}\times\{1,\ldots,n\}$ of the entries of $W$ can be observed. Here $\cP_\Omega: \Re^{m\times n} \to \Re^{m\times n}$ is the orthogonal projection operator  defined by
\begin{eqnarray}
  \label{P-omega}
  \cP_{\Omega}(X) = \left\{\begin{array}{llll}
    X_{ij} \quad &\textup{if}\; (i,j) \in \Omega, \\[5pt]
    0 \quad &\textup{otherwise}.
  \end{array}\right.
\end{eqnarray}
Again, problem \eqref{eq-rpca-g} satisfies Assumption  \ref{main-assumption}. In \cite{tao2011recovering}, Tao and Yuan tested one of the equivalent forms of  problem \eqref{eq-rpca-g}. In the numerical section, we will see other interesting examples.

%Here, we focus on the special case that the functions  $\theta_i$ and $\varphi_j$ are all
%quadratic functions of the following form
%\begin{eqnarray*}
%  \theta_i(y_i) = \frac{1}{2}\inprod{y_i}{\cP_i y_i} -\inprod{b_i}{y_i} ,  \;i=1,
%  \ldots,p, \qquad \varphi_j(z_j) = \frac{1}{2}\inprod{z_j}{\cQ_j z_j}-\inprod{d_j}{z_j}, \;j=1,\ldots,q,
%\end{eqnarray*}
%where each $\cP_i$ and $\cQ_j$ are given self-adjoint positive semidefinite linear operators.

For notational convenience,  let $\cY := \cY_1\times\cY_2\times,\ldots,\cY_p,$ $\cZ:=\cZ_1\times\cZ_2\times,\ldots,\cZ_q$.
 We write $y \equiv (y_1,y_2,\ldots,y_p)\in \cY $ and $z \equiv (z_1,z_2,\ldots,z_q)\in\cZ$. Define the linear map $\cA:\cX \rightarrow \cY$ such that its adjoint is given by \[\cA^*y = \sum_{i=1}^p\cA_i^*y_i\quad \forall y\in \cY.\]
Similarly, we
  define the linear map $\cB:\cX\rightarrow \cZ$ such that
  its adjoint is given by \[
  \cB^*z = \sum_{j=1}^q\cB_j^*z_j\quad \forall z\in \cZ.\]
  Additionally,  let $\theta(y):=\sum_{i=1}^p \theta_i(y_i),$ $y\in \cY$ and
  $\varphi(z) := \sum_{j=1}^q \varphi_j(z_j)$, $z\in \cZ$.
Now we can rewrite (\ref{eq-dd-M1}) in the following  compact form:
\begin{eqnarray}
  \begin{array}{rlll}
    \min & f(u)+\theta(y)+g(v)+\varphi(z)    \\[5pt]
   \mbox{s.t.} & \cF^*u + \cA^* y  + \cG^* v + \cB^* z  = c. \\[5pt]
   \end{array}
   \label{eq-dd-4}
\end{eqnarray}

Problem \eqref{eq-dd-M1} can be view as a special case of the following block-separable convex optimization problem:
\begin{eqnarray}
 \min \left\{ {{\sum}}_{i=1}^n\phi_i(w_i)  \mid {\sum}_{i=1}^n\cH^*_i w_i = c \right\},
  \label{eq-ADMMq}
\end{eqnarray}
where for each $i\in \{1, \ldots, n\}$, $ \cW_i$ is a  finite dimensional real Euclidean space equipped with an inner product $\inprod{\cdot}{\cdot}$ and its induced norm $\norm{\cdot}$,  $\phi_i: \cW_i \to (-\infty,+\infty]$ is a closed proper convex function, $\cH_i:   \cX \to \cW_i$ is a  linear map and  $c \in \cX$ is given. Note that when we rewrite problem \eqref{eq-dd-M1} in terms of \eqref{eq-ADMMq}, the quadratic structure in \eqref{eq-dd-M1} is hidden in the sense that each $\phi_i$ will be treated equally. However, this special quadratic structure will be thoroughly exploited in our search for an efficient yet  simple   ADMM-type method
with guaranteed convergence.

Let $\sigma > 0$ be a given parameter. The augmented Lagrangian function for (\ref{eq-ADMMq}) is defined by
\begin{eqnarray}
\cL_\sigma(w_1,\ldots,w_n;x): = {\sum}_{i=1}^n\phi_i(w_i) +  \inprod{x}{ {\sum}_{i=1}^n\cH^*_i w_i - c}+ \frac{\sigma}{2} \norm{ {\sum}_{i=1}^n\cH^*_i w_i  - c}^2\nonumber
\label{eq-ADM-q-AL}
\end{eqnarray}
for $w_i\in \cW_i$, $i=1, \ldots, n$ and $x\in \cX.$ Choose any initial points $w_i^0 \in {\rm dom}(\phi_i) $, $i=1, \ldots, q$ and $x^0 \in \cX$. The classical augmented Lagrangian  method consists of the following iterations:
\begin{eqnarray}
 (w_1^{k+1}, \ldots, w_n^{k+1}) &= & \textup{argmin}
\;\cL_\sigma(w_1,\ldots,w_n;x^k),
\label{eq-auglagl-w}
\\[5pt]
 x^{k+1} &=&\displaystyle  x^k +\tau \sigma\left ({\sum}_{i=1}^n\cH^*_i w_i^{k+1}  -c\right),
\label{eq-auglagl-x}
\end{eqnarray}
where $ \tau\in(0,2)$ guarantees the convergence. Due to the non-separability of the quadratic penalty term in $\cL_{\sigma}$, it is generally a challenging task to solve the joint minimization problem \eqref{eq-auglagl-w} exactly or approximately with high accuracy. To overcome this difficulty, one may consider the following
$n$-block alternating direction methods of multipliers (ADMM):

\begin{eqnarray}\label{admmlb}
 w_1^{k+1} &=& \textup{argmin}
  \;\cL_\sigma(w_1,w_2^k\ldots,w_n^k;x^k), \nonumber
\\[5pt]
& \vdots & \nonumber
\\[5pt]
 w_i^{k+1} &=& \textup{argmin} \;\cL_\sigma(w_1^{k+1},\ldots, w_{i-1}^{k+1}, w_i, w_{i+1}^k, \ldots, w^k_n;x^k), \nonumber
 \\[5pt]
& \vdots &\\[5pt] \nonumber
  w_n^{k+1} &=& \textup{argmin} \; \cL_\sigma(w_1^{k+1},\ldots,   w_{n-1}^{k+1},  w_n;x^k), \nonumber
\\[5pt]
 x^{k+1} &=& \displaystyle  x^k +\tau \sigma\left ({{\sum}_{i=1}^n}\cH^*_i w_i^{k+1}  -c\right).
 \nonumber
\end{eqnarray}
The above $n$-block ADMM   is  an direct extension of the ADMM
 for solving the following
  $2$-block convex optimization problem
\begin{eqnarray}
 \min \left\{ \phi_1(w_1)+\phi_2(w_2) \mid \cH_1^* w_1 + \cH^*_2 w_2   = c \right\}.
  \label{eq-ADMM2}
\end{eqnarray}
The convergence of $2$-block ADMM has already been extensively studied in \cite{glowinski1975sur,gabay1976dual,Glowinski1980lectures,fortin1983augmentedlag,Gabay1983299,
eckstein1992douglas} and references therein.
 %Endowed with simplicity and efficiency, recently the $2$-block ADMM has regained its popularity in several applications such as mathematical %imaging science, signal processing and machine learning.
%Despite the outstanding practical performance,
However, the convergence of the $n$-block ADMM has been
ambiguous for a long time.
Fortunately this ambiguity has been addressed very recently in \cite{chen2013direct} where
%The above ambiguity is not well addressed until very recently,
Chen, He, Ye, and Yuan  showed that the direct extension of the ADMM to the case of a $3$-block convex optimization problem is not necessarily convergent. On the other hand, the $n$-block ADMM with $\tau \ge 1$ often works very well in practice and this fact poses a big challenge if one attempts to develop new ADMM-type algorithms which have convergence guarantee but with competitive numerical efficiency and iteration simplicity as the $n$-block ADMM.
%There is always a dream to design a simple and efficient yet
%convergent ADMM-type method for multi-block case.

Recently, there is exciting progress in this active research area.
Sun, Toh and Yang \cite{styang2014} proposed a convergent semi-proximal ADMM (PADMM3c) for convex programming problems of three separable blocks in the objective function with the third part being linear. One distinctive
feature of algorithm PADMM3c is that it requires only  an inexpensive extra step, compared to the 3-block ADMM,  but yields a convergent and faster algorithm. Extensive numerical tests on the doubly non-negative SDP problems with equality and/or inequality constraints demonstrate that PADMM3c can have superior numerical efficiency over the directly extended ADMM. This opens up the possibility of designing an efficient and convergent ADMM type
method for  solving multi-block convex optimization problems.
% {However, there are still some limitations \red{to overcome when PADMM3c is applied}   to solve quadratic conic programming problems \eqref{eq-d-qsdp-r}, \eqref{eq-dlsn-miao-r} and \eqref{eq-rpca-g} as well as our general convex optimization model \eqref{eq-dd-M1}. In fact, in order to apply PADMM3c to the aforementioned problems, one \red{needs} to introduce additional slack variables to the underlying problems and group all the unknown variables wisely. Though properly scaled, these auxiliary variables
%may still slow down the convergence of the resulting algorithm. Furthermore, the need to explicitly  factorize a possibly large dimensional symmetric positive definite matrix may be too restrictive in some situations.}
Inspired by the
aforementioned work,  in this paper we shall propose a Schur complement based semi-proximal ADMM (SCB-SPADMM) to efficiently solve the convex  quadratic conic programming problems to medium accuracy. The development of our algorithm  is  based on the simple yet elegant idea of the Schur complement and the convenient  convergence results of the semi-proximal ADMM given in the appendix of  \cite{fazel2013hankel}. Our primary motivation for designing the proposed SCB-SPADMM is to generate a good initial point quickly to warm-start locally fast convergent method such as the semismooth Newton-CG method used in \cite{zhao2010newton,ystoh2014} for solving linear SDP though the method proposed here is definitely of its own interest.

The remaining parts of this paper are organized as follows. In the next section, we  present a Schur complement based semi-proximal augmented Lagrangian method (SCB-SPALM) to solve a 2-block convex optimization problem where the second function $g$ is quadratic and then show the relation between our SCB-SPALM and the generic 2-block semi-proximal ADMM (SPADMM). In section 3, we
propose our main algorithm SCB-SPADMM for solving the general convex model \eqref{eq-dd-M1}. Our main convergence results are  presented in this section. Section 4 is devoted to the
implementation and numerical experiments of using our SCB-SPADMM to solve convex quadratic conic programming problems and the various extensions. We conclude our paper in the final section.

\medskip
{\bf Notation.} Define the spectral (or operator) norm of a given linear operator $\cT$ by $\norm{\cT} := \sup_{\norm{w} =1} \norm{\cT w}.$ For any $w \in \cU,$ we let
$$\textup{Prox}_f(w) := \textup{argmin}_{u}\; f(u) + \frac{1}{2}\norm{u - w}^2.$$

%%%%%%%%%%%%%%%%%%%%%%%%%%%%%%
\section{A Schur complement based  semi-proximal augmented Lagrangian method}
Before we introduce our approach for the multi-block case, we need to  consider the convex optimization problem with the following 2-block separable structure
\begin{eqnarray}
  \begin{array}{rlll}
    \min & f(u)+g(v)    \\[5pt]
   \mbox{s.t.} & \cF^*u  + \cG^* v = c, \\[5pt]
   \end{array}
   \label{eq-dd-2}
\end{eqnarray}
where $f:\cU \rightarrow (-\infty,+\infty]$ and $g:\cV \rightarrow (-\infty,+\infty]$ are closed proper convex functions, $\cF:\cX\rightarrow \cU$ and $\cG: \cX\rightarrow \cV$ are given
linear maps. The dual of problem (\ref{eq-dd-2}) is given by
\begin{eqnarray}
   \min \left\{\inprod{c}{x} + f^*(s)+g^*(t)\mid  \cF x + s=0, \; \cG x + t=0 \right\}.
  \label{eq-ADMM2-nonlinear-primal}
\end{eqnarray}
Let $\sigma > 0$ be given. The augmented Lagrangian function associated with (\ref{eq-dd-2}) is given as follows:
\begin{eqnarray}
\cL_{\sigma}(u,v;x) &=& f(u) + g(v)
+\inprod{x}{\cF^*u +\cG^*v-c}
+ \frac{\sigma}{2}\norm{\cF^*u+\cG^*v-c}^2.
\label{ALM-2}
\end{eqnarray}

The semi-proximal ADMM proposed in \cite{fazel2013hankel}, when  applied to (\ref{eq-dd-2}),  has the following template. Since the proximal terms added here are allowed to be positive semidefinite, the corresponding method is  referred to  as semi-proximal ADMM instead of proximal ADMM as in \cite{fazel2013hankel}.

\centerline{\fbox{\parbox{\textwidth}{
{\bf Algorithm SPADMM: \bf{A generic 2-block semi-proximal ADMM for solving (\ref{eq-dd-2}).}}
\\
Let $\sigma >0$ and $\tau\in(0,\infty)$ be given parameters. Let $\cT_f$ and
$\cT_g$ be given self-adjoint positive semidefinite, not necessarily positive definite, linear operators defined on $\cU$ and $\cV$, respectively. Choose $(u^0,v^0,x^0)\in\mbox{dom}(f)\times\mbox{dom}(g)\times\cX.$ For $k=0,1,2,...$, perform the $k$th iteration as follows:
\begin{description}
\item [Step 1.] Compute
  \begin{eqnarray}
     u^{k+1} =\mbox{argmin}_{u} \;\cL_{\sigma}(u,v^k;x^k) + \frac{\sigma}{2}\norm{u - u^k}_{\cT_f}^2.
     \label{padm2u}
  \end{eqnarray}
\item [Step 2.] Compute
  \begin{eqnarray}
     v^{k+1} =\mbox{argmin}_{v} \;\cL_{\sigma}(u^{k+1},v;x^k) + \frac{\sigma}{2}\norm{v - v^k}_{\cT_g}^2.
     \label{padm2v}
  \end{eqnarray}
\item [Step 3.] Compute
  \begin{eqnarray}
     x^{k+1} = x^k + \tau\sigma(\cF^*u^{k+1} + \cG^*v^{k+1} - c).
     \label{padm2x}
  \end{eqnarray}
\end{description}
}}}
\bigskip
In the above 2-block semi-proximal ADMM for solving  \eqref{eq-dd-2}, the presence of $\cT_f$ and $\cT_g$ can help to guarantee the existence of solutions for the subproblems \eqref{padm2u} and \eqref{padm2v}.
In addition, they play important roles in ensuring the  boundedness  of the two generated sequences $\{y^{k+1}\}$ and $\{z^{k+1}\}$. Hence, these two proximal terms are preferred. The choices of $\cT_f$ and $\cT_g$ are very much problem dependent. The general principle is that both $\cT_f$ and $\cT_g$ should be as small as possible while $y^{k+1}$ and $z^{k+1}$ are still relatively easy to compute.

Let $\partial f$ and $\partial g$ be the subdifferential mappings of $f$ and $g$, respectively. Since both  $\partial f$ and $\partial g$  are maximally monotone, there exist two self-adjoint and positive semidefinite operators  $\Sigma_f$ and $\Sigma_g $    such that  for all $u,  \tilde{u}\in {\rm dom}(f)$, $\xi\in \partial f(u)$,
and $\tilde{\xi}\in \partial f(\tilde{u})$,
\begin{equation}\label{monosub1}
 \langle \xi-\tilde{\xi}, u -\tilde{u}\rangle \geq
 \|u -\tilde{u}\|^2_{\Sigma_f}
\end{equation}
and for all $v,  \tilde{v} \in {\rm {dom}}(g)$, $\zeta\in \partial
g(v )$, and $ \tilde{\zeta}\in \partial g(\tilde{v})$,
\begin{equation}\label{monosub2}
 \langle \zeta-\tilde{\zeta}, v -\tilde{v}\rangle \geq
 \|v -\tilde{v}\|^2_{\Sigma_g}.
\end{equation}

For the convergence of the 2-block semi-proximal ADMM, we need the following assumption.

\begin{assumption}\label{assumption:CQ2}
 There exists $(\hat u,\hat v)\in {\rm ri}({\rm dom}\,f\times{\rm dom}\,g) $ such that $\cF^* \hat u + \cG^* \hat v=c$.
\end{assumption}

\begin{theorem}\label{thempadm}  Let
 $\Sigma_f$ and $\Sigma_g $    be   the   self-adjoint and positive semidefinite operators defined by (\ref{monosub1}) and (\ref{monosub2}), respectively.
Suppose  that the solution set of problem \eqref{eq-dd-2} is nonempty and that  Assumption \ref{assumption:CQ2}  holds.   Assume that
 % both $\Sigma_f + T_f +\lambda A^*A $ and $\Sigma_g +T_g+ \lambda  B^*B $ are  positive definite.
$\cT_f$ and $\cT_g$ are chosen such that the
sequence $\{(u^k,v^k,x^k)\}$  generated by Algorithm SPADMM is well defined. % Let $(\bar u, \bar v)$ be any  optimal solution to problem \eqref{eq-dd-2} and $\bar{x}$ be any optimal solution to problem \eqref{eq-ADMM2-nonlinear-primal}, respectively.
%For $k=0,1, 2, \ldots, $ denote
% \[
% u_e^k: =u^k -\bar u, \quad v_e^k: =v^k -\bar v \quad {\rm and}\quad x_e^k: = x ^k -\bar {x}.
% \]
 Then, under the condition either
(a)  $\tau\in (0, (1+\sqrt{5}\,)/2)$ or (b)  $\tau \ge (1+\sqrt{5}\, )/2$ but $\sum_{k=0}^\infty (\|\cG^*(v^{k+1}-v^k)\|^2 + \tau^{-1} \|\cF^* u^{k+1} +\cG^* v^{k+1}-c\|^2) <\infty$,
 the following results hold:
\begin{enumerate}
%\item [{\rm (i)}] The sequence $\{\|x^{k+1}_e\|^2+\|v^{k+1}_e\|^2 _{(\sigma^{-1} \Sigma_g + \cT_g + \cG \cG^*)}+\|u^{k+1}_e\|^2 _{(\sigma^{-1}\Sigma_f +\cT_f +  \cF\cF^*)}\}$  is bounded.
%and
%    \[
%    \lim_{k\rightarrow \infty}\left( \|\lambda^{k+1}-\lambda^k \| + \|y^{k+1}-y^k \|_{\cT_g}+ \|x^{k+1}-x^k \|_{\cT_f}+ \|\cN^* (y^{k+1}- y^{k})\|\right)=0.
%    \]
 \item [{\rm (i)}]   If $(u^\infty, v^\infty, x^\infty)$ is an accumulation point of $\{(u^k, v^k, x^k)\}$, then $(u^\infty, v^\infty)$ solves problem \eqref{eq-dd-2} and $x^\infty$  solves \eqref{eq-ADMM2-nonlinear-primal}, respectively.
     % and it holds that
%     \[
%     \lim_{k\rightarrow \infty}\left(  \|x^{k+1}_e\|^2+\|v^{k+1}_e\|^2 _{(\sigma^{-1}\Sigma_g +\cT_g + \cG\cG^*)}+\|u^{k+1}_e\|^2 _{(\sigma^{-1}\Sigma_f +\cT_f +  \cF \cF^*)}\right)=0,
%     \]
%     where in the definition of  $(u_e^k, v^k_e, x_e^k)$, the point $(\bar u, \bar v, \bar x)$ is replaced by $(u^\infty, v^\infty, x^\infty)$.

   \item [{\rm (ii)}] If both $\sigma^{-1}\Sigma_f +\cT_f +  \cF\cF^*$ and $\sigma^{-1}\Sigma_g + \cT_g +\cG \cG^*$ are positive definite, then  the  sequence
$ \{(u^k, v^k, x^k)\}$, which  is automatically well defined, converges to a unique limit, say, $(u^\infty, v^\infty, x^\infty)$  with $(u^\infty, v^\infty)$ solving problem \eqref{eq-dd-2} and $x^\infty$ solving  \eqref{eq-ADMM2-nonlinear-primal}, respectively.
 \item [{\rm (iii)}]  When the $u$-part disappears, the corresponding results in parts (i)--(ii) hold under the condition either $\tau \in (0,2)$ or $\tau \ge2$ but $\sum_{k=0}^\infty \|\cG^* v^{k+1} -c\|^2<\infty$.
\end{enumerate}
\end{theorem}

\begin{remark}\label{rem:convergencePADMM2}
The conclusions of Theorem \ref{thempadm}  follow essentially from the results given in  \cite[Theorem B.1]{fazel2013hankel}. See \cite{styang2014} for more detailed discussions.
\end{remark}

Next, we shall pay particular attention to the case when  $g$ is a quadratic function:
\begin{eqnarray}
g(v) = \frac{1}{2}\inprod{v}{\Sigma_g v} -\inprod{b}{v}, \quad v \in \cV,
\label{q-gv}
\end{eqnarray}
where  $\Sigma_g$ a self-adjoint positive semidefinite linear operator defined on $\cV$ and $b\in \cV$ is a given vector. Problem (\ref{eq-dd-2}) now takes the form of
\begin{eqnarray}
  \begin{array}{rlll}
    \min & f(u)+ \frac{1}{2}\inprod{v}{\Sigma_g v} -\inprod{b}{v}    \\[5pt]
   \mbox{s.t.} & \cF^*u  + \cG^* v = c. \\[5pt]
   \end{array}
   \label{eq-dd-2q}
\end{eqnarray}
The dual of problem (\ref{eq-dd-2q}) is given by
\begin{eqnarray}
   \min \left\{\inprod{c}{x} + f^*(s)+ g^*(t)\mid  \cF x + s=0, \; \cG x + t=0 \right\}.
  \label{pg-eq-ADMM2-nonlinear-primal}
\end{eqnarray}
In order to solve subproblem \eqref{padm2v} in Algorithm \rm{SPADMM}, we need to solve a linear system with the linear operator given by $\sigma^{-1}\Sigma_g + \cG\cG^*$. Hence, an appropriate proximal term should be chosen such that \eqref{padm2v} can be solved efficiently. Here, we choose $\cT_g$ as follows.
Let $\cE_g : \cV \to \cV$ be a self-adjoint positive definite linear operator  such that it is a majorization of $\sigma^{-1}\Sigma_g + \cG\cG^*$, i.e.,
\[
 \cE_g \succeq \sigma^{-1}\Sigma_g + \cG\cG^*.
\]
We  choose $\cE_g $ such that its inverse
can be computed at a moderate cost. Define
 \begin{eqnarray}
  \cT_g:=\cE_g - \sigma^{-1}\Sigma_g - \cG\cG^* \succeq 0.
  \label{Tg}
\end{eqnarray}
Note that for numerical efficiency, we  need  the self-adjoint positive semidefinite linear operator $\cT_g$  to be as small as possible. In order to fully exploit the structure of the quadratic function $g$, we add, instead of a naive proximal term, a proximal term based on
the Schur complement as follows.
For a given $\cT_f\succeq 0$, we define the self-adjoint positive semidefinite linear operator
 \begin{eqnarray}
 \widehat{\cT}_f := \cT_f + \cF\cG^*\cE_g^{-1}\cG\cF^*. \label{TFhat}
 \end{eqnarray}

For later developments,
here we state a proposition which uses the Schur complement condition
for establishing the positive definiteness of a linear operator.
\begin{prop}\label{eqvi-psd-2}
It holds that
\begin{eqnarray*}\label{psd-f-g-2}
    \cW := \left(
       \begin{array}{c} \cF \\  \cG \end{array} \right)
 \left(
       \begin{array}{c} \cF \\  \cG \end{array} \right)^*
     + \sigma^{-1}\left(
                    \begin{array}{cc}
                      \Sigma_{f} &  \\
                       & \Sigma_{g} \\
                    \end{array}
                  \right)
      + \left(      \begin{array}{cc}
                      \widehat{\cT}_{f} &  \\
                      & \cT_{g} \\
                  \end{array}
       \right)\succ 0
\; \Leftrightarrow \;
    \cF\cF^* + \sigma^{-1}\Sigma_f + \cT_f \succ 0.
\end{eqnarray*}
\end{prop}
{\bf Proof.} We have that
\begin{eqnarray*}
\cW &=& \left(
               \begin{array}{cc}
          \cF\cF^*+\sigma^{-1}
     \Sigma_{f} + \widehat{\cT}_{f} & \cF\cG^* \\
            \cF^*\cG & \cG\cG^* + \sigma^{-1}\Sigma_{g} + \cT_{g} \\
       \end{array}
        \right).
\end{eqnarray*}
 Since $\cE_{g} = \cG\cG^* + \sigma^{-1}\Sigma_g + \cT_{g} \succ 0,$ by the Schur complement condition for ensuring  the  positive definiteness of linear operators, we have
$\cW \succ 0$ if and only if
 \begin{eqnarray*}
  & \cF\cF^* + \sigma^{-1}\Sigma_{f} + \widehat{\cT}_{f} - \cF\cG^*\cE_{g}^{-1}\cG\cF^* \succ 0.
\end{eqnarray*}
By \eqref{TFhat}, we  know that the conclusion of  this proposition holds.
 \qed \\

Now, we can propose our Schur complement based semi-proximal augmented Lagrangian method (SCB-SPALM) to solve (\ref{eq-dd-2q})
with a specially chosen proximal term involving $\widehat{\cT}_f$ and $\cT_g$.

\medskip
\centerline{\fbox{\parbox{\textwidth}{
{\bf Algorithm SCB-SPALM: \bf{A Schur complement based semi-proximal augmented Lagrangian method for solving (\ref{eq-dd-2q}).}} \\
Let $\sigma >0$ and $\tau\in(0,\infty)$ be given parameters.    Choose $(u^0,v^0,x^0)\in\mbox{dom}(f)\times\cV\times\cX.$ For $k=0,1,2,...$,
perform the $k$th iteration as follows:
\begin{description}
\item [Step 1.] Compute
  \begin{eqnarray}
     (u^{k+1},v^{k+1}) =\mbox{argmin}_{u,v} \;\cL_{\sigma}(u,v;x^k) + \frac{\sigma}{2}\norm{u - u^k}_{\widehat{\cT}_f}^2 + \frac{\sigma}{2}\norm{v-v^k}^2_{\cT_g}.
     \label{PALM-uv}
  \end{eqnarray}
\item [Step 2.] Compute
  \begin{eqnarray}
     x^{k+1} = x^k + \tau\sigma(\cF^*u^{k+1} + \cG^*v^{k+1} - c).
     \label{PALM-x}
  \end{eqnarray}
\end{description}
}}}
\medskip
Note that problem \eqref{PALM-uv} in Step 1 is well defined if the
the linear operator $\cW$ defined in Proposition \ref{eqvi-psd-2} is
positive definite, or equivalently, if $\cF\cF^*+\sig^{-1}\Sigma_f +\cT_f \succ 0$.
Also, note that in the context of the convex optimization problem \eqref{eq-dd-2q}, Assumption \ref{assumption:CQ2} is reduced to the following:
\begin{assumption}\label{assumption:CQ2-qg}
 There exists $(\hat u,\hat v)\in {\rm ri}({\rm dom}\,f)\times\cV $ such that $\cF^* \hat u + \cG^* \hat v=c$.
\end{assumption}
Now, we are ready to establish our convergence results for Algorithm SCB-SPALM for
 solving \eqref{eq-dd-2q}.
\begin{theorem}\label{thempalm}  Let
 $\Sigma_f$, $\Sigma_g$ and $\cT_g$ be three self-adjoint and positive semidefinite operators defined by (\ref{monosub1}), \eqref{q-gv} and \eqref{Tg}, respectively.
Suppose  that the solution set of problem \eqref{eq-dd-2q} is nonempty and that  Assumption \ref{assumption:CQ2-qg}  holds.   Assume that
 % both $\Sigma_f + T_f +\lambda A^*A $ and $\Sigma_g +T_g+ \lambda  B^*B $ are  positive definite.
$\cT_f$ is chosen such that the
sequence $\{(u^k,v^k,x^k)\}$  generated by Algorithm SCB-SPALM is well defined.
%Let $(\bar u, \bar v)$ be any  optimal solution to problem \eqref{eq-dd-2q} and $\bar{x}$ be any optimal solution to problem \eqref{pg-eq-ADMM2-nonlinear-primal}, respectively.
%For $k=0,1, 2, \ldots, $ denote
% \[
% u_e^k: =u^k -\bar u, \quad v_e^k: =v^k -\bar v \quad {\rm and}\quad x_e^k: = x ^k -\bar {x}.
% \]
Then, under the condition either
(a) $\tau \in (0,2)$ or (b) $\tau \ge2$ but $\sum_{k=0}^\infty \|\cF^* u^{k+1} + \cG^* v^{k+1} -c\|^2<\infty$,
 the following results hold:
\begin{enumerate}
%\item [{\rm (i)}] The sequence $\{\|x^{k+1}_e\|^2+\|v^{k+1}_e\|^2 _{(\sigma^{-1} \Sigma_g + \cT_g + \cG \cG^*)}+\|u^{k+1}_e\|^2 _{(\sigma^{-1}\Sigma_f +\cT_f +  \cF\cF^*)}\}$  is bounded.
%and
%    \[
%    \lim_{k\rightarrow \infty}\left( \|\lambda^{k+1}-\lambda^k \| + \|y^{k+1}-y^k \|_{\cT_g}+ \|x^{k+1}-x^k \|_{\cT_f}+ \|\cN^* (y^{k+1}- y^{k})\|\right)=0.
%    \]
 \item [{\rm (i)}]   If $(u^\infty, v^\infty, x^\infty)$ is an accumulation point of $\{(u^k, v^k, x^k)\}$, then $(u^\infty, v^\infty)$ solves problem \eqref{eq-dd-2q} and $x^\infty$  solves \eqref{pg-eq-ADMM2-nonlinear-primal}, respectively.
     %and it holds that
%     \[
%     \lim_{k\rightarrow \infty}\left(  \|x^{k+1}_e\|^2+\|v^{k+1}_e\|^2 _{(\sigma^{-1}\Sigma_g +\cT_g + \cG\cG^*)}+\|u^{k+1}_e\|^2 _{(\sigma^{-1}\Sigma_f +\cT_f +  \cF \cF^*)}\right)=0,
%     \]
%     where in the definition of  $(u_e^k, v^k_e, x_e^k)$, the point $(\bar u, \bar v, \bar x)$ is replaced by $(u^\infty, v^\infty, x^\infty)$.

 \item [{\rm (ii)}] If $\sigma^{-1}\Sigma_f +\cT_f +  \cF\cF^*$ is positive definite, then  the  sequence
$ \{(u^k, v^k, x^k)\}$, which is automatically well defined, converges to a unique limit, say, $(u^\infty, v^\infty, x^\infty)$  with $(u^\infty, v^\infty)$ solving problem \eqref{eq-dd-2q} and $x^\infty$ solving  \eqref{pg-eq-ADMM2-nonlinear-primal}, respectively.
 %\item [{\rm (iv)}]  When the $u$-part disappears, the corresponding results in parts (i)--(iii) hold under the condition either $\tau \in (0,2)$ or $\tau \ge2$ but $\sum_{k=0}^\infty \|\cG^* v^{k+1} -c\|^2<\infty$ .
\end{enumerate}
\end{theorem}
 {\bf Proof.} By combining Theorem \ref{thempadm} and Proposition \ref{eqvi-psd-2}, one can prove the results of this theorem directly.
 \qed

\bigskip
The relationship between Algorithm SCB-SPALM and Algorithm SPADMM for solving (\ref{eq-dd-2q}) will be revealed in the next proposition.

Let $\delta_g:\cU\times\cV\times\cX \rightarrow \cU$ be an auxiliary linear function associated with \eqref{eq-dd-2q} defined by
\begin{eqnarray}
  \delta_g(u,v,x)\;:=\;\cF\cG^*\cE_g^{-1}(b-\cG x - \Sigma_g v + \sigma\cG(c -\cF^*u-\cG^*v )).
  \label{dg}
\end{eqnarray}
Let  $\ubar\in \cU$, $\vbar\in \cV$, $\xbar\in\cX$ and $c \in \cX$ be given. Denote \[
\cbar \;: =\; c - \cF^*\ubar -\cG^*\vbar \quad {\rm and} \quad
\deltabar_g \;:=\;\delta_g(\ubar,\vbar,\xbar)
= \cF\cG^*\cE_g^{-1}(b-\cG \xbar - \Sigma_g \vbar + \sigma\cG \cbar).\]
Let  $(u^+,v^+)\in \cU\times \cV$ be defined  by
\begin{eqnarray}
  (u^+,v^+) = \mbox{argmin}_{u,v}\; \cL_\sigma(u,v;\xbar) +\frac{\sigma}{2}\norm{u-\ubar}_{\widehat{\cT}_f}^2
  + \frac{\sigma}{2}\norm{v-\vbar}_{\cT_g}^2.
  \label{prox2}
\end{eqnarray}

\begin{prop}\label{prop:prox2}
 Let $\alphabar : = \sigma^{-1}b + \cT_g \vbar + \cG(c- \sigma^{-1}\xbar)$. Define $v^\prime \in \cV$ by
   \begin{equation}\label{eq:vprime}
   \vprime = \textup{argmin}_{v}\; \cL_\sigma(\ubar,v;\xbar) + \frac{\sigma}{2}\norm{v - \vbar}_{\cT_g}^2
=  \cE_g^{-1}(\alphabar-\cG\cF^*\ubar).\end{equation}
 The optimal solution $(u^+,v^+)$ to problem (\ref{prox2}) is generated exactly by the following procedure
  \begin{eqnarray}
 \left\{ \begin{array}{lcl}
  u^+ &=& \textup{argmin}_{u}\; \cL_\sigma(u,\vbar;\xbar) + \inprod{\deltabar_g}{u} +
  \frac{\sigma}{2}\norm{u - \ubar}_{\cT_f}^2,
\\[5pt]
 v^+ &=& \textup{argmin}_{v}\; \cL_\sigma(u^+,v;\xbar) + \frac{\sigma}{2}\norm{v - \vbar}_{\cT_g}^2\; =\;  \cE_g^{-1}(\alphabar-\cG\cF^* u^+).
 \end{array} \right.
 \label{prox-m1}
 \end{eqnarray}
Furthermore,   $(u^+,v^+)$ can also be obtained  by the following equivalent procedure
  \begin{eqnarray}
 \left\{ \begin{array}{lcl}
  u^+ &=& \textup{argmin}_{u}\; \cL_\sigma(u,\vprime;\xbar) +
  \frac{\sigma}{2}\norm{u - \ubar}_{\cT_f}^2,
\\[5pt]
 v^+ &=& \textup{argmin}_{v}\; \cL_\sigma(u^+,v;\xbar) + \frac{\sigma}{2}\norm{v - \vbar}_{\cT_g}^2\; =\;  \cE_g^{-1}(\alphabar-\cG\cF^* u^+).
 \end{array} \right.
 \label{prox-m2}
\end{eqnarray}
\end{prop}
{\bf Proof.}
First we show that  the equivalence between   (\ref{prox2}) and  (\ref{prox-m1}).
Define
$$\widetilde{\cL}_\sig(u,v;\xbar) := \cL_{\sigma}(u,v;\xbar) + \frac{\sigma}{2}\norm{u-\ubar}_{\widehat{\cT}_f}^2
    +  \frac{\sigma}{2}\norm{v-\vbar}_{\cT_g}^2, \quad (u,v)\in \cU\times \cV.
$$
By simple algebraic manipulations, we have that
\begin{eqnarray}
 \widetilde{\cL}_\sig(u,v;\xbar)  = f(u) +
\frac{\sig}{2}\norm{u-\ubar}^2_{\widehat{\cT}_f} + \phi(u,v)  -\frac{1}{2\sig}\norm{\xbar}^2,
\label{eq-widetilde-cL-2}
\end{eqnarray}
where
\begin{eqnarray*}
 \phi(u,v) &=& g(v) + \frac{\sig}{2}\norm{\cF^* u + \cG^* v+\sig^{-1} \xbar-c}^2 +
\frac{\sig}{2}\norm{v-\vbar}_{\cT_g}^2
\nn \\[5pt]
&=& \frac{\sig}{2}\left(
 \inprod{v}{\cE_g v} + 2\inprod{v}{\cG\cF^* u -\alphabar}
+\norm{\cF^* u +\sig^{-1}\xbar-c}^2 + \norm{\vbar}_{\cT_g}^2
\right)
\end{eqnarray*}
with $\alphabar $ as defined in the proposition.
For any given $u\in \cU$, let
\begin{eqnarray*}
  v(u) \; :=\; \mbox{argmin}_{v\in \cV}\; \phi(u,v)
  \;=\; \cE_g^{-1}(\alphabar - \cG\cF^*u).
\end{eqnarray*}
Then by using the fact that
$\min_{v} \frac{1}{2}\inprod{v}{\cE_g v} + \inprod{q}{v} = -\frac{1}{2}\inprod{q}{\cE_g^{-1} q}$ for any $q\in \cV$, we have that
\begin{eqnarray*}
& &  \hspace{-0.7cm} \phi(u,v(u)) \;=\;
\frac{\sig}{2}\Big( -\inprod{\cG\cF^* u -\alphabar}{\cE_g^{-1}(\cG\cF^* u -\alphabar)}
+\norm{\cF^* u +\sig^{-1}\xbar-c}^2 + \norm{\vbar}_{\cT_g}^2
\Big)
\\[5pt]
& =&\frac{\sig}{2}\Big(\inprod{u}{(\cF\cF^*-\cF\cG^*\cE_g^{-1}\cG\cF^*) u}
 + 2\inprod{u}{\cF (\cG^*\cE_g^{-1}\alphabar+\sig^{-1}\xbar-c)}
\Big) + \kappa_0,
\end{eqnarray*}
where $\kappa_0 = \frac{\sig}{2}(\norm{\sig^{-1}\xbar-c}^2 + \norm{\vbar}_{\cT_g}^2 -
 \norm{\alphabar}_{\cE_g^{-1}}^2)$. Let
\begin{eqnarray*}
\kappa_1 : = \kappa_0 +\frac{\sig}{2}\norm{\cG\cF^*\ubar}_{\cE_g^{-1}}^2-\frac{1}{2\sig}\norm{\xbar}^2 =
 -\inprod{c}{\xbar} + \frac{\sig}{2}(\norm{c}^2 + \norm{\cG\cF^*\ubar}_{\cE_g^{-1}}^2 + \norm{\vbar}_{\cT_g}^2
-\norm{\alphabar}_{\cE_g^{-1}}^2).
\end{eqnarray*}
From \eqref{eq-widetilde-cL-2}, we have that  for any given $u\in \cU$,
\begin{eqnarray}
&&\hspace{-0.7cm}  \widetilde{\cL}_\sig(u,v(u);\xbar) \;=\; f(u) +
\frac{\sig}{2}\norm{u-\ubar}^2_{\cT_f}
+ \frac{\sig}{2}\norm{\cG\cF^*(u-\ubar)}^2_{\cE_g^{-1}}
+ \phi(u,v(u))  -\frac{1}{2\sig}\norm{\xbar}^2
 \nn \\[5pt]
  &=& f(u)+  \frac{\sig}{2}\norm{u-\ubar}^2_{\cT_f}
+\sig \inprod{u}{\cF (\cG^*\cE_g^{-1}\alphabar+\sig^{-1}\xbar-c) - \cF\cG^*\cE_g^{-1}\cG\cF^*\ubar}
 + \frac{\sigma}{2}\inprod{u}{\cF\cF^*u}  + \kappa_1
 \nn \\[5pt]
&=&  f(u)+  \frac{\sig}{2}\norm{u-\ubar}^2_{\cT_f}
+ \inprod{u}{\deltabar_g} + \inprod{u}{\cF(\xbar+\sig(\cG^*\vbar-c))}
 + \frac{\sigma}{2}\inprod{u}{\cF\cF^*u}  + \kappa_1
\nn \\[5pt]
&=&  \cL_\sig(u,\vbar;\xbar)   + \inprod{u}{\deltabar_g}
 +  \frac{\sig}{2}\norm{u-\ubar}^2_{\cT_f}
+\kappa_2,
  \label{eq-widetilde-cL-3}
\end{eqnarray}
where $\kappa_2 =  \kappa_1 -g(\vbar)
-\frac{\sig}{2}\norm{\cG^*\vbar-c}^2-\inprod{\xbar}{\cG^*\vbar-c}$.
Note that with some manipulations, we can show that the
constant term
$$
 \kappa_2 \;=\; \frac{\sig}{2}\norm{\cG\cF^*\ubar}_{\cE_g^{-1}}^2
-\frac{\sig}{2}\norm{\cE_g \vbar-\alphabar}_{\cE_g^{-1}}^2.
$$
Now, we have that
\begin{eqnarray*}
 \min_{u\in\cU,v\in\cV} \widetilde{\cL}_\sig(u,v;\xbar) =
\min_{u\in\cU} \Big(\min_{v\in\cV} \widetilde{\cL}_\sig(u,v;\xbar)\Big)
= \min_{u\in\cU}  \widetilde{\cL}_\sig(u,v(u) ;\xbar),
\end{eqnarray*}
where $ \widetilde{\cL}_\sig(u,v(u) ;\xbar)$ satisfies \eqref{eq-widetilde-cL-3}.
From here, the equivalence between   (\ref{prox2}) and  (\ref{prox-m1}) follows.

Next, we prove the equivalence between
  (\ref{prox-m1}) and (\ref{prox-m2}). Note that, the first minimization problem in (\ref{prox-m2}) can be equivalently recast
  as
\begin{eqnarray*}
  0 \in \partial{f}(u^+) + \cF \xbar + \sigma\cF(\cF^* u^+ + \cG^* \vprime  - c) + \sigma\cT_f(u^+ - \ubar),
\end{eqnarray*}
which, together with the definition of $\vprime$ given in \eqref{eq:vprime}, is equivalent to
\begin{eqnarray}
  0 \in \partial{f}(u^+) + \cF \xbar  + \sigma\cF(\cF^* u^+ -c
+ \cG^*\cE_g^{-1}(\alphabar - \cG\cF^*\ubar)  ) + \sigma\cT_f(u^+ - \ubar).
  \label{opu-prox-m2}
\end{eqnarray}
The condition (\ref{opu-prox-m2}) can be reformulated as
\begin{eqnarray*}
  0 \in \partial{f}(u^+) + \cF \xbar + \sigma\cF(\cF^* u^+ + \cG^*\vbar - c)
 + \sigma\cF\cG^*\cE_g^{-1}(\alphabar- \cG\cF^*\ubar-\cE_g \vbar)+\sigma\cT_f(u^+ - \ubar).
  \end{eqnarray*}
Thus, we have
  \begin{eqnarray}
   0 \in  \partial{f}(u^+) + \cF \xbar  + \sigma\cF(\cF^* u^+ + \cG^*\vbar - c)
+ \deltabar_g +\sigma\cT_f(u^+ - \ubar),
\end{eqnarray}
which can  equivalently be rewritten as
\begin{eqnarray*}
 u^+ =\mbox{argmin}_{u} \;\cL_{\sigma}(u,\vbar;\xbar) + \inprod{\deltabar_g}{u} + \frac{\sigma}{2}\norm{u - \ubar}_{\cT_f}^2.
\end{eqnarray*}
 The equivalence between
  (\ref{prox-m1}) and (\ref{prox-m2}) then follows. This completes the
proof of this proposition.
\qed

\vskip 15 true pt
%Let $\delta_g:\cU\times\cV\times\cX \rightarrow \cU$ be the function defined by
%\begin{eqnarray}
%  \delta_g(u,v,x):=\cF\cG^*\cE_g^{-1}(b-\cG x - \Sigma_g v + \sigma\cG(c -\cF^*u-\cG^*v )).
%  \label{dg}
%\end{eqnarray}

\begin{prop}\label{prop:equi-scb-palm}
Let $\delta_g^k := \delta_g(u^k,v^k,x^k)$ for $k=0,1,2,...$. We have that $u^{k+1}$ and $v^{k+1}$ obtained by  Algorithm SCB-SPALM
for solving \eqref{eq-dd-2q} can be generated exactly according to the following procedure:
  \begin{eqnarray} \label{padm2lu}
  \left\{\begin{array}{lcl}
     u^{k+1} &=&\textup{argmin}_{u} \;\cL_{\sigma}(u,v^k;x^k) + \inprod{\delta_g^k}{u} + \frac{\sigma}{2}\norm{u - u^k}_{\cT_f}^2,  \\[5pt]
    v^{k+1} &=& \textup{argmin}_{v} \;\cL_{\sigma}(u^{k+1},v;x^k) + \frac{\sigma}{2}\norm{v - v^k}_{\cT_g}^2, \\[5pt]
     x^{k+1} &=& x^k + \tau\sigma(\cF^*u^{k+1} + \cG^*v^{k+1} - c).
     \end{array}
     \right.
  \end{eqnarray}
\end{prop}
{\bf Proof.} The conclusion  follows directly from \eqref{prox-m1} in Proposition \ref{prop:prox2}.
\qed\bigskip

%Based on Proposition \ref{prop:prox2}, Algorithm PALM for solving (\ref{eq-dd-2q}) can be equivalently recast as follows.
%[{\bf make it a proposition.}]
%
%\centerline{\fbox{\parbox{\textwidth}{
%{\bf Algorithm \rm{II}: \bf{An equivalent form of PALM for sovling (\ref{eq-dd-2q}).}}\\
%Let $\sigma >0$ and $\tau\in(0,\infty)$ be given initial parameters. Let $\delta_g$ be defined by (\ref{dg}). Choose $(u^0,v^0,x^0)\in\mbox{dom}(f)\times\cV\times\cX.$ For $k=0,1,2,...$, generate $u^{k+1}$ and $v^{k+1}$ according to the following iteration:
%\begin{description}
%  \item [Step 1.] Compute
%  \begin{eqnarray}
%     u^{k+1} =\mbox{argmin}_{u} \;\cL_{\sigma}(u,v^k;x^k) + \inprod{\delta_g^k}{u} + \frac{\sigma}{2}\norm{u - u^k}_{\cT_f}^2,
%     \label{padm2lu}
%  \end{eqnarray}
%  where
%  \begin{eqnarray}
%  \delta_g^k := \delta_g(u^k,v^k,x^k).
%  \end{eqnarray}
%  \item [Step 2.] Compute
%  \begin{eqnarray}
%    v^{k+1} = \mbox{argmin}_{v} \;\cL_{\sigma}(u^{k+1},v;x^k) + \frac{\sigma}{2}\norm{v - v^k}_{\cT_g}^2.
%  \end{eqnarray}
%  \item [Step 3.] Compute
%  \begin{eqnarray}
%     x^{k+1} = x^k + \tau\sigma(\cF^*u^{k+1} + \cG^*v^{k+1} - c).
%  \end{eqnarray}
%\end{description}
%}}}
%\bigskip
\begin{remark}\label{rem:eq-palm}
{(i)} Note that comparing to (\ref{padm2u}) in Algorithm {SPADMM},
the first subproblem of (\ref{padm2lu})
has an extra linear term $\inprod{\delta_g^k}{\cdot}$.
It is this linear term that allows us to design a convergent SPADMM for solving multi-block convex optimization problems.
\\[5pt]
{(ii)} The linear term $\inprod{\delta_g^k}{\cdot}$
will vanish if $\Sigma_g = 0$, $\cE_g = \cG\cG^*\succ 0 $ and a proper starting point $(u^0,v^0,x^0)$ is chosen. Specifically, if we
choose $x^0\in \cX$ such that $\cG x^0 = b$ and $(u^0,v^0)\in {\rm dom}(f)\times\cV$ such that $v^0 = \cE_g^{-1}\cG(c - \cF^* u^0)$, then it holds that $\cG x^k = b$ and
$v^k = \cE_g^{-1}\cG(c-\cF^* u^k)$, which imply that
$\delta_g^k = 0$.
\\[5pt]
{(iii)} Observe that  when
 $\cT_f$ and $\cT_g$ are chosen to be $0$ in \eqref{padm2lu}, apart from the range of $\tau$,
 our Algorithm SCB-SPALM differs from the classical 2-block ADMM for solving problem
\eqref{eq-dd-2q} only in the linear term
$\inprod{\delta_g^k}{\cdot}$. This shows that the classical 2-block ADMM for solving problem \eqref{eq-dd-2q} has an unremovable deviation from the augmented Lagrangian method. This may explain why even when ADMM type methods suffer from slow local convergence, the latter can still enjoy fast local convergence.
\end{remark}

 In the following, we compare
  our Schur complement based proximal term
$\frac{\sig}{2}\norm{u-u^k}_{\widehat{\cT}_f}^2 + \frac{\sig}{2}\norm{v-v^k}^2_{\cT_g}$
used to derive the scheme \eqref{padm2lu} for solving \eqref{eq-dd-2q} with the
following proximal term which allows one to update $u$ and $v$ simultaneously:
\begin{eqnarray}
\frac{\sigma}{2}(\norm{(u,v) - (u^k,v^k)}^2_\cM +\norm{u-u^k}_{\cT_f}^2 + \norm{v-v^k}_{\cT_g}^2)
\quad \mbox{with} \quad
\cM = \left(
                 \begin{array}{cc}
                 \cD_1  & -\cF\cG^* \\
                -\cG\cF^*  & \cD_2  \\
              \end{array}
           \right)    \succeq 0,
\label{eq-cM}
\end{eqnarray}
where   $\cD_1:\cU\rightarrow\cU$ and $\cD_2:\cV\rightarrow\cV$ are two self-adjoint positive semidefinite linear operators satisfying
\[
\cD_1 \succeq \sqrt{ (\cF\cG^*) (\cF\cG^*)^*}\quad {\rm and}\quad  \cD_2 \succeq \sqrt{ (\cG\cF^*) (\cG\cF^*)^*}\, .
\]
A common naive choice will be $\cD_1 = \lambda_{\max}\cI_1$ and $\cD_2 = \lambda_{\max}\cI_2$ where $\lambda_{\max} = \norm{\cF\cG^*}_2,$
$\cI_1 :\cU \rightarrow\cU$ and $\cI_2 :\cV \rightarrow\cV$ are identity maps.
Simple calculations show that the resulting semi-proximal augmented Lagrangian method generates $(u^{k+1},v^{k+1},x^{k+1})$ as follows:
  \begin{eqnarray}
  \left\{\begin{array}{lcl}
     u^{k+1} &=&\textup{argmin}_{u} \;\cL_{\sigma}(u,v^k;x^k) + \frac{\sigma}{2}\norm{u - u^k}_{\cD_1 + \cT_f}^2, \label{s-nadmm} \\[5pt]
    v^{k+1} &=& \textup{argmin}_{v} \;\cL_{\sigma}(u^k,v;x^k) + \frac{\sigma}{2}\norm{v - v^k}_{\cD_2+ \cT_g}^2, \\[5pt]
     x^{k+1} &=& x^k + \tau\sigma(\cF^*u^{k+1} + \cG^*v^{k+1} - c).
     \end{array}
     \right.
  \end{eqnarray}
 To ensure that the subproblems in \eqref{s-nadmm} are well defined, we may require the following sufficient conditions to hold:
 \[\sigma^{-1}\Sigma_f + \cT_f + \cF\cF^* + \cD_1  \succ 0 \quad \mbox{and} \quad
 \sigma^{-1}\Sigma_g + \cT_g + \cG\cG^* + \cD_2  \succ 0.\]
Comparing the proximal terms used in \eqref{PALM-uv} and \eqref{eq-cM}, we can easily see
that the difference  is:
\begin{eqnarray*}
 \norm{u-u^k}^2_{\cF\cG^*\cE_g^{-1}\cG\cF^*} \quad \mbox{vs.} \quad
  \norm{(u,v)-(u^k,v^k)}_{\cM}^2.
\end{eqnarray*}
To simplify the comparison, we assume that
\[
\cD_1 = \sqrt{ (\cF\cG^*) (\cF\cG^*)^*}\quad {\rm and}\quad  \cD_2 = \sqrt{ (\cG\cF^*) (\cG\cF^*)^*}\, .
\]
% $\Sigma_g=0$ and $\cE_g = \cG\cG^*\succ 0$.
 By rescaling the equality constraint in \eqref{eq-dd-2q} if necessary, we may also assume that $\|\cF\|= 1$.
Now, we have that \[
\cF\cG^*\cE_g^{-1}\cG\cF^*
\preceq \cF \cF^*
\]
and
\[
  \norm{u-u^k}^2_{\cF\cG^*\cE_g^{-1}\cG\cF^*} \;\leq\; \norm{u-u^k}^2_{\cF\cF^*}\le  \norm{u-u^k}^2.
\]
In contrast, we have
\[
\begin{array}{lcl}
\norm{(u,v)-(u^k,v^k)}_{\cM}^2& \le &2\left (\norm{u-u^k}^2_{\cD_1} + \norm{v-v^k}^2_{\cD_2}\right)
\\[5mm]
&\le & 2\norm {\cF \cG^*} \left(\norm{u-u^k}^2 + \norm{v-v^k}^2\right)
\\[5mm]
&\le & 2\norm {\cG} \left(\norm{u-u^k}^2 + \norm{v-v^k}^2\right),
\end{array}
\]
which is larger than the former upper bound $\norm{u-u^k}^2$ if $\|\cG\|\ge 1/2$.
Thus we can conclude safely that the proximal term
$\norm{u-u^k}^2_{\cF\cG^*\cE_g^{-1}\cG\cF^*}$ can be potentially  much smaller
than $ \norm{(u,v)-(u^k,v^k)}^2_{\cM}$ unless $\|\cG\|$ is very small.

 The above mentioned  upper bounds   difference is of course due to the fact that the   SCB semi-proximal augmented Lagrangian method  takes advantage of the fact that $g$ is assumed to be  a convex quadratic function.  However, the key difference lies in the fact
that \eqref{s-nadmm} is a splitting version of the semi-proximal augmented Lagrangian method with a Jacobi type decomposition, whereas Algorithm SCB-SPALM is a splitting version of semi-proximal augmented Lagrangian method with a Gauss-Seidel type decomposition.
It is this fact that provides us with the key idea to design Schur complement based
proximal terms for multi-block convex optimization problems in the next section.

\section{A Schur complement based semi-proximal ADMM}

In this section, we focus on the problem
%\eqref{eq-dd-M1}
\begin{eqnarray}
  \begin{array}{rlll}
    \min & f(u)+\sum_{i=1}^p \theta_i(y_i) + g(v) + \sum_{j=1}^q \varphi_j(z_j)    \\[5pt]
   \mbox{s.t.} & \cF^*u + \sum_{i=1}^p \cA_i^*y_i  + \cG^* v + \sum_{j=1}^q \cB_j^* z_j  = c  \\[5pt]
   \end{array}
   \label{eq-dd-M}
\end{eqnarray}
with  all $\theta_i$ and $\varphi_j$  being assumed to be convex quadratic functions:
\begin{eqnarray*}
  \theta_i(y_i) = \frac{1}{2}\inprod{y_i}{\cP_i y_i} -\inprod{b_i}{y_i} ,  \;i=1,
  \ldots,p, \qquad \varphi_j(z_j) = \frac{1}{2}\inprod{z_j}{\cQ_j z_j}-\inprod{d_j}{z_j}, \;j=1,\ldots,q,
\end{eqnarray*}
where $\cP_i$ and $\cQ_j$ are given self-adjoint positive semidefinite linear operators.
The dual of (\ref{eq-dd-M}) is given by
\begin{eqnarray}
   \max \big\{ -\inprod{c}{x}-f^*(-\cF x)-\sum_{i=1}^p \theta^*_i(-\cA_i x) - g^*( -\cG x )-\sum_{j=1}^q \varphi^*_j(-\cB_j x)  \big\},
   \label{eq-pp-Mori}
\end{eqnarray}
which can equivalently be written as
\begin{eqnarray}
  \begin{array}{rlll}
    \min & \inprod{c}{x} + f^*(s)+\sum_{i=1}^p \theta^*_i(r_i) + g^*(t) + \sum_{j=1}^q \varphi^*_j(w_j)    \\[5pt]
   \mbox{s.t.} & \cF x + s =0, \quad \cA_i x + r_i = 0, \quad i=1,\ldots,p, \\[5pt]
   &\cG x + t = 0, \quad \cB_j x + w_j = 0, \quad j=1,\ldots,q.
   \end{array}
   \label{eq-pp-M}
\end{eqnarray}
%where  $\theta_i^*$ and $\varphi_j^*$ are given as follows if $\cP_i$ and $\cQ_j$ are invertible:
%\begin{eqnarray*}
%  \theta_i^*(r_i) = \frac{1}{2}\inprod{r_i+b_i}{\cP^{-1}_i(r_i+b_i)},  \qquad
%   \varphi^*_j(w_j) = \frac{1}{2}\inprod{w_j+d_j}{\cQ_j^{-1} (w_j+d_j)}
%\end{eqnarray*}

For $i=1,\ldots,p$, let $\cE_{\theta_i} $ be a self-adjoint positive definite linear operator  on $\cY_i$ such that it is a majorization of $\sigma^{-1}\cP_i + \cA_i\cA_i^*$, i.e.,
\[
 \cE_{\theta_i} \succeq \sigma^{-1}\cP_i + \cA_i\cA_i^*.
\]
We  choose $\cE_{\theta_i} $ in a way that its inverse
can be computed at a moderate cost. Define
 \begin{eqnarray}
  \cT_{\theta_i} := \cE_{\theta_i} - \sigma^{-1}\cP_i - \cA_i\cA_i^* \succeq 0,\quad i=1,\ldots,p.
  \label{Ttai}
\end{eqnarray}
Note that for numerical efficiency, we need the
 self-adjoint positive semidefinite linear operator $\cT_{\theta_i}$  to be as small as possible
for each $i$.
Similarly, for $j=1,\ldots,q$, let $\cE_{\varphi_j}$ be  a self-adjoint positive definite linear operator  on $\cZ_j$ that majorizes $\sigma^{-1}\cQ_j + \cB_j\cB_j^*$ in a way  that
$\cE_{\varphi_j}^{-1}$ can be computed relatively easily. Denote
 \begin{eqnarray}
  \cT_{\varphi_j}:=\cE_{\varphi_j}- \sigma^{-1}\cQ_j - \cB_j\cB_j^* \succeq 0,\quad j=1,\ldots,q.
  \label{Tvarphi}
\end{eqnarray}
Again, we need the self-adjoint positive semidefinite linear operator $\cT_{\varphi_j}$
to be as small as possible for each $j$.

For notational convenience, we define
\[
y_{\leq i} := (y_1,y_2,\ldots,y_i), \quad y_{\geq i} := (y_i,y_{i+1},\ldots,y_p),\; i=0,\ldots,p+1
\] with the convention that
$y_0 = y_{p+1} = y_{\leq0} = y_{\geq p+1}=\emptyset.$ For $i=1,\ldots,p,$ define the linear operator $\cA_{\leq i}: \cX \rightarrow \cY$ by
\begin{eqnarray*}
\left( \begin{array}{c}
                        \cA_1x \\
                        \cA_2x \\
                        \vdots \\
                        \cA_ix \\
                      \end{array}
                    \right)
\;\equiv\; \cA_{\leq i} x
\;:=\;
 \cA_1x\times\cA_2x\ldots\times\cA_ix \quad \forall x\in \cX.
 \end{eqnarray*}
 In a similar manner, we can define $z_{\leq j}, z_{\geq j}$ for $j=0,\ldots,q+1$ and
 define the linear operator $\cB_{\leq j}$ for $j=1,\ldots, q.$ Note that by definition, we have
 $y= y_{\leq p}$, $z= z_{\leq q}$, $\cA = \cA_{\leq p}$ and $\cB=\cB_{\leq q}$.
%For $1\leq i\leq p$, let $y_{\leq i} := (y_1,y_2,\ldots,y_i) \in \cY_1\times\cY_2\ldots\times\cY_i$, $y_{\geq i} := (y_i,y_{i+1},\ldots,y_p)$ and $y_{\leq 0}=y_{\geq p+1} = \emptyset,$ let linear operator $\cA_{\leq i}: \cX \rightarrow \cY_1\times\cY_2\ldots\times\cY_i$ be
%\[\cA_{\leq i}x :=\left(
%                      \begin{array}{c}
%                        \cA_1x \\
%                        \cA_2x \\
%                        \vdots \\
%                        \cA_ix \\
%                      \end{array}
%                    \right) :=
% (\cA_1x,\cA_2x,\ldots,\cA_ix);\]
%for $1\leq j\leq q$, let $z_{\leq j} := (z_1,z_2,\ldots,z_j) \in \cZ_1\times\cZ_2\ldots\times\cZ_j$ , $z_{\geq j} := (z_j,z_{j+1},\ldots,z_q)$ and $z_{\leq 0}=z_{\geq q+1} = \emptyset,$ let linear operator $\cB_{\leq j}: \cX \rightarrow \cZ_1\times\cZ_2\ldots\times\cZ_j$ be
%\[\cB_{\leq i}x :=\left(
%                      \begin{array}{c}
%                        \cB_1x \\
%                        \cB_2x \\
%                        \vdots \\
%                        \cB_jx \\
%                      \end{array}
%                    \right) :=
% (\cB_1x,\cB_2x,\ldots,\cB_ix).\]

Define the affine function $\Gamma: \cU\times\cY\times\cV\times\cZ \rightarrow \cX$ by
\begin{equation}\label{eq:GammaMapping}
\Gamma(u,y,v, z) := \cF^*u+\cA^*y+\cG^*v+\cB^*z-c\quad \forall \, (u,y, v,z)\in \cU\times \cY\times \cV\times \cZ.
\end{equation}
Let $\sigma > 0$ be given. The augmented Lagrangian function associated with
\eqref{eq-dd-M}
 is given as follows:
\begin{eqnarray}
\cL_{\sigma}(u,y,v,z;x) = f(u) + \theta(y) + g(v) +\varphi(z)
+\inprod{x}{\Gamma(u,y,v, z)}
+ \frac{\sigma}{2}\norm{\Gamma(u,y,v, z)}^2 \label{auglagfunction}
\end{eqnarray}
where $\theta(y) = \sum_{i=1}^p \theta_i(y_i)$ and
$\varphi(z) = \sum_{j=1}^q \varphi_j(z_j)$.

 Now we are ready to present our SCB-SPADMM (Schur complement based semi-proximal alternating direction method of multipliers) algorithm for solving (\ref{eq-dd-M}).

\bigskip
\centerline{\fbox{\parbox{\textwidth}{
{\bf Algorithm SCB-SPADMM: A Schur complement based SPADMM for solving (\ref{eq-dd-M}).}\\
Let $\sigma >0$ and $\tau\in(0,\infty)$ be given  parameters. Let $\cT_f$ and $\cT_g$ be
given self-adjoint positive semidefinite operators defined on $\cU$ and $\cV$ respectively. Choose $(u^0, y^0, v^0, z^0, x^0)\in\mbox{dom}(f)\times\cY\times\mbox{dom}(g)
\times\cZ\times\cX.$
For $k=0,1,2,...$, generate $(u^{k+1}, y^{k+1},v^{k+1},z^{k+1})$ and $x^{k+1}$
according to the following iteration.
\begin{description}
  \item [Step 1.] Compute for $i=p,  \ldots, 1,$
  \begin{eqnarray}
     \overline{y}_i^k \;=\; \mbox{argmin}_{y_i}
\;\cL_{\sigma}(u^k,(y_{\leq i-1}^k,y_i,\overline{y}_{\geq i+1}^{k}),v^k,z^k;x^k) + \frac{\sigma}{2}\norm{y_i - y_i^k}_{\cT_{\theta_i}}^2,
     \label{padmMyi}
  \end{eqnarray}
where  $\cT_{\theta_i}$ is defined as in (\ref{Ttai}). Then
  compute
   \begin{eqnarray}
    u^{k+1} = \mbox{argmin}_{u} \;\cL_{\sigma}(u,\overline{y}^k,v^k,z^k;x^k) + \frac{\sigma}{2}\norm{u - u^k}_{\cT_f}^2.
    \label{padmMu}
  \end{eqnarray}
  \item [Step 2.] Compute for $i=1,\ldots,p,$
  \begin{eqnarray}
     y_i^{k+1} \;=\; \mbox{argmin}_{y_i} \;\cL_{\sigma}(u^{k+1},(y_{\leq i-1}^{k+1},y_i,\overline{y}_{\geq i+1}^k),v^k,z^k;x^k) + \frac{\sigma}{2}\norm{y_i - y_i^k}_{\cT_{\theta_i}}^2.
  \end{eqnarray}
  \item [Step 3.] Compute for $j=q,\ldots,1$,
  \begin{eqnarray}
      \overline{z}_j^k \;=\; \mbox{argmin}_{z_j} \;\cL_{\sigma}(u^{k+1},y^{k+1},v^k,(z_{\leq j-1}^k,z_j,\overline{z}_{\geq j+1}^{k});x^k)
+ \frac{\sigma}{2}\norm{z_j - z_j^k}_{\cT_{\varphi_j}}^2,
  \end{eqnarray}
 where  $\cT_{\varphi_j}$ is defined as in (\ref{Tvarphi}). Then compute
  \begin{eqnarray}
    v^{k+1} = \mbox{argmin}_{v} \;\cL_{\sigma}(u^{k+1},y^{k+1},v,\overline{z}^{k};x^k) + \frac{\sigma}{2}\norm{v - v^k}_{\cT_g}^2.
    \label{padmMz1}
  \end{eqnarray}
  \item [Step 4.] Compute for $j=1,\ldots,q$,
   \begin{eqnarray}
      z_j^{k+1} \;=\; \mbox{argmin}_{z_j} \;\cL_{\sigma}(u^{k+1},y^{k+1},v^{k+1},(z_{\leq j-1}^{k+1},z_j,\overline{z}_{\geq j+1}^k);x^k) + \frac{\sigma}{2}\norm{z_j - z_j^k}_{\cT_{\varphi_j}}^2.
  \end{eqnarray}
  \item [Step 5.] Compute
  \begin{eqnarray}
    x^{k+1} = x^k + \tau \sigma(\cF^*u^{k+1} +  \cA^*y^{k+1} + \cG^*v^{k+1} +  \cB^*z^{k+1} -c).
    \label{padmMx}
  \end{eqnarray}
\end{description}
}}}
\bigskip
In order to prove the convergence of Algorithm SCB-SPADMM for solving (\ref{eq-dd-M}), we need first to study the relationship between SCB-SPADMM and the generic 2-block semi-proximal ADMM for solving a two-block convex optimization problem discussed in the previous section.

%We state a proposition here for later developments.
 Define for $l=1,\ldots, p,$
\begin{eqnarray}
&
f_1(u) := f(u), \quad  f_{l+1}(u,y_{\leq l}) := f(u) + \sum_{i=1}^l \theta_i(y_i)
\quad
 \forall\, (u,y_{\leq l}) \in
  \cU\times\cY_{\leq l},&
  \label{flp1}
\end{eqnarray}
where $\cY_{\leq l} = \cY_1\times\cY_2\times \ldots\times\cY_l$.
Similarly,  for $l=1,\ldots,q$, define
$\cZ_{\leq l} = \cZ_1\times\cZ_2\times \ldots\times\cZ_l$, and
\begin{eqnarray}
 & g_1(v) := g(v), \quad
g_{l+1}(v,z_{\leq l})\; :=\;  g(v) + \sum_{j=1}^l\varphi_{j}(z_j) \quad\forall\,
   (v,z_{\leq l})\in \cV\times\cZ_{\leq l}.&
  \label{glp1}
\end{eqnarray}
Denote $\cA_0^* \equiv \cF_1^* \equiv \cF^*$ and $\cB_0^* \equiv \cG_1^* \equiv \cG^*$. Let
\begin{eqnarray*}
 \cF_{i+1}^* = \Big(\cF^*, \cA_1^*,\ldots,\cA_i^*\Big), \quad i=1,\ldots,p, \qquad
\cG_{j+1}^* = \Big(\cG^*, \cB_1^*,\ldots,\cB_j^*\Big), \quad j=1,\ldots,q.
\end{eqnarray*}
Define the following  self-adjoint linear operators: $\widehat{\cT}_{f_1} := \cT_f + \cF_1\cA_1^*\cE_{\theta_1}^{-1}\cA_1\cF_1^*$,
\begin{eqnarray}
\widehat{\cT}_{f_i} := \left(
                  \begin{array}{cc}
                    \widehat{\cT}_{f_{i-1}} &  \\
                     & \cT_{\theta_{i-1}} \\
                  \end{array}
                \right) + \cF_{i}\cA_{i}^*\cE_{\theta_i}^{-1}\cA_{i}\cF_{i}^*,
                \qquad i=2,\ldots,p
  \label{Thi}
\end{eqnarray}
and $\widehat{\cT}_{g_1} := \cT_g + \cG_1\cB_1^*\cE_{\varphi_1}^{-1}\cB_1\cG_1^*$,
\begin{eqnarray}
\widehat{\cT}_{g_j} := \left(
                  \begin{array}{cc}
                    \widehat{\cT}_{g_{j-1}} &  \\
                     & \cT_{\varphi_{j-1}} \\
                  \end{array}
                \right) + \cG_{j}\cB_j^*\cE_{\varphi_j}^{-1}\cB_j\cG_{j}^*,
                \qquad j=2,\ldots,q.
\label{Tphi}
\end{eqnarray}
%Similarly, %define
%%\begin{eqnarray}
%%  \cT_{\widehat{f}_1} &:=& \cT_g + \cG\cB_1^*\cE_{\varphi_1}^{-1}\cB_1\cG^*
%%\end{eqnarray}
%%and
%for $j=1,2,\ldots,q$, define
%\begin{eqnarray}
%\cT_{g_j} &:=& \left(
%                  \begin{array}{cc}
%                    \cT_{g_{j-1}} &  \\
%                     & \cT_{\varphi_{j-1}} \\
%                  \end{array}
%                \right) + \cG_{j-1}\cB_j^*\cE_{\varphi_j}^{-1}\cB_j\cG_{j-1}^*.
%\label{Tphi}
%\end{eqnarray}
Let $(\vbar,\zbar,\xbar,c)\in \cV \times \cZ \times \cX \times \cX$ be given.
Denote
\[\cbar \; :=\;  c-\cG^*\vbar-\cB^*\zbar \quad {\rm and} \quad \gammabar \;:=\;
- \Gamma(\ubar,\ybar,\vbar,\zbar).
\]
Define
\begin{eqnarray}
\betabar _{p,j}  \; :=\;
 \cA_{j-1}\cA_p^*\cE_{\theta_p}^{-1}(b_p -\cA_p \xbar -\cP_p \ybar_p +
\sigma\cA_p\gammabar), \quad j=1, \ldots, p
\label{betapj}
\end{eqnarray}
and for  $i=p-1,  \ldots, 1$,
\begin{eqnarray}
\betabar _{i,j}\; : = \; \cA_{j-1}\cA_i^*\cE_{\theta_i}^{-1}\left (b_i -\sum _{k=i+1}^p\betabar_{k,i+1}  -\cA_i \xbar -\cP_i \ybar_i +
\sigma\cA_i\gammabar\right), \quad j=1, \ldots, i.
\label{betaij}
\end{eqnarray}
Let
\begin{eqnarray}\label{eq:linear term}
 \deltabar_\theta \; : = \; \sum_{i=1}^p\betabar_{i,1}.
\end{eqnarray}
We will show later in Proposition \ref{prop:proxaugM}
that $\deltabar_\theta$ is
the auxiliary linear term  associated with  problem (\ref{eq-dd-M}).
Recall that
\begin{eqnarray}
 \cL_{\sigma}(u,y,\vbar,\zbar;\xbar) = f(u) +  \theta(y) + g(\vbar) +\varphi(\zbar) +\inprod{\xbar}{\Gamma(u,y,\vbar,\zbar)}
+ \frac{\sigma}{2}\norm{\Gamma(u,y,\vbar,\zbar)}^2. \nn
\end{eqnarray}
For $i=p, \ldots,1,$  let  $y_i^\prime\in \cY_i$ be defined  by
\begin{eqnarray}
  \label{y-prime-i}
   y^\prime_i &:=& \textup{argmin}_{y_i} \;\cL_{\sigma}(\ubar,(\ybar_{\leq i-1},y_i,y^\prime_{\geq i+1}),\vbar,\zbar;\xbar) + \frac{\sigma}{2}\norm{y_i - \ybar_i}_{\cT_{\theta_i}}^2 \nn\\[5pt]
   &=& \cE_{\theta_i}^{-1}\big(\sigma^{-1}b_i - \sigma^{-1}\cA_i\xbar +
   \cT_{\theta_i}\ybar_i + \cA_i\cA_i^*\ybar_i - \cA_i\Gamma(\ubar,(\ybar_{\leq i-1},\ybar_i,y^\prime_{\geq i+1}),\vbar,\zbar)\big)
\end{eqnarray}
with the convention $y^\prime_{p+1} =\emptyset$.
Define $(u^+, y^+) \in \cU\times \cY$  by
\begin{eqnarray}
   (u^+,y^+) := \mbox{argmin}_{u,y}\; \cL_{\sigma}(u,y,\vbar,\zbar;\xbar)
  +\frac{\sigma}{2}\norm{(u,y_{\leq p-1}) - (\ubar,\ybar_{\leq p-1})}_{\widehat{\cT}_{f_p}}^2
  + \frac{\sigma}{2}\norm{y_p-\ybar_p}_{\cT_{\theta_p}}^2. \label{prox-cM1}
\end{eqnarray}
The following proposition about  two other equivalent procedures for computing $(u^+, y^+)$ is the key ingredient for our algorithmic developments.
The idea of proving this proposition is very simple:  use   Proposition \ref{prop:prox2} repeatedly though the proof itself is rather lengthy due to the multi-layered nature of the   problems involved.
For \eqref{prox-cM1}, we first express $y_p$ as a function
of $(u, y_{\leq p-1})$ to obtain a problem involving only
$(u,y_{\leq p-1})$, and from the resulting problem, express $y_{p-1}$ as  a
function of $(u,y_{\leq p-2})$ to get another problem involving only
$(u,y_{\leq p-2})$.  We continue this way until we get a
problem involving only $(u,y_1)$.

\begin{prop}\label{prop:proxaugM}
The optimal solution $(u^+,y^+)$ defined by \eqref{prox-cM1} can be obtained exactly by
  \begin{eqnarray}
 \left\{ \begin{array}{lcl}
  u^+ &=& \textup{argmin}_{u}\; \cL_\sigma(u,\ybar,\vbar,\zbar;\xbar) +  \inprod{ \deltabar_{\theta}}{u} +
  \frac{\sigma}{2}\norm{u - \ubar}_{\cT_f}^2,
\\[5pt]
  y_i^+ &=& \textup{argmin}_{y_i} \;\cL_{\sigma}(u^+,(y_{\leq i-1}^+,y_i, y^\prime_{\geq i+1}),\vbar,\zbar;\xbar) + \frac{\sigma}{2}\norm{y_i - \ybar_i}_{\cT_{\theta_i}}^2, \quad i=1,\ldots,p,
 \end{array} \right.
 \label{padmmMl}
 \end{eqnarray}
where the auxiliary linear term  $\deltabar_\theta$ is defined by (\ref{eq:linear term}). Furthermore, $(u^+,y^+)$ can also be generated by the following equivalent procedure
\begin{eqnarray}
 \left\{ \begin{array}{lrl}
 %y^\prime_i &=& \textup{argmin}_{y_i} \;\cL_{\sigma}(\ubar,\ybar_{\leq i-1},y_i,y^\prime_{\geq i+1},\vbar,\zbar;\xbar) + \frac{\sigma}{2}\norm{y_i - \ybar_i}_{\cT_{\theta_i}}^2, \quad i=p,p-1,\ldots,1,\\[5pt]
  u^+ &=& \textup{argmin}_{u}\; \cL_\sigma(u,y^\prime,\vbar,\zbar;\xbar) +
  \frac{\sigma}{2}\norm{u - \ubar}_{\cT_f}^2,
\\[5pt]
 y_i^+ &=& \textup{argmin}_{y_i} \;\cL_{\sigma}(u^+,(y_{\leq i-1}^+,y_i,
y^\prime_{\geq i+1}),\vbar,\zbar;\xbar) + \frac{\sigma}{2}\norm{y_i - \ybar_i}_{\cT_{\theta_i}}^2, \quad i=1,\ldots,p.
 \end{array} \right.
 \label{padmmMy}
 \end{eqnarray}
\end{prop}

{\bf Proof.}  We will separate our proof into two parts and for each part  we prove our conclusions  by induction.

{\bf Part one.}  In this part we show that  $(u^+,y^+)$ defined by \eqref{prox-cM1} can be obtained exactly by (\ref{padmmMl}). For the case $p=1$, this follows directly from
  Proposition \ref{prop:prox2}.

 Assume that the equivalence   between \eqref{prox-cM1} and \eqref{padmmMl}
holds for all $p \leq l$.  We need to show   that for $p=l+1$,
 this equivalence also  holds. For this purpose,
we consider the following optimization problem with respect to $(u,y_{\leq l})$ and $y_{l+1}$:
\begin{eqnarray}
    \begin{array}{rlll}
    \min & f_{l+1}(u,y_{\leq l}) + \theta_{l+1}(y_{l+1}) + g(\vbar) + \varphi(\zbar) \\[5pt]
   \mbox{s.t.} & \cF_{l+1}^*(u,y_{\leq l})  + \cA_{l+1}^*y_{l+1} = \cbar.
   \label{prop-eq-l1}
   \end{array}
\end{eqnarray}
The augmented Lagrangian function associated with problem \eqref{prop-eq-l1} is given by
\begin{eqnarray}
{\cL}^{l+1}_{\sigma}((u,y_{\leq l}),y_{l+1};\vbar,\zbar,x)
&=& f_{l+1}(u,y_{\leq l}) + \theta_{l+1}(y_{l+1}) + g(\vbar) +\varphi(\zbar) \nn\\[5pt]
&&+\inprod{x}{\Gamma(u,y,\vbar,\zbar)} + \frac{\sigma}{2}\norm{\Gamma(u,y,\vbar,\zbar)}^2.
\end{eqnarray}
We denote the vector $\delta_{\theta_{l+1}}$ as the
auxiliary linear term associated with problem \eqref{prop-eq-l1} by
\begin{eqnarray}
  \delta_{\theta_{l+1}} \; :=\;  \cF_{l+1}\cA^*_{l+1}\cE_{\theta_{l+1}}^{-1}
  (b_{l+1} - \cA_{l+1}\xbar - \cP_{l+1}\ybar_{l+1} + \sigma\cA_{l+1}\gammabar).
\end{eqnarray}
Note that by the definition of $\cF_{l+1}$ and $p=l+1$, we have
 \begin{eqnarray*}
& \inprod{\delta_{\theta_{p}}}{(u,y_{\leq l})}
\;=\; \inprod{\betabar_{p,1}}{u} + \sum_{j=1}^l \inprod{\betabar_{p,j+1}}{y_j} &
 \end{eqnarray*}
with $\betabar_{p,j}$, $j=1,\ldots,l+1$, defined as in \eqref{betapj}.

By noting that $ {\cL}^{l+1}_{\sigma}((u,y_{\leq l}),y_{l+1};\vbar,\zbar,\xbar)= {\cL}_{\sigma}(u,y_{\leq l},y_{l+1}, \vbar,\zbar;\xbar)$,
we can rewrite problem \eqref{prox-cM1} for $p=l+1$ equivalently as
\begin{eqnarray}
  ((u^{+},y^{+}_{\leq l}),y^{+}_{l+1})
\;=\;\mbox{argmin}
\left\{ \begin{array}{l} {\cL}^{l+1}_{\sigma}((u,y_{\leq l}),y_{l+1};\vbar,\zbar,\xbar) + \frac{\sigma}{2}\norm{(u,y_{\leq l}) - (\ubar,\ybar_{\leq l})}_{\widehat{\cT}_{f_{l+1}}}^2   \\[5pt]
 + \frac{\sigma}{2}\norm{y_{l+1}-\ybar_{l+1}}_{\cT_{\theta_{l+1}}}^2
\end{array}\right\}.
  \label{prop-eq-l1p}
\end{eqnarray}
Then, from  Proposition \ref{prop:prox2}, we know that problem \eqref{prop-eq-l1p} is equivalent to
\begin{eqnarray}
   (u^+,y^+_{\leq l})  &=&
\textup{argmin}_{(u,y_{\leq l})}\;
\left\{ \begin{array}{l}
 {\cL}^{l+1}_\sigma((u,y_{\leq l}),\ybar_{l+1};\vbar,\zbar,\xbar) + \inprod{\delta_{\theta_{l+1}}}{(u,y_{\leq l})}
\\[5pt]
   +\frac{\sigma}{2}\norm{(u,y_{\leq l-1}) - (\ubar,\ybar_{\leq l-1})}_{\widehat{\cT}_{f_l}}^2
  + \frac{\sigma}{2}\norm{y_{l} - \ybar_{l}}_{\cT_{\theta_l}}^2
\end{array}\right\},
\label{padmmM1uy}
\\[5pt]
y_{l+1}^{+} &=& \textup{argmin}_{y_{l+1}} \;{\cL}^{l+1}_{\sigma}((u^+,y_{\leq l}^+),y_{l+1};\vbar,\zbar,\xbar) + \frac{\sigma}{2}\norm{y_{l+1} - \ybar_{l+1}}_{\cT_{\theta_{l+1}}}^2.
 \label{padmmM1y}
\end{eqnarray}
By observing that ${\cL}^{l+1}_{\sigma}((u^+,y_{\leq l}^+),y_{l+1};\vbar,\zbar,\xbar) = \cL_{\sigma}(u^+,y_{\leq l}^+,y_{l+1},\vbar,\zbar;\xbar)$,
we know that  problem  \eqref{padmmM1y}  can  equivalently be rewritten   as
\begin{eqnarray}
  y_{l+1}^{+} = \textup{argmin}_{y_{l+1}} \;\cL_{\sigma}(u^+,y_{\leq l}^+,y_{l+1},\vbar,\zbar;\xbar) + \frac{\sigma}{2}\norm{y_{l+1} - \ybar_{l+1}}_{\cT_{\theta_{l+1}}}^2.
  \label{prop-eq-yp}
\end{eqnarray}
In order to apply our induction assumption to problem (\ref{padmmM1uy}),  we  need to construct a corresponding optimization problem.
Define  for $i=1,\ldots,l$,
\begin{eqnarray*}
&\widetilde{b}_i := b_i - \betabar _{p,i+1} \quad {\rm and}\quad  \widetilde{\theta_i}(y_i) := \theta_i(y_i)+ \inprod{\betabar _{p,i+1}}{y_i}
=\frac{1}{2}\inprod{y_i}{\cP_i y_i} - \inprod{\widetilde{b}_i}{y_i}
\quad \forall\, y_i \in \cY_i, &
\\[5pt]
&
\widetilde{f}_1(u) :=  f(u) + \inprod{\betabar _{p,1}}{u}, \quad
\widetilde{f}_{i+1}(u,y_{\leq i}) := \widetilde{f}_1(u)+
\sum_{j=1}^{i} \widetilde{\theta}_j(y_j)
 \quad \forall\, (u,y_{\leq i}) \in \cU\times\cY_{\leq i}. &
\end{eqnarray*}
We shall now consider the following  optimization problem with respect to $(u,y_{\le l})$:
\begin{eqnarray}
    \begin{array}{rlll}
    \min & \widetilde{f}_1(u) + \sum_{i=1}^l\widetilde{\theta}_{i}(y_{i}) +\theta_{l+1}(\ybar_{l+1}) + g(\vbar) + \varphi(\zbar) \\[5pt]
   \mbox{s.t.} & \cF^*u  + \cA_{\leq l}^*y_{\leq l} \;=\; \cbar - \cA_{l+1}^*\ybar_{l+1}.
   \label{prop-eq-l2}
   \end{array}
\end{eqnarray}
The augmented Lagrangian function associated with problem \eqref{prop-eq-l2} is defined by
\[
\begin{array}{lcl}
\widetilde{\cL}_{\sigma}(u,y_{\leq l};\ybar_{l+1},\vbar,\zbar,x)
& = &\widetilde{f}_1(u) + \sum_{i=1}^l\widetilde{\theta}_{i}(y_{i}) + \theta_{l+1}(\ybar_{l+1}) + g(\vbar) +\varphi(\zbar)
\\[8pt]
& & +\, \inprod{x}{\Gamma(u,(y_{\leq l},\ybar_{l+1}),\vbar,\zbar)}  +
 \frac{\sigma}{2}\norm{\Gamma(u,(y_{\leq l},\ybar_{l+1}),\vbar,\zbar)}^2.
 \end{array}
\]
Define
\[
\cT_{\widetilde{\theta_i}} \equiv  \cT_{{\theta_i}} \quad {\rm and} \quad \cT_{\widetilde {f}_i} \equiv \cT_{f_i},\quad i=1, \ldots, l.
\]
By using the definitions of $\widetilde{\theta}_i$ and $\widetilde{f}_i$, $i=1, \ldots, l$,  we have
\begin{equation}\label{eq:tilde}
\cE_{\widetilde{{\theta_i}}} \equiv \cE_{\theta_i}  \quad {\rm and} \quad \widehat{\cT}_{\widetilde{f}_i} \equiv  \widehat{\cT}_{{f}_i}, \quad i=1, \ldots, l.
\end{equation}
Therefore,    problem \eqref{padmmM1uy}   can  equivalently be rewritten as
\begin{eqnarray}
(u^+,y^+_{\leq l})  &=& \textup{argmin}_{(u,y_{\leq l})}
\left\{ \begin{array}{l}
\widetilde{\cL}_\sigma(u,y_{\leq l};\ybar_{l+1},\vbar,\zbar,\xbar)
\\[5pt]
+\frac{\sigma}{2}\norm{(u,y_{\leq l-1}) - (\ubar,\ybar_{\leq l-1})}_{\widehat{\cT}_{\widetilde{f}_l}}^2
  + \frac{\sigma}{2}\norm{y_{l} - \ybar_{l}}_{\cT_{\widetilde{\theta}_l}}^2
\end{array}\right\}.
\label{uyplus}
\end{eqnarray}
Define
\[
\widetilde{\beta} _{l,j}\; : = \; \cA_{j-1}\cA_l^*\cE_{\widetilde{\theta}_l}^{-1}(\widetilde{b}_l -\cA_l \xbar -\cP_l \ybar_l +
\sigma\cA_l\gammabar), \quad j=1, \ldots, l
\]
and for  $i=l-1, l-2, \ldots, 1$,
\begin{eqnarray*}
 & \widetilde{\beta} _{i,j}\; : =\; \cA_{j-1}\cA_i^*\cE_{\widetilde{\theta}_i}^{-1}\left (\widetilde{b}_i -\sum _{k=i+1}^l\widetilde{\beta}_{k,i+1}-\cA_i \xbar -\cP_i \ybar_i
+
\sigma\cA_i\gammabar\right), \quad j=1, \ldots, i. &
\end{eqnarray*}
The
 auxiliary linear term  ${\delta}_{\widetilde \theta}$ associated with problem (\ref{uyplus}) is given by
\begin{eqnarray}\label{eq:linear term-wd}
  &
 {\delta}_{\widetilde \theta} \; := \;
\sum _{i=1}^l\widetilde{\beta}_{i,1}.
&
\end{eqnarray}
We will show that for $ i=l,l-1,\ldots,1$,
\begin{eqnarray}
\widetilde{\beta}_{i,j} = \betabar_{i,j} \quad \forall\, j=1,\ldots,i.
\label{r:wbetaij-eq-betaij}
\end{eqnarray}
First, by using (\ref{eq:tilde}), we have   for  $j=1,\ldots,l$ that
\begin{eqnarray*}
\widetilde{\beta}_{l,j}  &=& \cA_{j-1}\cA_{l}^*\cE_{\widetilde{\theta}_{l}}^{-1}(\widetilde{b}_{l} -\cA_{l}\xbar -\cP_{l}\ybar_{l}  + \sigma\cA_{l}\gammabar) \\[5pt]
&=& \cA_{j-1}\cA_{l}^*\cE_{\theta_{l}}^{-1}(b_{l} -\betabar_{l+1,l+1}-\cA_{l}\xbar  -\cP_{l}\ybar_{l} + \sigma\cA_{l}\gammabar) \;=\; \betabar_{l,j}.
\end{eqnarray*}
That is, (\ref{r:wbetaij-eq-betaij}) holds for $i=l$ and $j=1,\ldots,l$.
Now assume that we have proven  $\widetilde{\beta}_{i,j} = \betabar_{i,j}$ for all
$i \geq k+1$ with $k+1\le l$  and $j = 1,\ldots,i$.
We shall next prove that \eqref{r:wbetaij-eq-betaij} holds for  $i=k$ and $j=1,\ldots, k$. Again, by using (\ref{eq:tilde}), we have   for $j = 1,\ldots,k$ that
\begin{eqnarray*}
  \widetilde{\beta}_{k,j}  &=& \cA_{j-1}\cA_{k}^*\cE_{\widetilde{\theta}_{k}}^{-1}
\Big(\widetilde{b}_{k}- \mbox{$\sum_{s=k+1}^l$} \widetilde{\beta}_{s,k+1} -\cA_{k}\xbar -\cP_{k}\ybar_{k}  + \sigma\cA_{k}\gammabar \Big) \\
&=& \cA_{j-1}\cA_{k}^*\cE_{\theta_{k}}^{-1}
\Big(b_{k} - \betabar_{p,k+1} - \mbox{$\sum_{s=k+1}^{l}$}
\betabar_{s,k+1} -\cA_{k}\xbar -\cP_{k}\ybar_{k}  + \sigma\cA_{k}\gammabar\Big) \\
&=& \cA_{j-1}\cA_{k}^*\cE_{\theta_{k}}^{-1}\Big(b_{k} - \mbox{$\sum_{s=k+1}^{l+1}$}
\betabar_{s,k+1}-\cA_{k}\xbar-\cP_{k}\ybar_{k}  + \sigma\cA_{k}\gammabar\Big)
\;\;=\;\; \betabar_{k,j},
\end{eqnarray*}
which, shows that \eqref{r:wbetaij-eq-betaij} holds for  $i=k$ and $j=1,\ldots, k$.
Thus, \eqref{r:wbetaij-eq-betaij} is proven.

For $i=l,l-1,\ldots,1,$  define $\widetilde{y}^\prime_i\in \cY_i$ by
\begin{eqnarray}
  \label{wty-prime-i}
   \widetilde{y}^\prime_i &:=& \textup{argmin}_{y_i} \;\widetilde{\cL}_{\sigma}(\ubar,(\ybar_{\leq i-1},y_i,\widetilde{y}^\prime_{\geq i+1});\ybar_{l+1},\vbar,\zbar,\xbar) + \frac{\sigma}{2}\norm{y_i - \ybar_i}_{\cT_{\widetilde{\theta}_i}}^2, \nn \\[5pt]
   &=& \cE_{\widetilde{\theta}_i}^{-1}\big(\sigma^{-1}\widetilde{b}_i - \sigma^{-1}\cA_i\xbar +
   \cT_{\widetilde{\theta}_i}\ybar_i + \cA_i\cA_i^*\ybar_i - \cA_i\Gamma(\ubar,(\ybar_{\leq i-1},\ybar_i,\widetilde{y}^\prime_{\geq i+1},\ybar_{l+1}),\vbar,\zbar)\big),
\end{eqnarray}
where we use the convention $\widetilde{y}^\prime_{l+1} = \emptyset.$
We will prove that
\begin{eqnarray}\label{wyp-eq-yp}
\widetilde{y}^\prime_{i} = y^\prime_i \quad \forall\, i=l,l-1,\ldots,1.
\end{eqnarray}
We first calculate
\begin{eqnarray}\label{yp-byp}
y^\prime_{l+1} - \ybar_{l+1} &=& \cE_{\theta_{l+1}}^{-1}(\sigma^{-1}b_{l+1} - \sigma^{-1}\cA_{l+1}\xbar + \cT_{\theta_{l+1}}\ybar_{l+1} + \cA_{l+1}\cA^*_{l+1} \ybar_{l+1} +\cA_{l+1}\gammabar - \cE_{\theta_{l+1}}\ybar_{l+1}) \nn\\[5pt]
&=& \cE_{\theta_{l+1}}^{-1}(\sigma^{-1}b_{l+1} - \sigma^{-1}\cA_{l+1}\xbar-\sigma^{-1}\cP_{l+1}\ybar_{l+1} + \cA_{l+1}\gammabar),
\end{eqnarray}
which, together with the definitions of $\beta_{p,i}$ in \eqref{betapj}, implies
\begin{eqnarray}
\label{betapip1}
  \cA_i\cA_{l+1}^*(y^\prime_{l+1} - \ybar _{l+1}) = \sigma^{-1}\beta_{p,i+1} \quad \forall\, i=0,\ldots,l.
\end{eqnarray}
Now, by using (\ref{eq:tilde}), \eqref{betapip1} and the definitions of $\widetilde{y}^\prime_l$ and $y^\prime_l$, we have
\begin{eqnarray*}
y^\prime_l - \widetilde{y}^\prime_l &=& \cE_{\theta_l}^{-1}\big(\sigma^{-1}\beta_{p,l+1} + \cA_l\cA_{l+1}^*(\ybar_{l+1} - y^\prime_{l+1})\big) \\[5pt]
&=& \cE_{\theta_l}^{-1}(\sigma^{-1}\beta_{p,l+1} - \sigma^{-1}\beta_{p,l+1}) = 0.
\end{eqnarray*}
That is, \eqref{wyp-eq-yp} holds for $i=l$. Now assume that we have proven  $\widetilde{y}^\prime_i = y^\prime_i$ for all
$i \geq k+1$ with $k+1\le l$. We shall next prove that \eqref{wyp-eq-yp} holds for  $i=k$.
Again, by using the definitions of $\widetilde{y}^\prime_k$ and $y^\prime_k$ and noting
\[\Gamma(\ubar,(\ybar_{\leq k},\widetilde{y}^\prime_{\geq k+1}, \ybar_{l+1}),\vbar,\zbar) -
\Gamma(\ubar,(\ybar_{\leq k},y^\prime_{\geq k+1}),\vbar,\zbar) = \cA_{l+1}^*(\ybar_{l+1} - y^\prime_{l+1}),\] we obtain that
\begin{eqnarray*}
  y^\prime_k - \widetilde{y}^\prime_k &=& \cE_{\theta_k}^{-1}\big(\sigma^{-1}(b_k - \widetilde{b}_k) + \cA_k\cA_{l+1}^*(\ybar_{l+1} - y^\prime_{l+1})\big) \\[5pt]
  &=& \cE_{\theta_k}^{-1}\big(\sigma^{-1}\beta_{p,k+1} + \cA_k\cA_{l+1}^*(\ybar_{l+1} - y^\prime_{l+1})\big) \\[5pt]
&=& \cE_{\theta_l}^{-1}(\sigma^{-1}\beta_{p,k+1} - \sigma^{-1}\beta_{p,k+1}) = 0,
\end{eqnarray*}
which, shows that \eqref{wyp-eq-yp} holds for $i=k.$ Thus, \eqref{wyp-eq-yp} holds.
%\begin{eqnarray}
%\widetilde{\beta}^{-}_{i,j} = \beta_{i,j}^-, \quad i=l,l-1,\ldots,1,\quad \mbox{and}\quad j=1,2,\ldots,i.
%\label{wbetaij-eq-betaij}
%\end{eqnarray}
% For $i=p-1,p-2,\ldots,1$ and $j=1,\ldots,i$, define the corresponding auxiliary linear terms associated with problem (\ref{prop-eq-l2}) by
%\begin{eqnarray}
% \widetilde{\delta}^{\,-}_{\theta_i,j} &:=& \cA_j\cA_i^*\cE_{\theta_i}^{-1}(\tilde{b}_{i} -\cA_{i}x^--\cP_iy_i^- - \sum_{s = i+1}^{p-1} \widetilde{\delta}^{\,-}_{\theta_s,i} -
%\sigma\cA_i\gamma^-).
%  \label{betaij}
%\end{eqnarray}
%In particular, when $i=p-1$ and $j=0,1,\ldots,p-2$, we have
%\begin{eqnarray}
% \widetilde{\delta}^{\,-}_{\theta_{p-1},j} &:=& \cA_j\cA_{p-1}^*\cE_{\theta_{p-1}}^{-1}(\tilde{b}_{p-1} -\cA_{p-1}x^--\cP_{p-1}y_{p-1}^- -
%\sigma\cA_{p-1}\gamma^-).
%  \label{betaij}
%\end{eqnarray}
%Denote  $\widetilde{\delta}_{\theta}^{\,-} \equiv  \sum_{i=1}^{p-1}\widetilde{\delta}_{\theta_i,0}^{\,-}$.

By applying our induction assumption  to problem (\ref{uyplus}), we obtain equivalently that
\begin{eqnarray}
    u^+ &=& \textup{argmin}_{u}\;
\widetilde{\cL}_\sigma(u,\ybar_{\leq l};\ybar_{l+1},\vbar,\zbar,\xbar) +  \inprod{ {\delta}_{\widetilde{\theta}}}{u} +
  \frac{\sigma}{2}\norm{u - \ubar}_{\cT_f}^2,\label{prop-eq-l2-u}
\\[5pt]
  y_i^+  &=& \textup{argmin}_{y_i} \;\widetilde{\cL}_{\sigma}(u^+,(y_{\leq i-1}^+,y_i,\widetilde{y}^\prime_{\geq i+1}); \ybar_{l+1},\vbar,\zbar,\xbar)
   + \frac{\sigma}{2}\norm{y_i - \ybar_i}_{\cT_{\theta_i}}^2,\quad i=1,\ldots,l,
\label{prop-eq-l2-y}
\end{eqnarray}
where we use the facts that  $\cT_{\widetilde{f}_1} =\cT_f$ and $\cT_{\widetilde{\theta}_i} =\cT_{\theta_i}$ for $i=1, \ldots, l$.
By combining \eqref{r:wbetaij-eq-betaij} and the definitions of $\deltabar_\theta$ and  ${\delta}_{\widetilde \theta}$ defined in \eqref{eq:linear term} and
\eqref{eq:linear term-wd}, respectively, we derive that
\begin{eqnarray}\label{eq:delta-equivalence}
&  \deltabar_\theta =\sum _{i=1}^{l}\betabar_{i,1} + \betabar_{l+1,1}= \sum _{i=1}^l\widetilde{\beta}_{i,1} + \betabar_{l+1,1}
={\delta}_{\widetilde \theta}+\betabar_{l+1,1}.&
\end{eqnarray}
By direct calculations,
\begin{eqnarray}\label{wcl-eq-cl-p-A}
&&\hspace{-0.7cm}
\widetilde{\cL}_\sigma(u,\ybar_{\leq l};\ybar_{l+1},\vbar,\zbar,\xbar) = \cL_\sigma(u,\ybar,\vbar,\zbar;\xbar) + \inprod{\betabar_{l+1,1}}{u} +
\sum_{i=1}^l \inprod{\betabar_{l+1,i+1}}{\ybar_i}.
\end{eqnarray}
 Using \eqref{wyp-eq-yp}, \eqref{betapip1} and the definition of $\widetilde{\cL}_\sigma$, we have for $ i=1,\ldots,l$ that
\begin{eqnarray}\label{wcl-eq-cl-p}
&&
\widetilde{\cL}_{\sigma}(u^+,(y_{\leq i-1}^+,y_i,\widetilde{y}^\prime_{\geq i+1});\ybar_{l+1},\vbar,\zbar,\xbar) - \cL_{\sigma}(u^+,(y_{\leq i-1}^+,y_i,y^\prime_{\geq i+1}),\vbar,\zbar;\xbar) \nn\\[5pt]
&=& \widetilde{\cL}_{\sigma}(u^+,(y_{\leq i-1}^+,y_i,y^\prime_{ i+1},\ldots,y^\prime_l);\ybar_{l+1},\vbar,\zbar,\xbar) - \cL_{\sigma}(u^+,(y_{\leq i-1}^+,y_i,y^\prime_{\geq i+1}),\vbar,\zbar;\xbar) \nn\\[5pt]
&=& \inprod{\beta_{p,i+1}}{y_i} + \inprod{\sigma\cA_i\cA_{l+1}^*(\ybar_{l+1} - y^\prime_{l+1})}{y_i} + c_i \nn\\[5pt]
&=& c_i,
 %+\inprod{\betabar_{l+1,1}}{u^+},
%\\[5pt]
%&& \hspace{4cm}\fbox{+$\inprod{\betabar_{l+1,i+1}}{y_i-\ybar_i}
%+\sum_{j=1}^{i-1} \inprod{\betabar_{l+1,j+1}}{y^+_j}
%+\sum_{j=i+1}^{l} \inprod{\betabar_{l+1,j+1}}{\ybar_j}$??}
\end{eqnarray}
where $c_i$ is a constant term given by
\begin{eqnarray*}
c_i&=&\inprod{\betabar_{l+1,1}}{u^+} + \sum^{i-1}_{j=1}\inprod{\betabar_{l+1,j+1}}{y^+_j}
+ \sum_{j=i+1}^{l} \inprod{\betabar_{l+1,j+1}}{y^\prime_j} \\[5pt]
&&+ \theta_{l+1}(\ybar_{l+1}) - \theta_{l+1}(y^\prime_{l+1}) + \inprod{\xbar}{\cA_{l+1}^*(\ybar_{l+1} - y^\prime_{l+1})} \\[5pt]
&& + \frac{\sigma}{2}\inprod{\cA_{l+1}^*(\ybar_{l+1} - y^\prime_{l+1})}{2(\cF^*u^+ + \cA_{\le i-1}^* y^+_{\le i-1} + \sum_{j=i+1}^l\cA_{j}^* y^\prime_{j} - \cbar) + \cA_{l+1}^*(\ybar_{l+1} + y^\prime_{l+1})}.
\end{eqnarray*}
Thus, by using  (\ref{eq:delta-equivalence}), \eqref{wcl-eq-cl-p-A} and \eqref{wcl-eq-cl-p} we know that
\eqref{prop-eq-l2-u} and \eqref{prop-eq-l2-y} can be rewritten as
 \begin{eqnarray*}
 \left\{ \begin{array}{l}
  u^+ = \textup{argmin}_{u}\; \cL_\sigma(u,\ybar,\vbar,\zbar;\xbar) + \inprod{\deltabar_{\theta}}{u}+
  \frac{\sigma}{2}\norm{u - \ubar}_{\cT_f}^2,
\\[5pt]
  y_i^+ = \textup{argmin}_{y_i} \;\cL_{\sigma}(u^+,(y_{\leq i-1}^+,y_i,
 y^\prime_{\geq i+1}),\vbar,\zbar;\xbar) + \frac{\sigma}{2}\norm{y_i - \ybar_i}_{\cT_{\theta_i}}^2, \quad i=1,\ldots,l,
 \end{array} \right.
 \label{padmmMll}
 \end{eqnarray*}
 which, together with \eqref{prop-eq-yp}, shows that
  the equivalence   between \eqref{prox-cM1} and \eqref{padmmMl}
holds for  $p = l+1$.     The proof of this part is completed.

\medskip
{\bf Part two.} In this part, we prove the equivalence between \eqref{padmmMl} and \eqref{padmmMy}. Again, for the case $p=1$, it follows directly from
  Proposition \ref{prop:prox2}.

Assume that  the equivalence  between \eqref{padmmMl} and \eqref{padmmMy}
holds for all $p \leq l$. We shall prove that  this equivalence also holds for $p=l+1$. Write
${f}_0(\cdot) \equiv f(\cdot) + \sum_{i=1}^l\inprod{\betabar_{i,1}}{\cdot}.$
Since ${f}_0$ differs from $f$ only with an extra linear term, we define $\cT_{{f_0}} \equiv \cT_f.$
In order to use Proposition \ref{prop:prox2},
we consider the following optimization problem with respect to $u$ and $ y_{l+1}$:
\begin{eqnarray}
    \begin{array}{rlll}
    \min & {f}_0(u) + \theta_{l+1}(y_{l+1}) + \sum_{i=1}^l \theta_i(\ybar_i) + g(\vbar) + \varphi(\zbar) \\[5pt]
   \mbox{s.t.} & \cF^*u  + \cA_{l+1}^*y_{l+1} = \cbar -\cA^*_{\leq l}\ybar_{\leq l}.
   \label{prop-eq-l3}
   \end{array}
\end{eqnarray}
The augmented Lagrangian function associated with problem \eqref{prop-eq-l3} is given as follows:
\begin{eqnarray*}
{\cL}_{\sigma}^0(u,y_{l+1};\ybar_{\leq l}, \vbar,\zbar,x) &=&
{f}_0(u) + \theta_{l+1}(y_{l+1}) + \mbox{$\sum_{i=1}^l$}
 \theta_i(\ybar_i)+ g(\vbar) +\varphi(\zbar) \nn\\[5pt]
&&+\inprod{x}{\Gamma(u,(\ybar_{\leq l},y_{l+1}),\vbar,\zbar)} + \frac{\sigma}{2}\norm{\Gamma(u,(\ybar_{\leq l},y_{l+1}),\vbar,\zbar)}^2.
\end{eqnarray*}
By observing that
 \begin{eqnarray*}
 &{\cL}^0_{\sigma}(u,\ybar_{l+1};\ybar_{\leq l}, \vbar,\zbar,\xbar)
=  \cL_\sigma(u,\ybar,\vbar,\zbar;\xbar) + \sum_{i=1}^l\inprod{ \betabar_{i,1}}{u} \quad {\rm and} \quad \cT_{{f}_0} \equiv \cT_f,&
\end{eqnarray*}
  we can rewrite the first subproblem in \eqref{padmmMl} as
 \begin{eqnarray}
   u^+ =
\mbox{argmin}_{u}\; {\cL}^0_{\sigma}(u,\ybar_{l+1};\ybar_{\leq l}, \vbar,\zbar,\xbar) + \inprod{ \betabar_{l+1,1}}{u} +\frac{\sigma}{2}\norm{u-\ubar}_{\cT_{{f}_0}}^2.
  \label{tilde-uv-l}
\end{eqnarray}

{} By using   the definition of $\yprime_{l+1}$ given in \eqref{y-prime-i}, we have
\begin{eqnarray}
\yprime_{l+1} = \cE_{\theta_{l+1}}^{-1}\big(\sigma^{-1}( b_{l+1}  - \cA_{l+1} \xbar) + \cT_{\theta_{l+1}}\ybar_{l+1} + \cA_{l+1}\cA_{l+1}^*\ybar_{l+1} + \cA_{l+1}\gammabar\big).
\label{by-lplus1}
\end{eqnarray}
Since
\begin{eqnarray*}
&{\cL}^0_{\sigma}(\ubar,y_{l+1};\ybar_{\leq l}, \vbar,\zbar,\xbar) = \cL_{\sigma}(\ubar,(\ybar_{\leq l}, y_{l+1}),\vbar,\zbar;\xbar) + \sum_{i=1}^l \inprod{\betabar_{i,1}}{\ubar},&
\end{eqnarray*}
the point $\yprime_{l+1}$   can be  rewritten equivalently as
 \begin{eqnarray} \yprime_{l+1}  =  \textup{argmin}_{y_{l+1}}\; {\cL}^0_{\sigma}(\ubar,y_{l+1};\ybar_{\leq l}, \vbar,\zbar,\xbar)
  + \frac{\sigma}{2}\norm{y_{l+1} - \ybar_{l+1}}^2_{\cT_{\theta_{l+1}}}. \label{padmmMbary}
  \end{eqnarray}
%\begin{eqnarray}
 % \yprime_{l+1} \;=\; \textup{argmin}_{y_{l+1}} \;\cL_{\sigma}(\ubar,\ybar_{\leq l}, y_{l+1},\vbar,\zbar;\xbar) + \frac{\sigma}{2}\norm{y_{l+1} - %\ybar_{l+1}}_{\cT_{\theta_{l+1}}}^2.
%\end{eqnarray}
Then, by applying Proposition \ref{prop:prox2} to problem  (\ref{prop-eq-l3}) with respect to $u$ and $ y_{l+1}$,
 we know that problem \eqref{tilde-uv-l} is  equivalent to
 \begin{eqnarray}
  u^+  =  \textup{argmin}_{u }\; {\cL}^0_{\sigma}(u,\yprime_{l+1};\ybar_{\leq l}, \vbar,\zbar,\xbar) +
  \frac{\sigma}{2}\norm{u - \ubar}_{\cT_{{f}_0}}^2.
 \label{padmmMup}
 \end{eqnarray}

In order to apply our  induction assumption to problem \eqref{padmmMup}, we need to consider the following optimization problem with respect to $(u,y_{\le l})$:
\begin{eqnarray}
    \begin{array}{rlll}
    \min & f(u) + \sum_{i=1}^l \theta_{i}(y_{i}) +\theta_{l+1}(\yprime_{l+1}) + g(\vbar) + \varphi(\zbar) \\[5pt]
   \mbox{s.t.} & \cF^*(u)  + \cA_{\leq l}^*y_{\leq l} = \cbar - \cA_{l+1}^*\yprime_{l+1}  .
   \label{prop-eq-l4}
   \end{array}
\end{eqnarray}
The augmented Lagrangian function associated with problem  \eqref{prop-eq-l4} is
given by
\begin{eqnarray}
\widehat{\cL}_{\sigma}(u,y_{\leq l};\yprime_{l+1},\vbar,\zbar,x)
&=&  f(u) + \mbox{$\sum_{i=1}^l$} \theta_{i}(y_i) + \theta_{l+1}(\yprime_{l+1}) + g(\vbar) +\varphi(\zbar) \nn\\[5pt]
&&+\inprod{x}{\Gamma(u,(y_{\leq l},\yprime_{l+1}),\vbar,\zbar)} + \frac{\sigma}{2}\norm{\Gamma(u,(y_{\leq l},\yprime_{l+1}),\vbar,\zbar)}^2.\nn
\end{eqnarray}
Define
\[\widehat{\gamma} \; := \; -\Gamma(\ubar,(\ybar_{\leq l},\yprime_{l+1}),\vbar,\zbar)
\quad {\rm and } \quad h_i \; := \; b_i - \cA_i \xbar - \cP_i \ybar_i, \quad i=1,\ldots,l.\]
 For problem \eqref{prop-eq-l4},  we define the following associated     terms
\[
\widehat{\beta} _{l,j}\; : = \; \cA_{j-1}\cA_l^*\cE_{\theta_l}^{-1}(h_l +
\sigma\cA_l\widehat{\gamma}), \quad j=1, \ldots, l
\]
and for  $i=l-1, l-2, \ldots, 1$,
\[
\widehat{\beta} _{i,j}\; : = \; \cA_{j-1}\cA_i^*\cE_{\theta_i}^{-1}\Big (h_i -
 \mbox{$\sum _{k=i+1}^l$} \widehat{\beta}_{k,i+1} +
\sigma\cA_i\widehat{\gamma}\Big), \quad j=1, \ldots, i.
\]
The auxiliary linear term $\widehat{\delta}$ associated with problem \eqref{prop-eq-l4}
is given by
\begin{eqnarray}\label{eq:linear term-p2}
  &\widehat{\delta} \;=\; \sum _{i=1}^l\widehat{\beta}_{i,1}.  &
\end{eqnarray}
We will show that, for $i=l,l-1,\ldots,1$,
\begin{eqnarray}
\widehat{\beta}_{i,j}\; = \;\betabar_{i,j} \quad \forall\,  j=1,\ldots,i.
\label{r:wbetaij-eq-betaij-p2}
\end{eqnarray}
Similar to what  we have done in  part one, we shall first prove that $\widehat{\beta}_{l,j} = \betabar_{l,j}$
for $j=1,2,\ldots,l$.
In fact, for  $j=1,\ldots,l$, we have
\begin{eqnarray*}
 \betabar_{l,j}  &=& \cA_{j-1}\cA_{l}^*\cE_{\theta_{l}}^{-1}(h_l - \betabar_{l+1,l+1}  +\sigma\cA_{l}\gammabar) \\[5pt]
&=& \cA_{j-1}\cA_{l}^*\cE_{\theta_{l}}^{-1}(h_l - \cA_l\cA_{l+1}^*\cE_{\theta_{l+1}}^{-1}(h_{l+1} + \sigma\cA_{l+1}\gammabar) + \sigma\cA_{l}\gammabar) \\[5pt]
&=& \cA_{j-1}\cA_{l}^*\cE_{\theta_{l}}^{-1}(h_{l} -
\sigma\cA_{l}\Gamma(\ubar,(\ybar_{\leq l},\yprime_{l+1}),\vbar,\zbar)) \\[5pt]
&=& \cA_{j-1}\cA_l^*\cE_{\theta_l}^{-1}(h_l +
\sigma\cA_l\widehat{\gamma}) = \widehat{\beta}_{l,j},
\end{eqnarray*}
where the  third equation  follows from \eqref{by-lplus1} and simple  calculations.
This shows that (\ref{r:wbetaij-eq-betaij-p2}) holds for  $i=l$ and $j=1,\ldots,l$.
Now we assume that   $\widehat{\beta}_{i,j} = \betabar_{i,j}$  for all $i \geq k+1$ with $k+1\le l$  and $j = 1,\ldots,i$. Next,    we shall  prove that  (\ref{r:wbetaij-eq-betaij-p2}) holds for
 $i=k$ and $j=1,\ldots, k$.
By direct calculations, we know for   $j = 1,\ldots,k$  that
\begin{eqnarray*}
\betabar_{k,j} &=& \cA_{j-1}\cA_{k}^*\cE_{\theta_{k}}^{-1}\Big(h_{k}
-\mbox{$\sum_{s=k+1}^{l+1}$} \betabar_{s,k} +\sigma\cA_{k}\gammabar\Big) \\[5pt]
&=& \cA_{j-1}\cA_{k}^*\cE_{\theta_{k}}^{-1}\Big(h_{k}
- \mbox{$\sum_{s=k+1}^{l}$} \widehat{\beta}_{s,k} -\betabar_{l+1,k} +\sigma\cA_{k}\gammabar\Big) \\[5pt]
& = & \cA_{j-1}\cA_{k}^*\cE_{\theta_{k}}^{-1}
\Big(h_{k}
- \mbox{$\sum_{s=k+1}^{l}$} \widehat{\beta}_{s,k} -\cA_{k}\cA_{l+1}^*\cE_{\theta_{l+1}}^{-1}(h_{l+1} + \sigma\cA_{l+1}\gammabar) +\sigma\cA_{k}\gammabar\Big) \\[5pt]
&=& \cA_{j-1}\cA_{k}^*\cE_{\theta_{k}}^{-1}
\Big(h_{k}
-\mbox{$\sum_{s=k+1}^{l}$} \widehat{\delta}_{\theta_{s},k}
-\sigma\cA_{k}\Gamma(\ubar,(\ybar_{\leq l},y^\prime_{l+1}),\vbar,\zbar)\Big)
\\[5pt]
&=& \cA_{j-1}\cA_{k}^*\cE_{\theta_{k}}^{-1}
\Big(h_{k}
-\mbox{$\sum_{s=k+1}^{l}$} \widehat{\delta}_{\theta_{s},k} +\sigma\cA_{k}\widehat{\gamma}\Big)
\;=\; \widehat{\beta}_{k,j},
\end{eqnarray*}
which, shows that \eqref{r:wbetaij-eq-betaij-p2} holds for  $i=k$ and $j=1,\ldots, k$.
Therefore, we have shown that \eqref{r:wbetaij-eq-betaij-p2} holds.

For $i=l,l-1,\ldots,1,$ define $\widehat{y}^\prime_i\in \cY_i$  as
\begin{eqnarray}
  \label{whty-prime-i}
   \widehat{y}^\prime_i &=& \textup{argmin}_{y_i} \;\widehat{\cL}_{\sigma}(\ubar,(\ybar_{\leq i-1},y_i,\widehat{y}^\prime_{\geq i+1});y^\prime_{l+1},\vbar,\zbar,\xbar) + \frac{\sigma}{2}\norm{y_i - \ybar_i}_{\cT_{ \theta_i}}^2 \nn \\[5pt]
   &=& \cE_{\theta_i}^{-1}\big(\sigma^{-1} b _i - \sigma^{-1}\cA_i\xbar +
   \cT_{ \theta _i}\ybar_i + \cA_i\cA_i^*\ybar_i - \cA_i\Gamma(\ubar,(\ybar_{\leq i-1},\ybar_i,\widehat{y}^\prime_{\geq i+1},y^\prime_{l+1}),\vbar,\zbar)\big),
\end{eqnarray}
where we use the convention $\widehat{y}^\prime_{l+1} = \emptyset.$
We will prove that
\begin{eqnarray}
  \label{whyp-eq-yp-p2}
  \widehat{y}^\prime_i = y^\prime_i  \quad \forall\, i=1,\ldots,l.
\end{eqnarray}
{}From \eqref{whty-prime-i}, we know that
\[\widehat{y}^\prime_l= \cE_{\theta_l}^{-1}\big(\sigma^{-1} b _l - \sigma^{-1}\cA_l\xbar +
   \cT_{ \theta _l}\ybar_l + \cA_l\cA_l^*\ybar_l - \cA_l\Gamma(\ubar,(\ybar_{\leq i-1},\ybar_l,y^\prime_{l+1}),\vbar,\zbar)\big),\]
which is exactly the same as $y^\prime_l$ defined   in \eqref{y-prime-i}. This shows that \eqref{whyp-eq-yp-p2} holds for $i=l.$ Now we assume that $\widehat{y}^\prime_i = y^\prime_i$ for all $i\ge k+1$ with $k+1\le l$. Next, we shall prove that \eqref{whyp-eq-yp-p2} holds for $i=k.$ Again, by using the definition of $\widehat{y}^\prime_k$ in \eqref{whty-prime-i} and the definition of $y^\prime_k$ in \eqref{y-prime-i}, we see that
\begin{eqnarray*}
\widehat{y}^\prime_k&=& \cE_{\theta_k}^{-1}\big(\sigma^{-1} b _k - \sigma^{-1}\cA_k\xbar +
   \cT_{ \theta _k}\ybar_k + \cA_k\cA_k^*\ybar_k - \cA_k\Gamma(\ubar,(\ybar_{\leq k-1},\ybar_k,\widehat{y}^\prime_{\geq k+1},y^\prime_{l+1}),\vbar,\zbar)\big) \\
   &=& \cE_{\theta_k}^{-1}\big(\sigma^{-1} b _k - \sigma^{-1}\cA_k\xbar +
   \cT_{ \theta _k}\ybar_k + \cA_k\cA_k^*\ybar_k - \cA_k\Gamma(\ubar,(\ybar_{\leq k-1},\ybar_k,y^\prime_{\geq k+1}),\vbar,\zbar)\big) \\
   &=& y^\prime_k.
   \end{eqnarray*}
Thus, \eqref{whyp-eq-yp-p2} is proven to be true.

By direct calculations, we obtain from  (\ref{eq:linear term-p2}) and \eqref{r:wbetaij-eq-betaij-p2}
 that
\begin{eqnarray}
&  {\cL}^0_{\sigma}(u,y^\prime_{l+1};\ybar_{\leq l}, \vbar,\zbar,\xbar)
- \widehat{\cL}_{\sigma}(u,\ybar_{\leq l};y^\prime_{l+1},\vbar,\zbar,\xbar)
\;=\; \sum_{i=1}^{l} \inprod{ \betabar_{i,1}}{u}
\;=\; \inprod{ \widehat{\delta}}{u}.&
  \label{hatl-eq-tildel-p2}
\end{eqnarray}
 By using  \eqref{hatl-eq-tildel-p2}  and  $\cT_{{f}_0} \equiv \cT_f$, we can reformulate problem  \eqref{padmmMup}  equivalently as
\begin{eqnarray}
u^+
 =  \textup{argmin}_{u}\; \widehat{\cL}_{\sigma}(u,\ybar_{\leq l};y^\prime_{l+1},\vbar,\zbar,\xbar) +  \inprod{ \widehat{\delta}}{u}+ \frac{\sigma}{2}\norm{u - \ubar}_{\cT_f}^2.
\label{padmmM2u}
\end{eqnarray}
Then, from our induction assumption we know that problem \eqref{padmmM2u} can be equivalently recast as
\begin{eqnarray}
\left\{\begin{array}{lcl}
 \widehat{y}^\prime_i &=& \textup{argmin}_{y_i} \;\widehat{\cL}_{\sigma}(\ubar,(\ybar_{\leq i-1},y_i,\widehat{y}^\prime_{\geq i+1});\yprime_{l+1},\vbar,\zbar,\xbar)
  + \frac{\sigma}{2}\norm{y_i - \ybar_i}_{\cT_{\theta_i}}^2, \quad i=l,l-1,\ldots,1,\label{hatl-eq-tildel-p22} \\[5pt]
  u^+ &=& \textup{argmin}_{u}\; \widehat{\cL}_\sigma(u,{\widehat{y}}^\prime_{\leq l};\yprime_{l+1},\vbar,\zbar,\xbar) +
  \frac{\sigma}{2}\norm{u - \ubar}_{\cT_f}^2.
  \end{array}\right.
 \end{eqnarray}
By using \eqref{whyp-eq-yp-p2} and observing  \[\widehat{\cL}_{\sigma}(u,y_{\leq l};\yprime_{l+1},\vbar,\zbar,\xbar)
\;=\; \cL_{\sigma}(u,y_{\leq l},\yprime_{l+1},\vbar,\zbar;\xbar), \]
  we know that  \eqref{hatl-eq-tildel-p22} is equivalent to
 \begin{eqnarray}
 \left\{\begin{array}{lcl}
 \yprime_i &=& \textup{argmin}_{y_i} \;\cL_{\sigma}(\ubar,(\ybar_{\leq i-1},y_i,
\yprime_{\geq i+1}),\vbar,\zbar;\xbar)
+ \frac{\sigma}{2}\norm{y_i - \ybar_i}_{\cT_{\theta_i}}^2, \quad i=l,l-1,\ldots,1,\nn
\\[5pt]
  u^+ &=& \textup{argmin}_{u}\; \cL_\sigma(u,(\yprime_{
  \leq l},\yprime_{l+1}),\vbar,\zbar;\xbar) +
  \frac{\sigma}{2}\norm{u - \ubar}_{\cT_f}^2,
  \end{array}\right.
 \end{eqnarray}
 which, together with \eqref{padmmMbary}, shows that
  the equivalence   between \eqref{padmmMl} and \eqref{padmmMy}
holds for  $p = l+1$.     This completes the proof to the second part of this proposition.   \qed
%The proof of this proposition is completed.
%Let $\tilde{\delta}_g := \widetilde{\cF}\widetilde{\cG}^*\cE_{\tilde{g}}^{-1}(\tilde{b}-\widetilde{{\cG}} x^- - \Sigma_{\tilde{g}} \tilde{v}^- + \sigma\widetilde{\cG} \tilde{c}^-),$ we have

%\medskip

%\bigskip
%\blue{\fbox{TKC: I have checked up to this point.}}

\bigskip
\begin{prop}\label{prop:equi-scb-padmm}
For  any $k\ge 0$, the point $(x^{k+1}, y^{k+1},v^{k+1}, z^{k+1})$ obtained by Algorithm SCB-SPADMM
for solving problem (\ref{eq-dd-M}) can be generated exactly according to the following iteration:
  \begin{eqnarray}
  \left\{ \begin{array}{l}
  (u^{k+1},y^{k+1}) = \textup{argmin}_{u,y}\; \cL_{\sigma}(u,y,v^k,z^k;x^k)
+\frac{\sigma}{2}\norm{(u,y_{\leq p-1}) - (u^k,y_{\leq p-1}^k)}_{\widehat{\cT}_{f_p}}^2
  + \frac{\sigma}{2}\norm{y_p-y_p^k}_{\cT_{\theta_p}}^2, \nn \\[5pt]
  (v^{k+1},z^{k+1}) = \textup{argmin}_{v,z}\; \cL_{\sigma}(u^{k+1},y^{k+1},v,z;x^k)
  +\frac{\sigma}{2}\norm{(v,z_{\leq q-1}) - (v^k,z_{\leq q-1}^k)}_{\widehat{\cT}_{g_q}}^2
  + \frac{\sigma}{2}\norm{z_q-z_q^k}_{\cT_{\varphi_q}}^2, \nn \\[5pt]
    x^{k+1} = x^k + \tau \sigma(\cF^*u^{k+1} +  \cA^*y^{k+1} + \cG^*v^{k+1} +  \cB^*z^{k+1} -c).
    \end{array}
    \right.
  \end{eqnarray}
\end{prop}
{\bf Proof.} The $(u^{k+1},y^{k+1})$ part directly follows from Proposition \ref{prop:proxaugM}. The conclusion for the $(v^{k+1},z^{k+1})$ part can be obtained in  similar arguments to the part about $(u^{k+1},y^{k+1})$.
Hence, the required result follows.
\qed\bigskip

Write $\Sigma_{f_1} \equiv \Sigma_f$ and $\Sigma_{g_1} \equiv \Sigma_g$.
Define
\[\Sigma_{f_i} := \left(
                                         \begin{array}{cc}
                                           \Sigma_{f_{i-1}} &  \\
                                            & \cP_{i-1} \\
                                         \end{array}
                                       \right), \quad i=2,\ldots,p+1
\]
and
\[
\Sigma_{g_j} := \left(
                                         \begin{array}{cc}
                                           \Sigma_{g_{j-1}} &  \\
                                            & \cQ_{j-1} \\
                                         \end{array}
                                       \right), \quad j=2,\ldots,q+1.
          \]
In order to prove the convergence of our algorithm SCB-SPADMM for solving problem \eqref{eq-dd-M},
we need the following proposition.

\begin{prop}\label{eqvi-psd}
 It holds that
  \begin{eqnarray}
    \cF_{p+1}\cF_{p+1}^* + \sigma^{-1}\Sigma_{f_{p+1}} + \left(
                                               \begin{array}{cc}
                                                 \widehat{\cT}_{f_p} &  \\
                                                  & \cT_{\theta_p} \\
                                               \end{array}
                                             \right)\succ 0 \Leftrightarrow
    \cF\cF^* + \sigma^{-1}\Sigma_f + \cT_f \succ 0, \label{psd-f}\\[5pt]
    \cG_{q+1}\cG_{q+1}^* + \sigma^{-1}\Sigma_{g_{q+1}} + \left(
                                               \begin{array}{cc}
                                                 \widehat{\cT}_{g_q} &  \\
                                                  & \cT_{\varphi_q} \\
                                               \end{array}
                                             \right)\succ 0 \Leftrightarrow
    \cG\cG^* + \sigma^{-1}\Sigma_g + \cT_g \succ 0. \label{pad-g}
  \end{eqnarray}
\end{prop}
{\bf Proof.} We only need to prove \eqref{psd-f} as  \eqref{pad-g} can be obtained in the similar manner. For $i=3,\ldots,p+1$, we have
\begin{eqnarray*}
\cF_i\cF_i^* + \sigma^{-1}\Sigma_{f_i} + \left(
                                               \begin{array}{cc}
                                                 \widehat{\cT}_{f_{i-1}} &  \\
                                                  & \cT_{\theta_{i-1}} \\
                                               \end{array}
                                             \right) = \left(
                                                         \begin{array}{cc}
                                                           \cF_{i-1}\cF_{i-1}^*+\sigma^{-1}
                                                           \Sigma_{f_{i-1}} + \widehat{\cT}_{f_{i-1}} & \cF_{i-1}\cA_{i-1}^* \\
                                                           \cA_{i-1}\cF^*_{i-1} & \cA_{i-1}\cA_{i-1}^* + \sigma^{-1}\cP_{i-1} + \cT_{\theta_{i-1}} \\
                                                         \end{array}
                                                       \right).
\end{eqnarray*}
 Since $\cE_{\theta_{i-1}} = \cA_{i-1}\cA_{i-1}^* + \sigma^{-1}\cP_{i-1} + \cT_{\theta_{i-1}} \succ 0$ for all $i\ge 3$,  by the Schur complement condition for ensuring  the  positive definiteness of linear operators, we have
\begin{eqnarray*}
  &\left(
  \begin{array}{cc}
  \cF_{i-1}\cF_{i-1}^*+\sigma^{-1}
  \Sigma_{f_{i-1}} + \widehat{\cT}_{f_{i-1}} & \cF_{i-1}\cA_{i-1}^* \\
  \cA_{i-1}\cF^*_{i-1} & \cE_{\theta_{i-1}} \\
  \end{array}
  \right) \succ 0 \\ &\Updownarrow \\
  & \cF_{i-1}\cF_{i-1}^* + \sigma^{-1}\Sigma_{f_{i-1}} + \widehat{\cT}_{f_{i-1}} - \cF_{i-1}\cA_{i-1}^*\cE_{\theta_{i-1}}^{-1}\cA_{i-1}\cF^*_{i-1} \succ 0 \\  &\Updownarrow \\
  &\cF_{i-1}\cF_{i-1}^* + \sigma^{-1}\Sigma_{f_{i-1}} + \left(
                                               \begin{array}{cc}
                                                 \widehat{\cT}_{f_{i-2}} &  \\
                                                  & \cT_{\theta_{i-2}} \\
                                               \end{array}
                                             \right)\succ 0.
\end{eqnarray*}
Therefore, by taking $i=3$, we obtain that
\begin{eqnarray*}
  \cF_{p+1}\cF_{p+1}^* + \sigma^{-1}\Sigma_{f_{p+1}} + \left(
                                               \begin{array}{cc}
                                                 \widehat{\cT}_{f_p} &  \\
                                                  & \cT_{\theta_p} \\
                                               \end{array}
                                             \right)\succ 0 \Leftrightarrow
  \cF_2\cF_2^* + \sigma^{-1}\Sigma_{f_2} + \left(
                                               \begin{array}{cc}
                                                 \widehat{\cT}_{f_1} &  \\
                                                  & \cT_{\theta_1} \\
                                               \end{array}
                                             \right) \succ 0.
\end{eqnarray*}
Note that
\begin{eqnarray*}
  \cF_2\cF_2^* + \sigma^{-1}\Sigma_{f_2} + \left(
                                               \begin{array}{cc}
                                                 \widehat{\cT}_{f_1} &  \\
                                                  & \cT_{\theta_1} \\
                                               \end{array}
                                             \right) = \left(
                                                         \begin{array}{cc}
                                                           \cF_{1}\cF_{1}^*+\sigma^{-1}
                                                           \Sigma_{f_{1}} + \widehat{\cT}_{f_{1}} & \cF_{1}\cA_{1}^* \\
                                                           \cA_{1}\cF^*_{1} & \cA_{1}\cA_{1}^* + \sigma^{-1}\cP_{1} + \cT_{\theta_{1}} \\
                                                         \end{array}
                                                       \right).
\end{eqnarray*}
 Since $\cE_{\theta_{1}} = \cA_{1}\cA_{1}^* + \sigma^{-1}\cP_{1} + \cT_{\theta_{1}} \succ 0,$ again  by the Schur complement condition for ensuring  the  positive definiteness of linear operators, we have
\begin{eqnarray*}
  &\left(
  \begin{array}{cc}
  \cF_{1}\cF_{1}^*+\sigma^{-1}
  \Sigma_{f_{1}} + \widehat{\cT}_{f_{1}} & \cF_{1}\cA_{1}^* \\
  \cA_{1}\cF^*_{1} & \cE_{\theta_{1}} \\
  \end{array}
  \right) \succ 0 \\ &\Updownarrow \\
  & \cF_{1}\cF_{1}^* + \sigma^{-1}\Sigma_{f_{1}} + \widehat{\cT}_{f_{1}} - \cF_{1}\cA_{1}^*\cE_{\theta_{1}}^{-1}\cA_{1}\cF^*_{1} \succ 0 \\  &\Updownarrow \\
  &\cF\cF^* + \sigma^{-1}\Sigma_{f} + \cT_f\succ 0.
\end{eqnarray*}
Thus, we have
\begin{eqnarray*}
  \cF_{p+1}\cF_{p+1}^* + \sigma^{-1}\Sigma_{f_{p+1}} + \left(
                                               \begin{array}{cc}
                                                 \cT_{f_p} &  \\
                                                  & \cT_{\theta_p} \\
                                               \end{array}
                                             \right) \succ 0 \Leftrightarrow \cF\cF^* + \sigma^{-1}\Sigma_f + \cT_f \succ 0.
\end{eqnarray*}
The proof of this proposition is completed.
\qed

Note that in the context of the multi-block convex optimization problem \eqref{eq-dd-M},
Assumption \ref{assumption:CQ2} takes the following form:
\begin{assumption}\label{assumption:CQM}
 There exists $(\hat u,\hat y, \hat v, \hat z)\in {\rm ri}({\rm dom}\,f) \times \cY \times{\rm ri}({\rm dom}\,g)\times\cZ $ such that $\cF^* \hat u + \cA^*\hat y + \cG^* \hat v + \cB^*\hat z=c$.
\end{assumption}

After all these preparations, we can finally state our main convergence theorem.

%\begin{theorem}\label{main:themscbpadmm} The sequence $\{(u^k,y^k,v^k,z^k, x^k)\}$ generated by   Algorithm SCB-SPADMM for solving
%(\ref{eq-dd-M}) converges to a unique limit solution point under standard assumptions.
%\end{theorem}

\begin{theorem}\label{main:themscbpadmm} Let
 $\Sigma_f$ and $\Sigma_g $    be   the   two self-adjoint and positive semidefinite operators defined by (\ref{monosub1}) and (\ref{monosub2}), respectively.
Suppose  that the solution set of problem \eqref{eq-dd-M} is nonempty and that  Assumption \ref{assumption:CQM}  holds.   Assume that
 % both $\Sigma_f + T_f +\lambda A^*A $ and $\Sigma_g +T_g+ \lambda  B^*B $ are  positive definite.
$\cT_f$ and $\cT_g$ are chosen such that the
sequence $\{(u^k,y^k,v^k,z^k,x^k)\}$  generated by Algorithm SCB-SPADMM is well defined.
Recall that $\cT_{\theta_i}$ is defined in (\ref{Ttai}) for $1\leq i \leq p$ and $\cT_{\varphi_j}$ is defined in (\ref{Tvarphi}) for $1\leq j \leq q$. %Let $(\bar u, \bar y, \bar v, \bar z, \bar x)$ be any  optimal solution to problem \eqref{eq-dd-M} and $\bar{x}$ be any optimal solution to problem \eqref{eq-pp-M}, respectively.
Then, under the condition either
(a)  $\tau\in (0, (1+\sqrt{5}\,)/2)$ or (b)  $\tau \ge (1+\sqrt{5}\, )/2$ but $\sum_{k=0}^\infty (\|\cG^*(v^{k+1}-v^k) + \cB^*(z^{k+1} - z^k) \|^2 + \tau^{-1} \|\cF^* u^{k+1} + \cA^*y^{k+1} +\cG^* v^{k+1} + \cB^* z^{k+1} -c\|^2) <\infty$,
 the following results hold:
\begin{enumerate}
%\item [{\rm (i)}] The sequence $\{\|x^{k+1}_e\|^2+\|v^{k+1}_e\|^2 _{(\sigma^{-1} \Sigma_g + \cT_g + \cG \cG^*)}+\|u^{k+1}_e\|^2 _{(\sigma^{-1}\Sigma_f +\cT_f +  \cF\cF^*)}\}$  is bounded.
%and
%    \[
%    \lim_{k\rightarrow \infty}\left( \|\lambda^{k+1}-\lambda^k \| + \|y^{k+1}-y^k \|_{\cT_g}+ \|x^{k+1}-x^k \|_{\cT_f}+ \|\cN^* (y^{k+1}- y^{k})\|\right)=0.
%    \]
 \item [{\rm (i)}]   If $(u^\infty, y^\infty, v^\infty, z^\infty, x^\infty)$ is an accumulation point of $\{(u^k, y^k, v^k, z^k, x^k)\}$, then $(u^\infty, y^\infty, v^\infty, z^\infty)$ solves problem \eqref{eq-dd-M} and $x^\infty$  solves \eqref{eq-pp-M}, respectively.
     %, and it holds that
%     \[
%     \lim_{k\rightarrow \infty}\left(  \|x^{k+1}_e\|^2+\|v^{k+1}_e\|^2 _{(\sigma^{-1}\Sigma_g +\cT_g + \cG\cG^*)}+\|u^{k+1}_e\|^2 _{(\sigma^{-1}\Sigma_f +\cT_f +  \cF \cF^*)}\right)=0,
%     \]
%     where in the definition of  $(u_e^k, v^k_e, x_e^k)$, the point $(\bar u, \bar v, \bar x)$ is replaced by $(u^\infty, v^\infty, x^\infty)$.

   \item [{\rm (ii)}] If both $\sigma^{-1}\Sigma_f +\cT_f +  \cF\cF^*$ and $\sigma^{-1}\Sigma_g + \cT_g +\cG \cG^*$ are positive definite, then  the  sequence
$ \{(u^k,y^k,v^k,z^k,x^k)\}$, which  is automatically well defined, converges to a unique limit, say, $(u^\infty, y^\infty, v^\infty, z^\infty, x^\infty)$  with $(u^\infty, y^\infty, v^\infty, z^\infty)$ solving problem \eqref{eq-dd-M} and $x^\infty$ solving  \eqref{eq-pp-M}, respectively.
 \item [{\rm (iii)}]  When the $u,y$-part disappears, the corresponding results in parts (i)--(ii) hold under the condition either $\tau \in (0,2)$ or $\tau \ge2$ but $\sum_{k=0}^\infty \|\cG^* v^{k+1} + \cB^* z^{k+1} -c\|^2<\infty$.
\end{enumerate}
\end{theorem}
{\bf Proof.} By combining Theorem \ref{thempadm} with Proposition \ref{prop:equi-scb-padmm} and Proposition \ref{eqvi-psd},
we can readily obtain the conclusions of this theorem. \qed

\begin{remark}\label{rem:scb-psadmm}
Our SCB-SPADMM algorithm actually provides a potentially efficient approach to handle
large-scale and dense linear constraints. When dealing with such difficult linear systems,
instead of being trapped with the possible convergence issues brought about by inexact solvers such as conjugate gradient methods, one can always first decompose the large systems into serval smaller pieces,
 and then apply our SCB-SPADMM algorithm to the decomposed problems. As a result, these smaller systems can always be handled by adding suitable proximal terms or by solving them exactly.
\end{remark}
 \section{Numerical experiments}\label{section:numerical results}
 We first examine the optimality condition for the general problem \eqref{eq-dd-M} and its dual \eqref{eq-pp-Mori}. Suppose  that the solution set of problem \eqref{eq-dd-M} is nonempty and that  Assumption \ref{assumption:CQM}  holds. Then in order that $(u^*,y^*,v^*,z^*)$ be an optimal solution for \eqref{eq-dd-M} and $x^*$ be an optimal solution for \eqref{eq-pp-Mori}, it is necessary and sufficient that $(u^*,y^*,v^*,z^*)$ and $x^*$ satisfy
 \begin{eqnarray}
 \left\{ \begin{array}{lll}
     \cF^*u + \sum_{i=1}^p \cA_i^*y_i  + \cG^* v + \sum_{j=1}^q \cB_j^* z_j  = c, \\
     f(u) + f^*(-\cF x) = \inprod{-\cF x}{u}, \quad \theta_i(y_i) + \theta^*_i(-\cA_i x) = \inprod{-\cA_i x}{y_i}, \quad i=1,\ldots,p, \\
     g(v) + g^*(-\cG x) = \inprod{-\cG x}{v}, \quad \varphi_i(z_i) + \varphi^*_i(-\cB_i x) = \inprod{-\cB_i x}{z_i}, \quad j=1,\ldots,q.
 \end{array}\right.
 \label{eq:op-M}
 \end{eqnarray}
We will measure the accuracy of an approximate solution based on the above optimality condition. If the given problem is properly scaled, the following relative residual is a natural choice to be used in our stopping criterion:
\begin{eqnarray}
  \eta = \max\{\eta_P, \eta_f, \eta_g, \eta_\theta, \eta_\varphi\},
  \label{eta-stoptol}
\end{eqnarray}
where
\begin{eqnarray*}
&& \eta_P = \frac{\norm{\cF^*u + \cA^*y  + \cG^* v + \cB^* z - c}}{1 + \norm{c}}, \quad
\eta_f = \frac{\norm{u - \textup{Prox}_f(u - \cF x)}}{1 + \norm{u}+\norm{\cF x}}, \quad \eta_g = \frac{\norm{v - \textup{Prox}_g(v - \cG x)}}{1 + \norm{v} + \norm{\cG x}}, \\[5pt]
&& \eta_\theta = \max_{i=1,\ldots,p} \frac{\norm{y_i - \textup{Prox}_{\theta_i}(y_i - \cA_i x)}}{1 + \norm{y_i} + \norm{\cA_i x}}, \quad \eta_\varphi = \max_{j = 1,\ldots,q} \frac{\norm{z_j - \textup{Prox}_{\varphi_j}(z_j - \cB_j x)}}{1 + \norm{z_j} + \norm{\cB_j x}}.
\end{eqnarray*}
Additionally, we compute the relative gap by
%\[\eta_g = \frac{\big(f(u)+\sum_{i=1}^p \theta_i(y_i) + g(v) + \sum_{j=1}^q \varphi_j(z_j)\big) - \big(\inprod{c}{x} + f^*(s)+\sum_{i=1}^p \theta^*_i(r_i) + g^*(t) + \sum_{j=1}^q \varphi^*_j(w_j) \big)}{1+|f(u)+\sum_{i=1}^p \theta_i(y_i) + g(v) + \sum_{j=1}^q \varphi_j(z_j)| + |\inprod{c}{x} + f^*(s)+\sum_{i=1}^p \theta^*_i(r_i) + g^*(t) + \sum_{j=1}^q \varphi^*_j(w_j)|}.\]
\[\eta_{\textup{gap}} = \frac{\textup{obj}_P-\textup{obj}_D}{1+|\textup{obj}_P|+|\textup{obj}_D|},\]
where $\textup{obj}_P := f(u)+\sum_{i=1}^p \theta_i(y_i) + g(v) + \sum_{j=1}^q \varphi_j(z_j)$ and
$\textup{obj}_D := \inprod{c}{x} + f^*(s)+\sum_{i=1}^p \theta^*_i(r_i) + g^*(t) + \sum_{j=1}^q \varphi^*_j(w_j)$. We test the following problem sets.
%Let $\epsilon>0$ be a given accuracy tolerance. We terminate all three algorithms SCB-PADMM, ADMM and ADMMGB when $\eta < \epsilon.$

%%%%%%%%%%%%%%%%%%%%%%%%%%
\subsection{Numerical results for convex quadratic SDP}\label{section:sqsdp}
 Consider the following QSDP problem
 \begin{eqnarray}
  \begin{array}{ll}
    \min &  \frac{1}{2} \inprod{X}{\cQ X} + \inprod{C}{X}  \\[5pt]
   \mbox{s.t.}
       &\cA_E X   =  b_E,\quad \cA_I X \ge b_I,
       \quad
       X \in \Sn_+\cap \cK
\end{array}
 \label{eq-qsdp-eg}
\end{eqnarray}
and its dual problem
 \begin{eqnarray}
  \begin{array}{rllll}
    \max &   -\delta_{\cK}^*(-Z)  + \inprod{b_I}{y_I} - \frac{1}{2}\inprod{\Xprime}{\cQ \Xprime} + \inprod{b_E}{y_E}    \\[5pt]
   \mbox{s.t.} & Z + \cA_I^*y_I - \cQ \Xprime + S + \cA_E^* y_E  = C, \quad  y_I\geq 0, \quad S \in \Sn_+\,.
   \end{array}
   \label{eq-d-qsdp-eg}
\end{eqnarray}
We use $\Xprime$ here to indicate the fact that $\Xprime$ can be different from the primal variable
$X$. Despite this fact, we have that at the optimal point, $\cQ X = \cQ\Xprime$. Since $\cQ$  is only assumed to be a self-adjoint positive semidefinite linear operator, the augmented Lagrangian function associated with \eqref{eq-d-qsdp-eg} may not be strongly convex with respect to $\Xprime$. Without further adding a proximal term, we propose the following strategy to rectify this difficulty.
Since $\cQ$ is positive semidefinite, $\cQ$ can be decomposed as $\cQ = \cB^*\cB$ for some linear map $\cB$. By introducing a new variable $\Xi = -\cB\Xprime$, the problem \eqref{eq-d-qsdp-eg} can be rewritten as follows:
 \begin{eqnarray}
  \begin{array}{rllll}
    \max &  -\delta_{\cK}^*(-Z) + \inprod{b_I}{y_I}-\frac{1}{2}\norm{\Xi}^2_F  + \inprod{b_E}{y_E}    \\[5pt]
   \mbox{s.t.} &  Z + \cA_I^*y_I + \cB^*\Xi + S + \cA_E^* y_E  = C, \quad y_I \ge 0,  \quad
   S \in \Sn_+\, .
   \end{array}
   \label{eq-d1-qsdp-eg}
\end{eqnarray}
Note that now the augmented Lagrangian function associated with \eqref{eq-d1-qsdp-eg} is strongly convex with respect to $\Xi$. Surprisingly, much to our delight,  we can update the iterations in our SCB-SPADMM without explicitly computing $\cB$ or $\cB^*$.
 Given $\overline{Z}, \ybar_I, \overline{S}, \bar{y}_E$ and $\overline{X}$,  denote
\begin{eqnarray*}
  \Xi^{+} := \textup{argmin}_{\Xi}\; \frac{1}{2}\norm{\Xi}^2 +\frac{\sigma}{2}\norm{\overline{Z} +\cA_I^*\ybar_I + \cB^*\Xi + \overline{S} + \cA_E^*\ybar_E - C + \sigma^{-1}\overline{X}}^2
  \;=\; -(\cI + \sigma\cB\cB^*)^{-1}\cB\overline{R},
\end{eqnarray*}
where $\overline{R} = \overline{X} + \sigma(\overline{Z} + \cA_I ^*\ybar_I + \overline{S} + \cA_E^*\ybar_E - C)$.
In updating the SCB-SPADMM iterations, we actually do not need
$\Xi^+$ explicitly, but only need $\Upsilon^{+} := -\cB^*\Xi^{+}$.
From the condition that $(\cI + \sig\cB\cB^*)(-\Xi^+) = \cB \overline{R}$,
we get $(\cI + \sig\cB^*\cB)(-\cB^*\Xi^+) = \cB^*\cB \overline{R}$,
hence we can compute $\Upsilon^+$ via $\cQ$:
 \begin{eqnarray*}
   \Upsilon^{+} %= \cQ\overline{R} -\sig\cQ (\cI + \sigma\cQ)^{-1}\cQ \overline{R}
= (\cI + \sigma\cQ)^{-1} (\cQ\overline{R}).
 \end{eqnarray*}
 In fact, $\Upsilon :=-\cB^*\Xi$ can be viewed as the shadow of $\cQ\Xprime$.
Meanwhile, for the function $\delta_{\cK}^*(-Z)$,
we have the following useful observation that for any $\lambda >0$,
\begin{eqnarray}
 Z^+ =\mbox{argmin}\; \delta_{\cK}^*(-Z) + \frac{\lambda}{2} \norm{Z-\overline{Z}}^2
     =\overline{Z} + \frac{1}{\lambda}\Pi_{\cK}(-\lambda \overline{Z}),
  \label{prox-gz}
\end{eqnarray}
where  (\ref{prox-gz}) follows from the following  Moreau decomposition:
\begin{eqnarray*}
  x = \mbox{Prox}_{\tau f^*}(x) + \tau\mbox{Prox}_{f/\tau}(x/\tau), \quad \forall\,  \tau >0.
\end{eqnarray*}

In our numerical experiments, we test QSDP problems without inequality constraints (i.e., $\cA_I$ and $b_I$ are vacuous). We consider first the linear operator $\cQ$  given by $\cQ(X) = \frac{1}{2}(BX + XB)$ for a given
matrix $B \in \Sn_+.$  Suppose that we have the eigenvalue decomposition $B = P\Lambda P^{T},$ where $\Lambda = \rm{diag}(\lambda)$ and $\lambda = (\lambda_1,\ldots,\lambda_n)^{T}$ is the vector of eigenvalues of B. Then
\begin{eqnarray*}
  \inprod{X}{\cQ X} = \frac{1}{2} \inprod{\widehat{X}}{\Lambda\widehat{X} + \widehat{X}\Lambda}
  = \frac{1}{2}\sum_{i=1}^n\sum_{j=1}^n \widehat{X}_{ij}^2 (\lambda_i + \lambda_j)
  = \sum_{i=1}^n\sum_{j=1}^n \widehat{X}_{ij}^2 H_{ij}^2
  = \inprod{X}{\cB^*\cB X},
\end{eqnarray*}
where $\widehat{X} = P^T X P$, $H_{ij} = \sqrt{\frac{\lambda_i  + \lambda_j}{2}}$,
$\cB X = H\circ(P^T X P)$ and $\cB^* \Xi = P(H\circ \Xi)P^T$.
In our numerical experiments, the matrix $B$ is a low rank random symmetric positive semidefinite matrix. Note that when $\textup{rank}(B) = 0$ and $\cK$ is a polyhedral cone, problem \eqref{eq-qsdp-eg} reduces to the SDP problem considered in \cite{styang2014}. In our experiments, we test both the cases where $\textup{rank}(B) = 5$ and $\textup{rank}(B) = 10$. All the linear constraints are extracted from the numerical test examples in \cite{styang2014} (Section 4.1). For instance, we construct QSDP-BIQ problem sets based on the formulation in \cite{styang2014} as follows:
\begin{eqnarray*}
     \begin{array}{lll}
    \min &   \frac{1}{2}\inprod{X}{\cQ X} + \frac{1}{2} \inprod{Q}{X_0} + \inprod{c}{x}  \\[5pt]
   \mbox{s.t.}
       & \textup{diag}(X_0) - x = 0, \quad \alpha = 1, \\[5pt]
       & X = \left(
               \begin{array}{cc}
                 X_0 & x \\
                 x^T & \alpha \\
               \end{array}
             \right)
       \in \Sn_+, \quad X\in \cK := \{X\in \Sn:\; X\ge 0\}.
       \end{array}
\end{eqnarray*}
In our numerical experiments, the test data for $Q$ and $c$ are taken from Biq Mac Library maintained by Wiegele, which is available at \url{http://biqmac.uni-klu.ac.at/biqmaclib.html}. In the same sprit, we construct test problems  QSDP-BIQ, QSDP-$\theta_+$, QSDP-QAP and QSDP-RCP.

%%%%%%%%%%%%%%%%%%%%%%%%%%%%
%\subsubsection{Numerical results for QSDP}
Here we compare our algorithm \SCBADMM\ with the directly extended {\sc Admm} (with step length $\tau = 1$) and the convergent alternating direction method with a Gaussian back substitution proposed in \cite{he2012alternating} (we call the method {\sc Admmgb} here and use
the parameter $\alpha = 0.99$ in the Gaussian back substitution step). We have implemented all the algorithms \SCBADMM,  {\sc Admm} and \ADMMGB\ in {\sc Matlab} version 7.13.
The numerical results reported later are obtained from  a PC with 24 GB memory and
2.80GHz quad-core CPU running on 64-bit Windows Operating System.

We measure the accuracy of an approximate optimal solution $(X,Z,\Xi,S,y_E)$ for QSDP \eqref{eq-qsdp-eg} and its dual \eqref{eq-d1-qsdp-eg} by using the following relative residual obtained from the general optimality condition \eqref{eq:op-M}:
\begin{eqnarray}
  \eta_{\textup{qsdp}} = \max\{\eta_P, \eta_D, \eta_Z, \eta_{S_1}, \eta_{S_2}\},
  \label{stop:sqsdp}
\end{eqnarray}
where
{\small
\begin{eqnarray*}
&& \eta_P = \frac{\norm{\cA_E X - b_E}}{1+\norm{b_E}},\quad \eta_D = \frac{\norm{Z + \cB^*\Xi  + S + \cA_E^* y_E - C}}{1 + \norm{C}}, \quad
\eta_{Z} = \frac{\norm{X - \Pi_{\cK}(X-Z)}}{1+\norm{X}+\norm{Z}}, \\[5pt]
&& \eta_{S_1} = \frac{|\inprod{S}{X}|}{1+\norm{S}+\norm{X}}, \quad \eta_{S_2} = \frac{\norm{X - \Pi_{\Sn_+}(X)}}{1+\norm{X}}.
\end{eqnarray*}}
%\blue{and for given $Z\in \{Z\in \Sn:\; \delta_{\cK}^*(-Z) < +\infty\}$, $W$ is chosen such that $\delta_{\cK}^*(-Z) = -\inprod{W}{Z}$.}
We terminate the solvers \SCBADMM, {\sc Admm} and \ADMMGB\ when $\eta_{\textup{qsdp}} < 10^{-6}$  with the maximum number of iterations set at 25000.

%%\eta_{C_{1}} = \frac{|\inprod{Z}{W-X}|}{1+\norm{Z}+\norm{W-X}},
%%\eta_{\cK_1} = \frac{\norm{ X - \Pi_{\cK}(X)}}{1 + \norm{X}},

Table \ref{table:sqsdp} reports detailed numerical results for \SCBADMM, {\sc Admm}\ and \ADMMGB\
 in solving some large scale QSDP problems. Here, we only list the results for the case of $\textup{rank}(B) = 10$, since we obtain similar results for the case of $\textup{rank}(B) = 5$. From the numerical results, one can observe that \SCBADMM\ is generally the fastest in terms of the computing time, especially when the problem size is large. In addition, we can see that \SCBADMM\ and {\sc Admm}\ solved all instances to the required accuracy, while \ADMMGB\ failed in certain cases.

\begin{figure}
\centering
\includegraphics[width=0.45\textwidth]{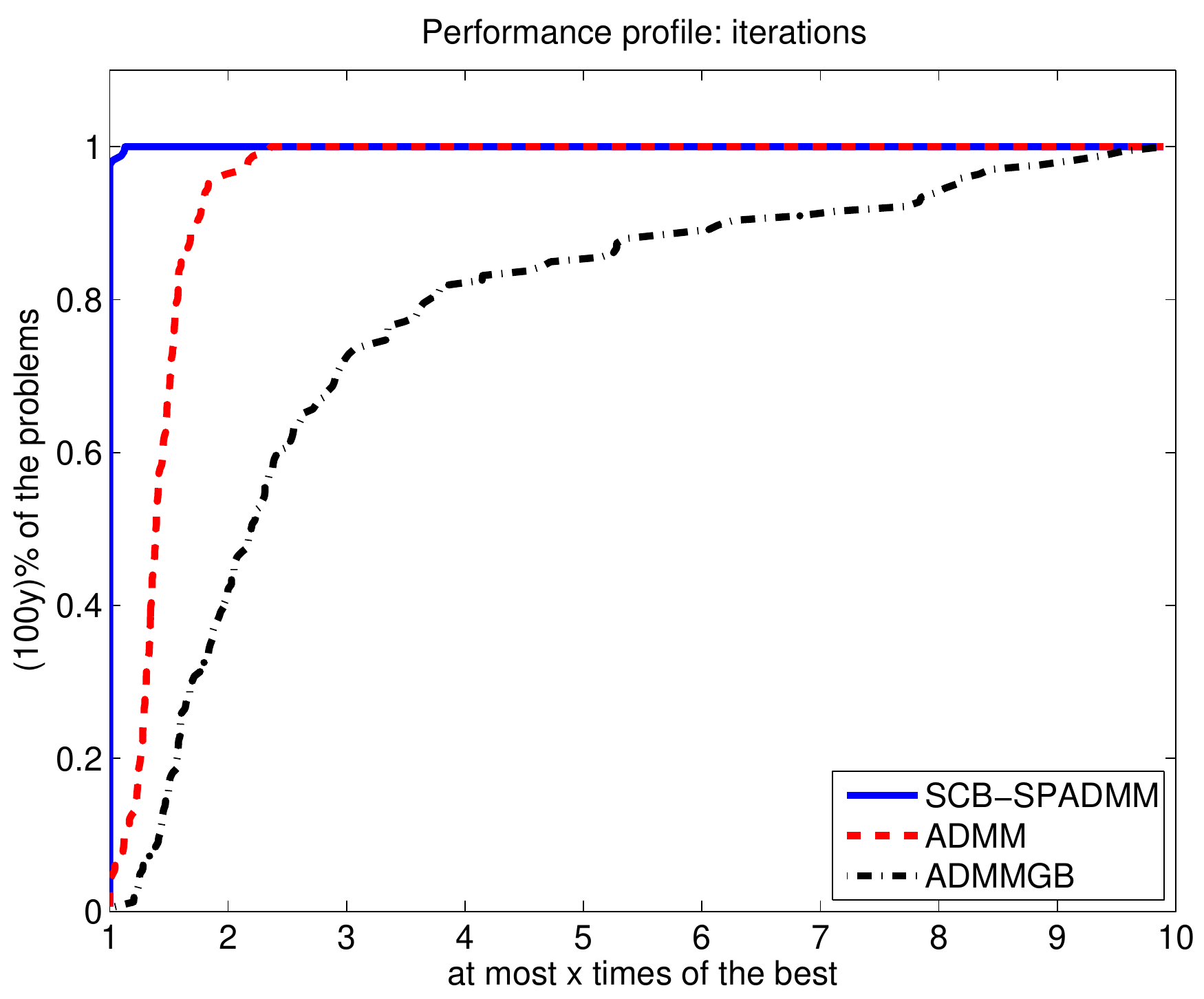}
\quad
\includegraphics[width=0.45\textwidth]{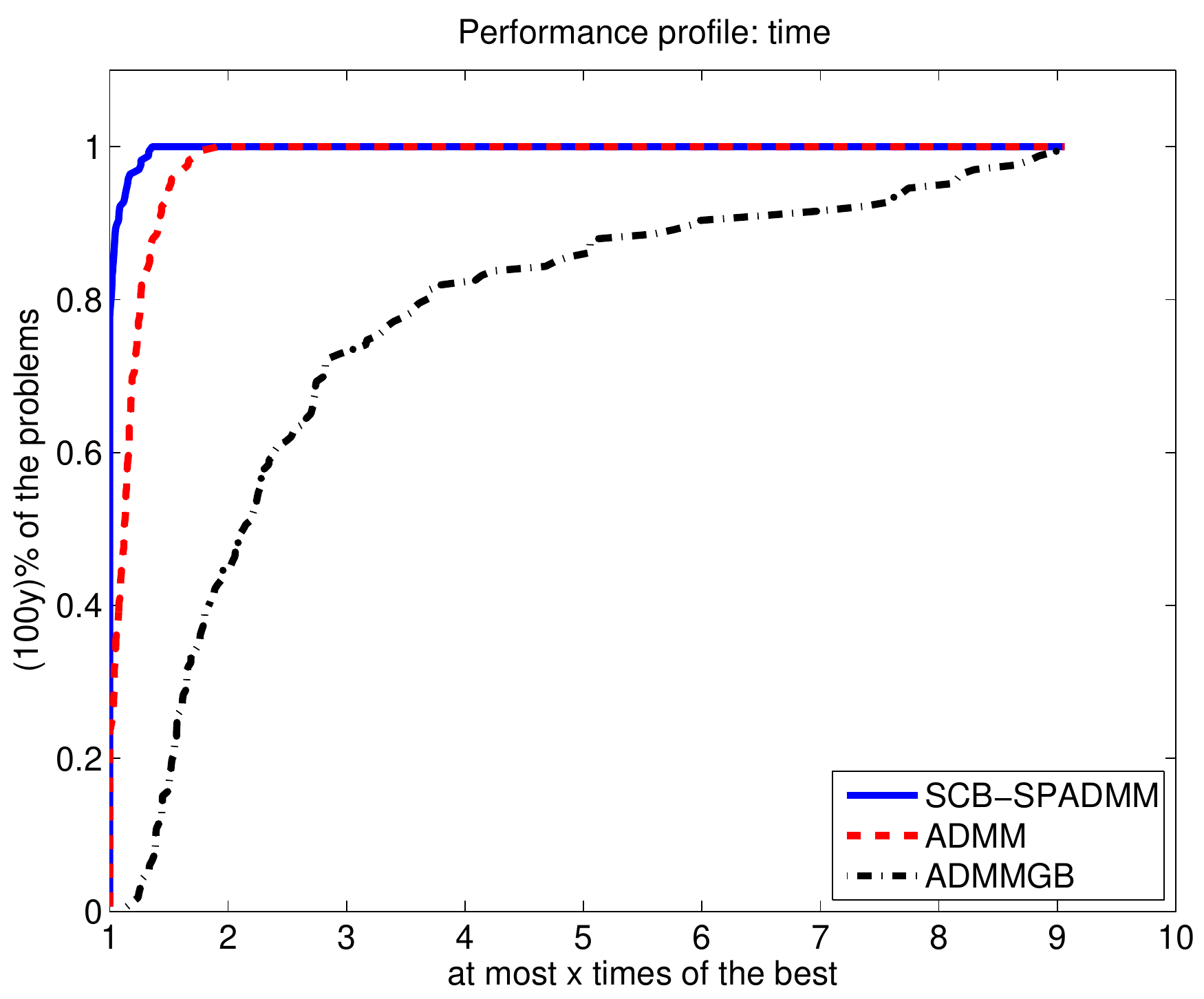}
\caption{Performance profiles of \SCBADMM, {\sc Admm}\ and \ADMMGB\ for the tested large scale QSDP.}
\label{figure-1}
\end{figure}

Figure \ref{figure-1} shows the performance profiles in terms of the number of iterations and computing time for \SCBADMM, {\sc Admm}\ and \ADMMGB, for all the tested large scale QSDP problems. We recall
that a point $(x,y)$ is in the performance profiles curve of a method if and only if it can solve $(100y)\%$ of all the tested problems no slower than $x$ times of any other methods.
We may observe that for the majority of the tested problems, \SCBADMM\ takes the least number of iterations. Besides, in terms of computing time, it can be seen that both \SCBADMM\ and {\sc Admm} outperform \ADMMGB\ by a significant margin, even though {\sc Admm} has no
convergence guarantee.

%%%%%%%%%%%%%%%%%%%%%%%%%%%%%%%
\subsection{ Numerical results for nearest correlation matrix (NCM) approximations}

In this subsection, we first consider the problem  of finding the nearest correlation matrix (NCM) to a given
matrix $G \in \Sn$:
\begin{eqnarray}
    \begin{array}{ll}
    \min  & \frac{1}{2}\norm{H\circ(X-G)}^2_F + \inprod{C}{X}  \\[5pt]
    \mbox{s.t.}    & \cA_E X \;=\; b_E, \quad
       X \in \Sn_+\cap \cK  \, ,
\end{array} \label{eq-egHNCM-F}
\end{eqnarray}
where $H\in \Sn$ is a nonnegative weight matrix, $\cA_E:\Sn \rightarrow \Re^{m_E}$ is a linear map, $G \in \Sn$, $C\in \Sn$ and $b_E\in \Re^{m_E}$ are given data, $\cK$ is a nonempty simple closed  convex set, e.g., $\cK =\{W\in\Sn:\; L\leq W\leq U\}$ with $L,U\in \Sn$ being given matrices. In fact, this is also an instance of the general model of problem \eqref{eq-qsdp-eg} with no inequality constraints, $\cQ X = H\circ H\circ X$ and $\cB X = H\circ X$. We place this special example of QSDP here since an extension will be considered next.
%By introducing slack variables $Y$ and $W$, we can reformulate problem \eqref{eq-egHNCM-F} as
%\begin{eqnarray}
%   \begin{array}{lll}
%    \min &   \frac{1}{2}\norm{Y}^2_F +\inprod{C}{X} + \delta_{\cK}(W)  \\[5pt]
%   \mbox{s.t.}
%       & H\circ(X-G) = Y, \quad \cA_E X =  b_E, \quad X = W, \quad X \in \Sn_+.
%       \end{array}
% \label{eq-p1-egHNCM-F}
%\end{eqnarray}
%The dual of problem  (\ref{eq-p1-egHNCM-F}) is given by
% \begin{eqnarray}
%  \begin{array}{rllll}
%    \max &  -\delta_{\cK}^*(-Z) -\frac{1}{2}\norm{\Xi}^2_F + \inprod{H\circ G}{\Xi} + \inprod{b_E}{y_E}    \\[5pt]
%   \mbox{s.t.} &  Z  + H\circ\Xi + S + \cA_E^* y_E  = C,  \quad
%   S \in \Sn_+.
%   \end{array}
%   \label{eq-d1-egHNCM-F}
%\end{eqnarray}

Now, let's consider an interesting variant of the above NCM problem:
\begin{eqnarray}
    \begin{array}{ll}
    \min  & \norm{H\circ(X-G)}_2 + \inprod{C}{X}  \\[5pt]
    \mbox{s.t.}    & \cA_E X \;=\; b_E,  \quad
     X \in \Sn_+\cap \cK \, .
\end{array} \label{eq-egHNCM}
\end{eqnarray}
Note, in \eqref{eq-egHNCM}, instead of the Frobenius norm, we use the spectral norm. By introducing a slack variable $Y$, we can reformulate problem \eqref{eq-egHNCM} as
\begin{eqnarray}
   \begin{array}{lll}
    \min &   \norm{Y}_2 +\inprod{C}{X}  \\[5pt]
   \mbox{s.t.}
       & H\circ(X-G) = Y, \quad \cA_E X =  b_E,  \quad X \in \Sn_+\cap \cK\, .
       \end{array}
 \label{eq-p1-egHNCM}
\end{eqnarray}
The dual of problem  (\ref{eq-p1-egHNCM}) is given by
 \begin{eqnarray}
  \begin{array}{rllll}
    \max &  -\delta_{\cK}^*(-Z) + \inprod{H\circ G}{\Xi} + \inprod{b_E}{y_E}    \\[5pt]
   \mbox{s.t.} &  Z  + H\circ\Xi + S + \cA_E^* y_E  = C, \quad \norm{\Xi}_*\le 1, \quad
   S \in \Sn_+\, ,
   \end{array}
   \label{eq-d1-egHNCM}
\end{eqnarray}
which is obviously equivalent to the following problem
 \begin{eqnarray}
  \begin{array}{rllll}
    \max &  -\delta_{\cK}^*(-Z) + \inprod{H\circ G}{\Xi}  + \inprod{b_E}{y_E}    \\[5pt]
   \mbox{s.t.} &  Z  + H\circ\Xi + S + \cA_E^* y_E  = C, \quad \norm{\Gamma}_*\le 1, \quad
   S \in \Sn_+\, , \\[5pt]
   & \cD^* \Gamma - \cD^* \Xi = 0,
   \end{array}
   \label{eq-d2-egHNCM}
\end{eqnarray}
where $\cD: \Sn \rightarrow \Sn$ is a nonsingular linear operator. Note that \SCBADMM\
can not be directly applied to solve the problem \eqref{eq-d1-egHNCM} while the equivalent
reformulation \eqref{eq-d2-egHNCM} fits our model nicely.

In our numerical test, matrix $\widehat{G}$ is the gene correlation matrix from \cite{li2010inexact}.
For testing purpose we perturb $\widehat{G}$ to
\begin{eqnarray*}
  G := (1 - \alpha)\widehat{G} + \alpha E,
\end{eqnarray*}
where $\alpha \in (0,1)$ and $E$ is a randomly generated symmetric matrix with entries in $[-1,1]$. We also set $G_{ii} = 1,\ i=1,\ldots,n.$  The weight matrix $H$ is generated
from a weight matrix $H_0$ used by a hedge fund company.
The matrix $H_0$ is a $93 \times 93$ symmetric matrix with all positive entries. It has about $24\%$ of the entries equal to
$10^{-5}$ and the rest are distributed in the interval $[2, 1.28\times 10^3].$ It has $28$ eigenvalues in
the interval $[-520, -0.04]$, $11$ eigenvalues in the interval $[-5\times 10^{-13}, 2\times 10^{-13}]$, and the
rest of $54$ eigenvalues in the interval $[10^{-4}, 2\times 10^4]$.
 The {\sc Matlab} code for generating the matrix $H$ is given by
%% \begin{lstlisting}[language=Matlab]
\begin{verbatim}
      tmp = kron(ones(25,25),H0); H = tmp(1:n,1:n); H = (H'+H)/2.
\end{verbatim}
%%\end{lstlisting}
The reason for using such a weight matrix is because the resulting problems
generated are more challenging to solve as opposed to a randomly generated
weight matrix.
Note that the matrices $G$ and $H$ are generated in the same way as in \cite{jiang2012inexact}.
For simplicity, we further set $C = 0$ and $\cK = \{X\in \Sn:\; X\ge -0.5\}$.

Generally speaking, there is no widely accepted stopping criterion for spectral norm H-weighted NCM problem \eqref{eq-p1-egHNCM}. Here, with reference to the general relative residue \eqref{eta-stoptol}, we measure the accuracy of an approximate optimal solution $(X, Z,\Xi,S,y_E)$ for spectral norm H-weighted NCM problem problem \eqref{eq-egHNCM} (equivalently \eqref{eq-p1-egHNCM}) and its dual \eqref{eq-d1-egHNCM} (equivalently \eqref{eq-d2-egHNCM}) by using the following relative residual derived from the general optimality condition \eqref{eq:op-M}:
\begin{eqnarray}
  \eta_{\textup{sncm}} = \max\{\eta_P, \eta_D, \eta_{Z}, \eta_{S_1}, \eta_{S_2}, \eta_{\Xi} \},
  \label{stop:sncm}
\end{eqnarray}
where
{\small
\begin{eqnarray*}
&& \eta_P = \frac{\norm{\cA_E X - b_E}}{1+\norm{b_E}},\quad \eta_D = \frac{\norm{Z + H\circ\Xi  + S + \cA_E^* y_E}}{1 + \norm{Z} + \norm{S}}, \quad
\eta_{Z} = \frac{\norm{X - \Pi_{\cK}(X-Z)}}{1+\norm{X}+\norm{Z}}, \\[5pt]
&& \eta_{S_{1}} = \frac{|\inprod{S}{X}|}{1+\norm{S}+\norm{X}}, \quad \eta_{S_2} = \frac{\norm{X - \Pi_{\Sn_+}(X)}}{1+\norm{X}}, \quad \eta_{\Xi} = \frac{\norm{\Xi - \Pi_{\{X\in\Re^{n\times n} \; :\norm{X}_* \le 1 \}}(\Xi - H\circ(X- G)}}{1+\norm{\Xi} + \norm{H\circ(X-G)}}.
\end{eqnarray*}}
%\blue{and for given $Z\in \{Z\in \Sn:\; \delta_{\cK}^*(-Z) < +\infty\}$, $W$ is chosen such that $\delta_{\cK}^*(-Z) = -\inprod{W}{Z}$. }
%\footnotetext{We already take $C=0$ into account.}
%%%%%%%%%%%%%%%%%%%%%%%%%%%%%%%%%%%%
%\subsubsection{Numerical results for NCM problems}

Firstly, numerical results for solving F-norm H-weighted NCM problems \eqref{eq-egHNCM} are reported. We compare all three algorithms, namely
\SCBADMM, \ADMM, \ADMMGB\ using the relative residue
\eqref{stop:sqsdp}. We terminate the solvers  when $\eta_{\textup{qsdp}} < 10^{-6}$  with the maximum number of iterations set at 25000.

In Table \ref{table:ncmf}, we report detailed numerical results for \SCBADMM, \ADMM\ and \ADMMGB\ in solving various instances of F-norm H-weighted NCM problem. As we can see from Table \ref{table:ncmf}, our \SCBADMM\ is certainly more efficient than the other two algorithms on most of the problems tested.

%\begin{center}
%\begin{figure}
%\includegraphics[width=0.9\textwidth,height=0.5\textwidth]{profile}
%\caption{Performance profiles of \SCBADMM, \ADMM\ and \ADMMGB.}
%\label{figure-1}
%\end{figure}
%\end{center}

%%%%%%%%%%%%%%%%%%%%%%%%%%%%
\begin{table}
\centering
\begin{scriptsize}
\caption{{\small The performance of \SCBADMM, \ADMM, \ADMMGB\ on Frobenius norm H-weighted NCM problems (dual of \eqref{eq-egHNCM-F}) (accuracy $= 10^{-6}$). In the table, ``scb'' stands for \SCBADMM\ and ``gb'' stands for \ADMMGB, respectively. The computation time is in the format of ``hours:minutes:seconds''.}}
\label{table:ncmf}
\begin{tabular}{| ccc | c |  c | c| c|}
\hline
   \mc{3}{|c|}{} &\mc{1}{c|}{} &\mc{1}{c|}{}&\mc{1}{c|}{}&\mc{1}{c|}{}\\[-5pt]
\mc{3}{|c|}{} & \mc{1}{c|}{iteration} &\mc{1}{c|}{$\eta_{\textup{qsdp}}$}&\mc{1}{c|}{$\eta_\textup{gap}$}&\mc{1}{c|}{time}\\[2pt] \hline
\mc{1}{|c}{problem} &\mc{1}{c}{$n_s$} &\mc{1}{c|}{$\alpha$}&\mc{1}{c|}{scb$|$admm$|$gb}&\mc{1}{c|}{scb$|$admm$|$gb}&\mc{1}{c|}{scb$|$admm$|$gb}&\mc{1}{c|}{scb$|$admm$|$gb}\\[2pt]
\hline
Lymph
	 &587 &0.10	 &263  $|$  522  $|$  696 	 &   9.9-7 $|$    9.9-7 $|$    9.9-7	 &   -4.4-7 $|$    -4.5-7 $|$    -4.0-7 	 &30 $|$ 53 $|$ 1:23\\[2pt]
	 &587 &0.05	 &264  $|$  356  $|$  592 	 &   9.9-7 $|$    9.9-7 $|$    9.9-7	 &   -3.9-7 $|$    -3.4-7 $|$    -3.0-7 	 &29 $|$ 35 $|$ 1:08\\[2pt]
 \hline
ER
	 &692 &0.10	 &268  $|$  355  $|$  711 	 &   9.9-7 $|$    9.9-7 $|$    9.9-7	 &   -5.1-7 $|$    -4.7-7 $|$    -4.2-7 	 &43 $|$ 51 $|$ 1:58\\[2pt]
	 &692 &0.05	 &226  $|$  293  $|$  603 	 &   9.9-7 $|$    9.9-7 $|$    9.9-7	 &   -4.2-7 $|$    -3.8-7 $|$    -3.3-7 	 &37 $|$ 43 $|$ 1:54\\[2pt]
 \hline
Arabidopsis
	 &834 &0.10	 &510  $|$  528  $|$  725 	 &   9.9-7 $|$    9.9-7 $|$    9.9-7	 &   -5.9-7 $|$    -5.3-7 $|$    -3.9-7 	 &2:11 $|$ 2:02 $|$ 3:03\\[2pt]
	 &834 &0.05	 &444  $|$  470  $|$  650 	 &   9.9-7 $|$    9.9-7 $|$    9.9-7	 &   -5.8-7 $|$    -5.2-7 $|$    -4.8-7 	 &1:51 $|$ 1:43 $|$ 2:44\\[2pt]
 \hline
Leukemia
	 &1255 &0.10	 &292  $|$  420  $|$  826 	 &   9.9-7 $|$    9.9-7 $|$    9.9-7	 &   -5.4-7 $|$    -5.3-7 $|$    -4.4-7 	 &3:13 $|$ 4:11 $|$ 9:13\\[2pt]
	 &1255 &0.05	 &251  $|$  408  $|$  670 	 &   9.9-7 $|$    9.7-7 $|$    9.6-7	 &   -5.4-7 $|$    -4.9-7 $|$    -4.0-7 	 &2:48 $|$ 4:03 $|$ 7:35\\[2pt]
 \hline
hereditarybc
	 &1869 &0.10	 &555  $|$  634  $|$  871 	 &   9.9-7 $|$    9.9-7 $|$    9.9-7	 &   -9.1-7 $|$    -9.1-7 $|$    -7.0-7 	 &17:39 $|$ 18:38 $|$ 28:01\\[2pt]
	 &1869 &0.05	 &530  $|$  626  $|$  839 	 &   9.9-7 $|$    9.9-7 $|$    9.9-7	 &   -8.7-7 $|$    -8.7-7 $|$    -5.2-7 	 &16:50 $|$ 18:15 $|$ 26:34\\[2pt]
 \hline
\end{tabular}
\end{scriptsize}
\end{table}
%\FloatBarrier

The rest of this subsection is devoted to the numerical results of the spectral norm H-weighted NCM problem \eqref{eq-egHNCM}. As mentioned before, \SCBADMM\ is applied to solve the problem \eqref{eq-d2-egHNCM} rather than \eqref{eq-d1-egHNCM}. We implemented all the algorithms for solving problem \eqref{eq-d2-egHNCM} using the relative residue \eqref{stop:sncm}.  We terminate the solvers  when $\eta_{\textup{sncm}} < 10^{-5}$  with the maximum number of iterations set at 25000. In Table \ref{table:ncms}, we report detailed numerical results for \SCBADMM, \ADMM\ and \ADMMGB\ in solving various instances of spectral norm H-weighted NCM problem. As we can see from Table \ref{table:ncms}, our {\SCBADMM} is much more efficient than the other two algorithms.

Observe that although there is no convergence guarantee, one may still apply the directly extended \ADMM\ with 4 blocks to the original dual problem \eqref{eq-d1-egHNCM} by adding a proximal term for the $\Xi$ part. We call this method {\sc Ladmm}.  Moreover, by using the same proximal strategy for $\Xi$, a convergent linearized alternating direction method with a Gausssian back substitution proposed in \cite{he2011linearized} (we call the method
{\sc Ladmmgb} here and use the parameter $\alpha = 0.99$ in the Gasussian back substitution step) can also be applied to the original problem \eqref{eq-d1-egHNCM}. We have also implemented
{\sc Ladmm} and
{\sc Ladmmgb} in {\sc Matlab}. Our experiments show that solving the problem \eqref{eq-d1-egHNCM} directly is much slower than solving the equivalent problem \eqref{eq-d2-egHNCM}. Thus,
the reformulation of \eqref{eq-d1-egHNCM} to \eqref{eq-d2-egHNCM} is in fact
advantageous for both \ADMM\ and \ADMMGB.
In Table \ref{table:ncmsP}, for the purpose of illustration we  list a couple of  detailed numerical results
on the performance of {\sc Ladmm} and {\sc Ladmmgb}.

%%%%%%%%%%%%%%%%%%%%%%%%%%%%%%%%%%%%%%

\begin{table}[t]
\begin{center}
\begin{scriptsize}
\caption{ {\small The performance of \SCBADMM, \ADMM, \ADMMGB\ on spectral norm H-weighted NCM problem \eqref{eq-d2-egHNCM}
(accuracy $= 10^{-5}$). In the table, ``scb'' stands for \SCBADMM\ and ``gb'' stands for
\ADMMGB, respectively. The computation time is in the format of ``hours:minutes:seconds''.}}
\label{table:ncms}
\begin{tabular}{| ccc | c |  c | c| c|}
\hline
 \mc{3}{|c|}{} &\mc{1}{c|}{} &\mc{1}{c|}{}&\mc{1}{c|}{}&\mc{1}{c|}{}\\[-5pt]
\mc{3}{|c|}{} & \mc{1}{c|}{iteration} &\mc{1}{c|}{$\eta_{\textup{sncm}}$}&\mc{1}{c|}{$\eta_\textup{gap}$}&\mc{1}{c|}{time}\\[2pt] \hline
\mc{1}{|c}{problem} &\mc{1}{c}{$n_s$} &\mc{1}{c|}{$\alpha$}&\mc{1}{c|}{scb$|$admm$|$gb}&\mc{1}{c|}{scb$|$admm$|$gb}&\mc{1}{c|}{scb$|$admm$|$gb}&\mc{1}{c|}{scb$|$admm$|$gb}\\[3pt]
\hline

Lymph
	 &587 &0.10	 &4110  $|$  6048  $|$  7131 	 &   {\blue{ 9.9-6}} $|$    {\blue{ 9.9-6}} $|$    {\red{ 1.0-5}}	 &   {\red{ -3.4-5}} $|$    {\red{ -2.8-5}} $|$    {\red{ -2.7-5}} 	 &13:21 $|$ 17:10 $|$ 21:43\\[2pt]
	 &587 &0.05	 &5001  $|$  7401  $|$  8101 	 &   {\blue{ 9.8-6}} $|$    {\blue{ 9.9-6}} $|$    {\blue{ 9.9-6}}	 &   {\red{ -2.0-5}} $|$    {\red{ -2.3-5}} $|$    {\blue{ -8.1-6}} 	 &19:41 $|$ 21:25 $|$ 25:13\\[2pt]
 \hline
ER
	 &692 &0.10	 &3251  $|$  4844  $|$  6478 	 &   {\blue{ 9.9-6}} $|$    {\blue{ 9.9-6}} $|$    {\red{ 1.0-5}}	 &   {\red{ -3.1-5}} $|$    {\red{ -2.6-5}} $|$    {\blue{ -6.0-6}} 	 &15:06 $|$ 19:30 $|$ 28:03\\[2pt]
	 &692 &0.05	 &4201  $|$  5851  $|$  7548 	 &   {\blue{ 9.3-6}} $|$    {\blue{ 9.8-6}} $|$    {\red{ 1.0-5}}	 &   {\red{ -3.5-5}} $|$    {\red{ -2.9-5}} $|$    {\red{ -3.4-5}} 	 &18:44 $|$ 23:46 $|$ 32:57\\[2pt]
 \hline
Arabidopsis
	 &834 &0.10	 &3344  $|$  6251  $|$  7965 	 &   {\blue{ 9.9-6}} $|$    {\blue{ 9.7-6}} $|$    {\red{ 1.0-5}}	 &   {\red{ -3.8-5}} $|$    {\red{ -2.0-5}} $|$    {\red{ -3.7-5}} 	 &23:20 $|$ 40:12 $|$ 54:31\\[2pt]
	 &834 &0.05	 &2496  $|$  3101  $|$  3231 	 &   {\blue{ 9.9-6}} $|$    {\blue{ 9.9-6}} $|$    {\red{ 1.0-5}}	 &   {\red{ -9.1-5}} $|$    {\red{ -4.3-5}} $|$    {\red{ -5.3-5}} 	 &17:03 $|$ 19:53 $|$ 21:56\\[2pt]
 \hline
Leukemia
	 &1255 &0.10	 &4351  $|$  6102  $|$  7301 	 &   {\blue{ 9.9-6}} $|$    {\blue{ 9.9-6}} $|$    {\red{ 1.0-5}}	 &   {\red{ -3.7-5}} $|$    {\red{ -3.3-5}} $|$    {\red{ -3.0-5}} 	 &  1:22:42 $|$   1:49:02 $|$   2:16:52\\[2pt]
	 &1255 &0.05	 &3957  $|$  5851  $|$  10151 	 &   {\blue{ 9.9-6}} $|$    {\blue{ 9.7-6}} $|$    {\blue{ 9.5-6}}	 &   {\red{ -7.2-5}} $|$    {\red{ -5.7-5}} $|$    {\red{ -1.1-5}} 	 &  1:18:19 $|$   1:44:47 $|$   3:26:08\\[2pt]
 \hline
\end{tabular}
\end{scriptsize}
\end{center}
\end{table}

\begin{table}[t]
\begin{center}
\begin{scriptsize}
\caption{{\small The performance of  {\sc Ladmm}, {\sc Ladmmgb} on spectral norm H-weighted NCM problem\eqref{eq-d1-egHNCM}
(accuracy $= 10^{-5}$). In the table, ``lgb'' stands for {\sc Ladmmgb}. The computation time is in the format of ``hours:minutes:seconds''.}}
\label{table:ncmsP}
\begin{tabular}{| ccc | c |  c | c| c|}
\hline
 \mc{3}{|c|}{} &\mc{1}{c|}{} &\mc{1}{c|}{}&\mc{1}{c|}{}&\mc{1}{c|}{}\\[-5pt]
\mc{3}{|c|}{} & \mc{1}{c|}{iteration} &\mc{1}{c|}{$\eta_{\textup{sncm}}$}&\mc{1}{c|}{$\eta_\textup{gap}$}&\mc{1}{c|}{time}\\[2pt] \hline
\mc{1}{|c}{problem} &\mc{1}{c}{$n_s$} &\mc{1}{c|}{$\alpha$}&\mc{1}{c|}{ladmm$|$lgb}&\mc{1}{c|}{ladmm$|$lgb}&\mc{1}{c|}{ladmm$|$lgb}&\mc{1}{c|}{ladmm$|$lgb}\\ \hline
\mc{1}{|c}{Lymph}
         &587    &0.10   &8401  $|$  25000       &   9.9-6 $|$    {\green{ 1.4-5}}         &   {\green{ -1.6-5}} $|$    {\green{ -2.1-5}}    &23:59 $|$   1:22:58\\[2pt]
\mc{1}{|c}{Lymph}
         &587    &0.05   &  13609 $|$  25000      &   9.9-6 $|$    {\green{ 2.3-5}}       &   {\green{ -1.6-5}} $|$    {\green{ -4.2-5}}          &39:29 $|$   1:18:50\\[2pt]
 \hline
\end{tabular}
\end{scriptsize}
\end{center}
\end{table}

%%%%%%%%%%%%%%%%%%%%%%%%%%%%%%%%%%%%%%

 \section{Conclusions}\label{section:final}
In this paper, we have proposed a Schur complement based convergent yet efficient  semi-proximal ADMM for solving convex programming problems, with a coupling linear equality constraint,
whose objective function is the sum of two proper closed convex functions plus an arbitrary number of convex quadratic or linear functions.  The ability of dealing with  an arbitrary number of convex quadratic or linear functions
 in the objective function makes the proposed algorithm very flexible in   solving various multi-block convex optimization problems.
By conducting numerical experiments on QSDP and its extensions, we have presented convincing numerical results to  demonstrate the superior performance of our
proposed SCB-SPADMM. As mentioned in the introduction, our primary motivation of introducing this SCB-SPADMM is to quickly generate a good initial point
so as to warm-start methods which have fast local convergence properties. For standard linear SDP and linear SDP with doubly nonnegative constraints, this has already been done by Zhao, Sun and Toh in \cite{zhao2010newton} and Yang, Sun and Toh in \cite{ystoh2014}, respectively.  Naturally, our next target is to extend the approach   of \cite{zhao2010newton,ystoh2014} to solve QSDP with an initial point generated by SCB-SPADMM. We will report our corresponding findings in subsequent works.

 \section{Acknowledgements}
The authors would like to thank Mr Liuqin Yang at National University of Singapore for
suggestions on the
numerical implementations of the algorithms described in the paper.

\bibliographystyle{siam}
\bibliography{SCADMMbib}
%\begin{thebibliography}{99}
%
%\bibitem{fazel2013hankel} M. Fazel, T. K. Pong, D. F. Sun, and P. Tseng, {\em Hankel matrix rank minimization with applications to system identification and realization,} SIAM Journal on Matrix Analysis and Applications 34, no. 3 (2013): 946-977.
%
%\bibitem{styang2014} Defeng Sun, Kim-Chuan Toh and Liuqin Yang, {\em A convergent proximal  alternating direction method
%of multipliers for conic programming with $4$-block constraints,} Technical Report, National University of Singapore, April 2014.
%
%\bibitem{wright2009robust} Wright, John, Arvind Ganesh, Shankar Rao, Y. Peng, and Y. Ma, {\em Robust principal component analysis: Exact recovery of corrupted low-rank matrices by convex optimization,} In Proc. of Neural Information Processing Systems, vol. 3. 2009.
%
%\end{thebibliography}
%where $\delta_\cK(\cdot)$ is the indicate function of $\cK$, i.e., $\delta_\cK(x) = 0$ if $x \in \cK$ and $\delta_\cK(x) = \infty$ if $x \notin \cK$.

%{|c|r|r||l||c||@{}c@{}|@{}r@{}|r|}

\begin{landscape}
\begin{scriptsize}
\begin{longtable}{| @{}c@{}c@{}c@{} | c |  c | c| c|}
\caption{The performance of \SCBADMM, \ADMM, \ADMMGB\ on QSDP-$\theta_+$, QSDP-QAP, QSDP-BIQ and QSDP-RCP problems
(accuracy $= 10^{-6}$). In the table, ``scb'' stands for \SCBADMM and ``gb'' stands for \ADMMGB, respectively. The computation time is in the format of ``hours:minutes:seconds''.}\label{table:sqsdp}
\\
\hline
 \mc{3}{|c|}{} &\mc{1}{c|}{} &\mc{1}{c|}{}&\mc{1}{c|}{}&\mc{1}{c|}{}\\[-5pt]
\mc{3}{|c|}{} & \mc{1}{c|}{iteration} &\mc{1}{c|}{$\eta_{\textup{qsdp}}$}&\mc{1}{c|}{$\eta_\textup{gap}$}&\mc{1}{c|}{time}\\[2pt] \hline
\mc{1}{|@{}c@{}}{problem} &\mc{1}{@{}c@{}}{$m_E;n_s$} &\mc{1}{@{}c@{}|}{rank(B)}&\mc{1}{c|}{scb$|$admm$|$gb}&\mc{1}{c|}{scb$|$admm$|$gb}&\mc{1}{c|}{scb$|$admm$|$gb}&\mc{1}{c|}{scb$|$admm$|$gb}\\ \hline
\endhead

\input{table.dat}

\end{longtable}
\end{scriptsize}
\end{landscape}

\end{document}

%% file: table.dat
theta6 
	 &4375  $;$  300 	 &10 	 &311  $|$  407  $|$  549 	 &   7.9-7 $|$    9.7-7 $|$    9.9-7	 &   {\green{ 2.1-6}} $|$    {\green{ -1.6-6}} $|$    -6.2-7 	 &08 $|$ 09 $|$ 14\\[2pt] 
 \hline 
theta62 
	 &13390  $;$  300 	 &10 	 &153  $|$  196  $|$  229 	 &   9.6-7 $|$    9.9-7 $|$    9.6-7	 &   -1.1-7 $|$    9.6-8 $|$    -4.5-7 	 &04 $|$ 05 $|$ 06\\[2pt] 
 \hline 
theta8 
	 &7905  $;$  400 	 &10 	 &314  $|$  384  $|$  616 	 &   9.5-7 $|$    9.6-7 $|$    9.5-7	 &   {\green{ 2.7-6}} $|$    {\green{ -1.3-6}} $|$    -5.4-7 	 &17 $|$ 18 $|$ 33\\[2pt] 
 \hline 
theta82 
	 &23872  $;$  400 	 &10 	 &158  $|$  179  $|$  234 	 &   9.5-7 $|$    9.7-7 $|$    9.9-7	 &   -3.7-8 $|$    -2.8-7 $|$    -8.2-7 	 &10 $|$ 09 $|$ 13\\[2pt] 
 \hline 
theta83 
	 &39862  $;$  400 	 &10 	 &200  $|$  177  $|$  219 	 &   9.3-7 $|$    9.6-7 $|$    9.4-7	 &   6.2-9 $|$    1.4-7 $|$    -1.2-7 	 &11 $|$ 09 $|$ 14\\[2pt] 
 \hline 
theta10 
	 &12470  $;$  500 	 &10 	 &329  $|$  439  $|$  614 	 &   9.0-7 $|$    8.5-7 $|$    9.7-7	 &   {\green{ -2.5-6}} $|$    {\green{ 1.5-6}} $|$    5.8-7 	 &27 $|$ 33 $|$ 50\\[2pt] 
 \hline 
theta102 
	 &37467  $;$  500 	 &10 	 &150  $|$  187  $|$  235 	 &   8.7-7 $|$    9.4-7 $|$    9.9-7	 &   6.4-7 $|$    2.9-7 $|$    -9.3-7 	 &15 $|$ 15 $|$ 21\\[2pt] 
 \hline 
theta103 
	 &62516  $;$  500 	 &10 	 &202  $|$  184  $|$  222 	 &   9.8-7 $|$    9.5-7 $|$    9.9-7	 &   -4.2-8 $|$    6.9-8 $|$    -1.6-7 	 &20 $|$ 15 $|$ 21\\[2pt] 
 \hline 
theta104 
	 &87245  $;$  500 	 &10 	 &181  $|$  181  $|$  242 	 &   9.4-7 $|$    9.5-7 $|$    9.9-7	 &   6.9-8 $|$    2.0-7 $|$    -2.8-7 	 &20 $|$ 15 $|$ 23\\[2pt] 
 \hline 
theta12 
	 &17979  $;$  600 	 &10 	 &343  $|$  441  $|$  703 	 &   9.9-7 $|$    8.3-7 $|$    9.9-7	 &   {\green{ 3.0-6}} $|$    {\green{ 1.4-6}} $|$    -8.8-7 	 &42 $|$ 48 $|$ 1:27\\[2pt] 
 \hline 
theta123 
	 &90020  $;$  600 	 &10 	 &204  $|$  205  $|$  213 	 &   9.7-7 $|$    9.8-7 $|$    9.9-7	 &   -9.1-8 $|$    6.6-8 $|$    -1.9-7 	 &29 $|$ 25 $|$ 31\\[2pt] 
 \hline 
san200-0.7-1 
	 &5971  $;$  200 	 &10 	 &2150  $|$  4758  $|$  5172 	 &   9.8-7 $|$    9.9-7 $|$    9.9-7	 &   {\blue{ 5.1-6}} $|$    {\green{ 2.0-6}} $|$    {\green{ -3.5-6}} 	 &15 $|$ 26 $|$ 36\\[2pt] 
 \hline 
sanr200-0.7 
	 &6033  $;$  200 	 &10 	 &177  $|$  223  $|$  280 	 &   9.6-7 $|$    9.7-7 $|$    9.7-7	 &   1.9-7 $|$    -6.0-8 $|$    1.7-8 	 &02 $|$ 02 $|$ 03\\[2pt] 
 \hline 
c-fat200-1 
	 &18367  $;$  200 	 &10 	 &2257  $|$  3027  $|$  3268 	 &   9.9-7 $|$    9.7-7 $|$    9.9-7	 &   {\green{ -2.6-6}} $|$    {\green{ -2.0-6}} $|$    {\green{ -2.2-6}} 	 &24 $|$ 26 $|$ 35\\[2pt] 
 \hline 
hamming-8-4 
	 &11777  $;$  256 	 &10 	 &2820  $|$  2945  $|$  3517 	 &   9.9-7 $|$    9.9-7 $|$    9.9-7	 &   -6.0-7 $|$    -6.4-7 $|$    {\green{ -1.1-6}} 	 &53 $|$ 49 $|$ 1:09\\[2pt] 
 \hline 
hamming-9-8 
	 &2305  $;$  512 	 &10 	 &3891  $|$  4980  $|$  5577 	 &   9.9-7 $|$    9.9-7 $|$    9.9-7	 &   {\green{ -3.4-6}} $|$    -5.8-7 $|$    9.9-7 	 &3:54 $|$ 4:12 $|$ 5:50\\[2pt] 
 \hline 
hamming-8-3-4 
	 &16129  $;$  256 	 &10 	 &202  $|$  220  $|$  294 	 &   4.8-7 $|$    8.9-7 $|$    9.8-7	 &   {\green{ 4.5-6}} $|$    5.9-7 $|$    2.2-7 	 &04 $|$ 04 $|$ 06\\[2pt] 
 \hline 
hamming-9-5-6 
	 &53761  $;$  512 	 &10 	 &436  $|$  535  $|$  684 	 &   8.5-7 $|$    8.7-7 $|$    9.6-7	 &   {\red{ 1.1-5}} $|$    {\green{ -1.7-6}} $|$    -1.6-7 	 &36 $|$ 37 $|$ 57\\[2pt] 
 \hline 
brock200-1 
	 &5067  $;$  200 	 &10 	 &198  $|$  210  $|$  291 	 &   9.7-7 $|$    9.4-7 $|$    9.8-7	 &   9.9-8 $|$    -2.9-7 $|$    -6.9-10 	 &02 $|$ 02 $|$ 03\\[2pt] 
 \hline 
brock200-4 
	 &6812  $;$  200 	 &10 	 &209  $|$  186  $|$  263 	 &   9.8-7 $|$    9.9-7 $|$    9.8-7	 &   1.2-7 $|$    -2.6-9 $|$    -1.1-7 	 &03 $|$ 02 $|$ 03\\[2pt] 
 \hline 
brock400-1 
	 &20078  $;$  400 	 &10 	 &168  $|$  217  $|$  275 	 &   9.0-7 $|$    9.6-7 $|$    9.7-7	 &   8.6-7 $|$    -4.9-8 $|$    6.2-9 	 &11 $|$ 10 $|$ 15\\[2pt] 
 \hline 
keller4 
	 &5101  $;$  171 	 &10 	 &669  $|$  909  $|$  963 	 &   9.9-7 $|$    9.9-7 $|$    9.9-7	 &   -1.3-8 $|$    4.6-9 $|$    -8.4-8 	 &06 $|$ 07 $|$ 09\\[2pt] 
 \hline 
p-hat300-1 
	 &33918  $;$  300 	 &10 	 &468  $|$  829  $|$  2501 	 &   9.9-7 $|$    9.9-7 $|$    8.3-7	 &   -8.7-7 $|$    2.1-7 $|$    {\green{ -1.0-6}} 	 &14 $|$ 20 $|$ 1:09\\[2pt] 
 \hline 
be250.1 
	 &251  $;$  251 	 &10 	 &4126  $|$  7439  $|$  25000 	 &   9.6-7 $|$    9.9-7 $|$    {\green{ 1.3-6}}	 &   -5.8-7 $|$    -8.6-7 $|$    -1.3-8 	 &59 $|$ 1:27 $|$ 5:41\\[2pt] 
 \hline 
be250.2 
	 &251  $;$  251 	 &10 	 &3604  $|$  6504  $|$  16322 	 &   9.8-7 $|$    9.9-7 $|$    9.9-7	 &   -4.9-7 $|$    -6.8-7 $|$    -7.4-9 	 &52 $|$ 1:18 $|$ 3:40\\[2pt] 
 \hline 
be250.3 
	 &251  $;$  251 	 &10 	 &3562  $|$  5712  $|$  8501 	 &   9.9-7 $|$    9.9-7 $|$    9.7-7	 &   -9.2-7 $|$    -9.4-7 $|$    9.3-7 	 &52 $|$ 1:08 $|$ 1:57\\[2pt] 
 \hline 
be250.4 
	 &251  $;$  251 	 &10 	 &4072  $|$  7668  $|$  25000 	 &   9.9-7 $|$    9.9-7 $|$    {\green{ 1.4-6}}	 &   {\green{ -2.1-6}} $|$    {\green{ 2.8-6}} $|$    -9.4-9 	 &57 $|$ 1:32 $|$ 5:41\\[2pt] 
 \hline 
be250.5 
	 &251  $;$  251 	 &10 	 &3210  $|$  4635  $|$  7406 	 &   9.9-7 $|$    9.9-7 $|$    9.9-7	 &   -8.6-7 $|$    -8.8-7 $|$    {\green{ 1.4-6}} 	 &46 $|$ 55 $|$ 1:41\\[2pt] 
 \hline 
be250.6 
	 &251  $;$  251 	 &10 	 &3250  $|$  5580  $|$  9812 	 &   9.9-7 $|$    9.9-7 $|$    9.9-7	 &   -2.8-7 $|$    -3.1-7 $|$    -3.6-7 	 &46 $|$ 1:05 $|$ 2:10\\[2pt] 
 \hline 
be250.7 
	 &251  $;$  251 	 &10 	 &3699  $|$  6562  $|$  13501 	 &   9.9-7 $|$    9.9-7 $|$    9.9-7	 &   -6.5-7 $|$    -3.8-7 $|$    5.4-9 	 &52 $|$ 1:17 $|$ 3:03\\[2pt] 
 \hline 
be250.8 
	 &251  $;$  251 	 &10 	 &3507  $|$  4712  $|$  7701 	 &   9.9-7 $|$    9.9-7 $|$    9.6-7	 &   -9.7-7 $|$    {\green{ -1.0-6}} $|$    5.1-7 	 &50 $|$ 56 $|$ 1:43\\[2pt] 
 \hline 
be250.9 
	 &251  $;$  251 	 &10 	 &3678  $|$  7292  $|$  21001 	 &   9.9-7 $|$    9.9-7 $|$    9.9-7	 &   -4.1-7 $|$    -7.2-7 $|$    -1.2-8 	 &53 $|$ 1:28 $|$ 4:57\\[2pt] 
 \hline 
be250.10 
	 &251  $;$  251 	 &10 	 &3305  $|$  5752  $|$  10500 	 &   9.9-7 $|$    9.9-7 $|$    9.9-7	 &   {\green{ -1.1-6}} $|$    -8.2-7 $|$    -3.7-8 	 &49 $|$ 1:06 $|$ 2:19\\[2pt] 
 \hline 
bqp100-1 
	 &101  $;$  101 	 &10 	 &1376  $|$  2134  $|$  3067 	 &   9.9-7 $|$    9.9-7 $|$    9.9-7	 &   2.6-7 $|$    -1.9-7 $|$    -5.1-7 	 &05 $|$ 06 $|$ 10\\[2pt] 
 \hline 
bqp100-2 
	 &101  $;$  101 	 &10 	 &3109  $|$  4319  $|$  7107 	 &   9.9-7 $|$    9.9-7 $|$    9.9-7	 &   -1.8-7 $|$    -7.2-7 $|$    -5.3-7 	 &10 $|$ 13 $|$ 22\\[2pt] 
 \hline 
bqp100-3 
	 &101  $;$  101 	 &10 	 &1751  $|$  2371  $|$  6276 	 &   9.9-7 $|$    9.9-7 $|$    9.9-7	 &   {\green{ -2.7-6}} $|$    {\green{ -3.1-6}} $|$    4.7-7 	 &06 $|$ 06 $|$ 20\\[2pt] 
 \hline 
bqp100-4 
	 &101  $;$  101 	 &10 	 &2646  $|$  3986  $|$  13901 	 &   9.9-7 $|$    9.9-7 $|$    9.1-7	 &   -4.0-7 $|$    -6.6-7 $|$    -3.3-8 	 &09 $|$ 11 $|$ 45\\[2pt] 
 \hline 
bqp100-5 
	 &101  $;$  101 	 &10 	 &1979  $|$  3001  $|$  6901 	 &   9.9-7 $|$    9.9-7 $|$    9.7-7	 &   -3.7-7 $|$    -1.5-7 $|$    1.7-8 	 &07 $|$ 08 $|$ 22\\[2pt] 
 \hline 
bqp100-6 
	 &101  $;$  101 	 &10 	 &1316  $|$  2083  $|$  2937 	 &   9.4-7 $|$    9.9-7 $|$    9.9-7	 &   1.1-7 $|$    3.3-7 $|$    -9.5-7 	 &05 $|$ 06 $|$ 11\\[2pt] 
 \hline 
bqp100-7 
	 &101  $;$  101 	 &10 	 &1787  $|$  2341  $|$  3664 	 &   9.9-7 $|$    9.9-7 $|$    9.9-7	 &   -5.5-7 $|$    -5.1-7 $|$    {\green{ -1.3-6}} 	 &06 $|$ 06 $|$ 12\\[2pt] 
 \hline 
bqp100-8 
	 &101  $;$  101 	 &10 	 &1820  $|$  3337  $|$  9612 	 &   9.9-7 $|$    9.9-7 $|$    9.9-7	 &   7.3-7 $|$    8.9-8 $|$    1.1-8 	 &06 $|$ 09 $|$ 32\\[2pt] 
 \hline 
bqp100-9 
	 &101  $;$  101 	 &10 	 &1948  $|$  4146  $|$  15901 	 &   9.9-7 $|$    9.9-7 $|$    9.9-7	 &   {\green{ -2.2-6}} $|$    -6.7-7 $|$    2.6-9 	 &07 $|$ 11 $|$ 52\\[2pt] 
 \hline 
bqp100-10 
	 &101  $;$  101 	 &10 	 &3207  $|$  5077  $|$  12101 	 &   9.9-7 $|$    9.9-7 $|$    9.9-7	 &   8.0-8 $|$    4.3-7 $|$    2.7-8 	 &10 $|$ 15 $|$ 38\\[2pt] 
 \hline 
bqp250-1 
	 &251  $;$  251 	 &10 	 &3931  $|$  5941  $|$  11758 	 &   9.6-7 $|$    9.9-7 $|$    9.9-7	 &   {\green{ -1.2-6}} $|$    {\green{ -1.5-6}} $|$    1.2-7 	 &57 $|$ 1:10 $|$ 2:39\\[2pt] 
 \hline 
bqp250-2 
	 &251  $;$  251 	 &10 	 &4007  $|$  5774  $|$  9704 	 &   9.5-7 $|$    9.9-7 $|$    9.9-7	 &   -6.6-7 $|$    -2.3-7 $|$    {\green{ -1.2-6}} 	 &57 $|$ 1:07 $|$ 2:11\\[2pt] 
 \hline 
bqp250-3 
	 &251  $;$  251 	 &10 	 &4112  $|$  5708  $|$  12202 	 &   9.9-7 $|$    9.9-7 $|$    9.9-7	 &   {\green{ -3.9-6}} $|$    3.8-8 $|$    {\green{ 3.0-6}} 	 &57 $|$ 1:05 $|$ 2:40\\[2pt] 
 \hline 
bqp250-4 
	 &251  $;$  251 	 &10 	 &3158  $|$  4290  $|$  9671 	 &   9.9-7 $|$    9.9-7 $|$    9.9-7	 &   -5.5-7 $|$    {\green{ -2.4-6}} $|$    {\green{ 4.5-6}} 	 &45 $|$ 52 $|$ 2:13\\[2pt] 
 \hline 
bqp250-5 
	 &251  $;$  251 	 &10 	 &4430  $|$  7349  $|$  22802 	 &   9.9-7 $|$    9.9-7 $|$    9.9-7	 &   {\green{ -2.0-6}} $|$    {\green{ 3.6-6}} $|$    -1.3-8 	 &1:02 $|$ 1:29 $|$ 5:13\\[2pt] 
 \hline 
bqp250-6 
	 &251  $;$  251 	 &10 	 &2871  $|$  5122  $|$  7801 	 &   9.9-7 $|$    9.9-7 $|$    9.9-7	 &   {\green{ -1.2-6}} $|$    {\green{ -1.3-6}} $|$    -2.5-7 	 &42 $|$ 1:01 $|$ 1:47\\[2pt] 
 \hline 
bqp250-7 
	 &251  $;$  251 	 &10 	 &3991  $|$  5570  $|$  11508 	 &   9.9-7 $|$    9.9-7 $|$    9.9-7	 &   {\green{ -2.2-6}} $|$    {\green{ -2.0-6}} $|$    {\green{ -2.7-6}} 	 &57 $|$ 1:04 $|$ 2:31\\[2pt] 
 \hline 
bqp250-8 
	 &251  $;$  251 	 &10 	 &2882  $|$  4008  $|$  5501 	 &   9.9-7 $|$    9.8-7 $|$    9.8-7	 &   -2.0-7 $|$    -7.1-7 $|$    {\green{ -1.0-6}} 	 &40 $|$ 45 $|$ 1:14\\[2pt] 
 \hline 
bqp250-9 
	 &251  $;$  251 	 &10 	 &4127  $|$  6279  $|$  11998 	 &   9.7-7 $|$    9.9-7 $|$    9.9-7	 &   -5.1-7 $|$    -3.9-7 $|$    {\green{ 3.8-6}} 	 &58 $|$ 1:11 $|$ 2:38\\[2pt] 
 \hline 
bqp250-10 
	 &251  $;$  251 	 &10 	 &3044  $|$  4185  $|$  7986 	 &   9.9-7 $|$    9.9-7 $|$    9.9-7	 &   -9.3-7 $|$    -7.5-7 $|$    {\green{ -2.5-6}} 	 &43 $|$ 48 $|$ 1:43\\[2pt] 
 \hline 
bqp500-1 
	 &501  $;$  501 	 &10 	 &6003  $|$  8391  $|$  13416 	 &   9.9-7 $|$    9.9-7 $|$    9.9-7	 &   -3.9-7 $|$    -7.3-7 $|$    -5.4-7 	 &6:01 $|$ 7:05 $|$ 13:34\\[2pt] 
 \hline 
bqp500-2 
	 &501  $;$  501 	 &10 	 &6601  $|$  10203  $|$  25000 	 &   9.7-7 $|$    9.9-7 $|$    {\green{ 3.4-6}}	 &   -4.2-7 $|$    -1.2-7 $|$    {\red{ 1.8-5}} 	 &6:52 $|$ 8:43 $|$ 25:23\\[2pt] 
 \hline 
bqp500-3 
	 &501  $;$  501 	 &10 	 &7450  $|$  10517  $|$  21140 	 &   9.9-7 $|$    9.9-7 $|$    9.9-7	 &   7.6-7 $|$    {\green{ -4.3-6}} $|$    {\green{ 1.1-6}} 	 &7:31 $|$ 8:46 $|$ 21:10\\[2pt] 
 \hline 
bqp500-4 
	 &501  $;$  501 	 &10 	 &7035  $|$  9903  $|$  23551 	 &   9.6-7 $|$    9.9-7 $|$    9.9-7	 &   -3.3-7 $|$    {\green{ -1.3-6}} $|$    {\green{ 2.6-6}} 	 &7:08 $|$ 8:12 $|$ 23:36\\[2pt] 
 \hline 
bqp500-5 
	 &501  $;$  501 	 &10 	 &6164  $|$  8406  $|$  20533 	 &   9.9-7 $|$    9.9-7 $|$    9.9-7	 &   -8.8-7 $|$    -4.8-7 $|$    {\green{ 2.8-6}} 	 &6:30 $|$ 7:04 $|$ 20:37\\[2pt] 
 \hline 
bqp500-6 
	 &501  $;$  501 	 &10 	 &6905  $|$  8659  $|$  25000 	 &   9.8-7 $|$    9.9-7 $|$    {\bf 1.4-4}	 &   -3.8-7 $|$    {\green{ -1.5-6}} $|$    {\bf -1.8-4} 	 &7:13 $|$ 7:30 $|$ 25:44\\[2pt] 
 \hline 
bqp500-7 
	 &501  $;$  501 	 &10 	 &6587  $|$  9038  $|$  18072 	 &   9.9-7 $|$    9.9-7 $|$    9.9-7	 &   -6.8-7 $|$    2.5-7 $|$    {\green{ 2.8-6}} 	 &6:41 $|$ 7:39 $|$ 18:13\\[2pt] 
 \hline 
bqp500-8 
	 &501  $;$  501 	 &10 	 &6300  $|$  8832  $|$  16496 	 &   9.9-7 $|$    9.9-7 $|$    9.9-7	 &   {\green{ 1.3-6}} $|$    {\green{ -1.6-6}} $|$    {\blue{ 5.8-6}} 	 &6:24 $|$ 7:17 $|$ 16:20\\[2pt] 
 \hline 
bqp500-9 
	 &501  $;$  501 	 &10 	 &6532  $|$  9015  $|$  18065 	 &   9.9-7 $|$    9.9-7 $|$    9.9-7	 &   9.9-7 $|$    -6.5-7 $|$    {\green{ -3.5-6}} 	 &6:39 $|$ 7:37 $|$ 18:10\\[2pt] 
 \hline 
bqp500-10 
	 &501  $;$  501 	 &10 	 &7199  $|$  9787  $|$  24119 	 &   9.9-7 $|$    9.9-7 $|$    9.9-7	 &   {\green{ -1.9-6}} $|$    {\green{ 2.1-6}} $|$    {\green{ -2.3-6}} 	 &7:09 $|$ 8:12 $|$ 24:15\\[2pt] 
 \hline 
gka1d 
	 &101  $;$  101 	 &10 	 &1600  $|$  2266  $|$  4068 	 &   9.8-7 $|$    9.9-7 $|$    9.7-7	 &   -4.2-7 $|$    -8.8-7 $|$    7.4-7 	 &06 $|$ 06 $|$ 13\\[2pt] 
 \hline 
gka2d 
	 &101  $;$  101 	 &10 	 &1903  $|$  3097  $|$  5601 	 &   9.9-7 $|$    9.9-7 $|$    9.3-7	 &   -5.9-7 $|$    -2.4-7 $|$    -3.8-8 	 &07 $|$ 09 $|$ 21\\[2pt] 
 \hline 
gka3d 
	 &101  $;$  101 	 &10 	 &2431  $|$  3101  $|$  5618 	 &   9.9-7 $|$    9.9-7 $|$    9.9-7	 &   -2.6-7 $|$    -3.8-7 $|$    1.7-8 	 &08 $|$ 09 $|$ 19\\[2pt] 
 \hline 
gka4d 
	 &101  $;$  101 	 &10 	 &2266  $|$  2787  $|$  6632 	 &   9.9-7 $|$    9.9-7 $|$    9.9-7	 &   2.3-7 $|$    -4.4-7 $|$    -1.9-8 	 &08 $|$ 09 $|$ 22\\[2pt] 
 \hline 
soybean-large-2 
	 &308  $;$  307 	 &10 	 &1267  $|$  1717  $|$  11208 	 &   9.9-7 $|$    9.9-7 $|$    9.9-7	 &   -5.8-8 $|$    -6.5-8 $|$    -7.9-8 	 &20 $|$ 23 $|$ 2:55\\[2pt] 
 \hline 
soybean-large-3 
	 &308  $;$  307 	 &10 	 &936  $|$  1362  $|$  9261 	 &   8.3-7 $|$    9.1-7 $|$    9.8-7	 &   -5.1-8 $|$    -5.7-8 $|$    -1.7-8 	 &17 $|$ 17 $|$ 2:29\\[2pt] 
 \hline 
soybean-large-4 
	 &308  $;$  307 	 &10 	 &1681  $|$  2132  $|$  13401 	 &   9.9-7 $|$    9.9-7 $|$    9.9-7	 &   -1.0-7 $|$    -1.0-7 $|$    -4.3-8 	 &29 $|$ 28 $|$ 3:49\\[2pt] 
 \hline 
soybean-large-5 
	 &308  $;$  307 	 &10 	 &834  $|$  1229  $|$  3937 	 &   9.9-7 $|$    9.9-7 $|$    9.9-7	 &   -3.2-8 $|$    -1.9-8 $|$    -2.3-8 	 &14 $|$ 18 $|$ 1:08\\[2pt] 
 \hline 
soybean-large-6 
	 &308  $;$  307 	 &10 	 &310  $|$  475  $|$  707 	 &   9.4-7 $|$    8.9-7 $|$    8.3-7	 &   -8.1-8 $|$    -5.8-8 $|$    -1.5-7 	 &05 $|$ 06 $|$ 12\\[2pt] 
 \hline 
soybean-large-7 
	 &308  $;$  307 	 &10 	 &1028  $|$  1327  $|$  3970 	 &   9.9-7 $|$    9.9-7 $|$    9.9-7	 &   -3.6-8 $|$    -6.3-8 $|$    -1.8-8 	 &19 $|$ 20 $|$ 1:12\\[2pt] 
 \hline 
soybean-large-8 
	 &308  $;$  307 	 &10 	 &782  $|$  1091  $|$  2901 	 &   9.8-7 $|$    9.9-7 $|$    8.9-7	 &   -3.7-8 $|$    -4.5-8 $|$    -1.0-8 	 &14 $|$ 15 $|$ 51\\[2pt] 
 \hline 
soybean-large-9 
	 &308  $;$  307 	 &10 	 &928  $|$  1187  $|$  4901 	 &   9.8-7 $|$    9.8-7 $|$    9.9-7	 &   1.1-7 $|$    -6.0-8 $|$    -1.7-8 	 &17 $|$ 19 $|$ 1:26\\[2pt] 
 \hline 
soybean-large-10 
	 &308  $;$  307 	 &10 	 &309  $|$  489  $|$  518 	 &   9.9-7 $|$    9.9-7 $|$    9.7-7	 &   2.0-7 $|$    3.1-7 $|$    1.4-7 	 &06 $|$ 07 $|$ 09\\[2pt] 
 \hline 
soybean-large-11 
	 &308  $;$  307 	 &10 	 &877  $|$  1605  $|$  1755 	 &   9.9-7 $|$    8.6-7 $|$    9.5-7	 &   -2.2-7 $|$    3.5-7 $|$    -2.6-7 	 &17 $|$ 23 $|$ 32\\[2pt] 
 \hline 
spambase-small-2 
	 &301  $;$  300 	 &10 	 &409  $|$  610  $|$  2792 	 &   8.8-7 $|$    9.5-7 $|$    9.0-7	 &   -3.1-7 $|$    -3.9-7 $|$    {\green{ -1.1-6}} 	 &06 $|$ 07 $|$ 40\\[2pt] 
 \hline 
spambase-small-3 
	 &301  $;$  300 	 &10 	 &476  $|$  665  $|$  1201 	 &   9.6-7 $|$    9.9-7 $|$    9.6-7	 &   7.8-9 $|$    -3.7-8 $|$    -3.3-8 	 &09 $|$ 08 $|$ 17\\[2pt] 
 \hline 
spambase-small-4 
	 &301  $;$  300 	 &10 	 &1305  $|$  1983  $|$  6073 	 &   9.9-7 $|$    9.9-7 $|$    9.9-7	 &   -4.5-9 $|$    6.6-9 $|$    -1.7-8 	 &20 $|$ 28 $|$ 1:36\\[2pt] 
 \hline 
spambase-small-5 
	 &301  $;$  300 	 &10 	 &608  $|$  819  $|$  868 	 &   8.5-7 $|$    9.8-7 $|$    9.9-7	 &   -7.3-7 $|$    -2.7-7 $|$    -1.4-7 	 &11 $|$ 11 $|$ 14\\[2pt] 
 \hline 
spambase-small-6 
	 &301  $;$  300 	 &10 	 &811  $|$  1198  $|$  1334 	 &   9.9-7 $|$    9.9-7 $|$    9.9-7	 &   -1.5-7 $|$    -2.0-7 $|$    -1.3-7 	 &14 $|$ 17 $|$ 23\\[2pt] 
 \hline 
spambase-small-7 
	 &301  $;$  300 	 &10 	 &849  $|$  1240  $|$  1359 	 &   9.9-7 $|$    9.9-7 $|$    9.9-7	 &   4.0-7 $|$    2.8-7 $|$    1.8-7 	 &15 $|$ 18 $|$ 25\\[2pt] 
 \hline 
spambase-small-8 
	 &301  $;$  300 	 &10 	 &1109  $|$  1244  $|$  1501 	 &   9.9-7 $|$    9.9-7 $|$    8.8-7	 &   7.1-8 $|$    9.3-8 $|$    7.6-8 	 &20 $|$ 18 $|$ 27\\[2pt] 
 \hline 
spambase-small-9 
	 &301  $;$  300 	 &10 	 &1090  $|$  1415  $|$  1440 	 &   9.9-7 $|$    9.7-7 $|$    9.9-7	 &   -1.7-7 $|$    2.9-8 $|$    -1.3-8 	 &20 $|$ 21 $|$ 27\\[2pt] 
 \hline 
spambase-small-10 
	 &301  $;$  300 	 &10 	 &1081  $|$  1341  $|$  1500 	 &   9.9-7 $|$    9.9-7 $|$    9.9-7	 &   1.7-7 $|$    1.5-7 $|$    -1.5-7 	 &20 $|$ 22 $|$ 27\\[2pt] 
 \hline 
spambase-small-11 
	 &301  $;$  300 	 &10 	 &1319  $|$  1482  $|$  1653 	 &   9.9-7 $|$    9.9-7 $|$    9.9-7	 &   -3.6-7 $|$    -8.3-7 $|$    -5.8-7 	 &25 $|$ 25 $|$ 31\\[2pt] 
 \hline 
spambase-medium-2 
	 &901  $;$  900 	 &10 	 &471  $|$  596  $|$  1201 	 &   9.9-7 $|$    9.9-7 $|$    8.9-7	 &   {\green{ -1.6-6}} $|$    {\green{ -1.3-6}} $|$    {\green{ -1.9-6}} 	 &1:42 $|$ 1:37 $|$ 4:01\\[2pt] 
 \hline 
spambase-medium-3 
	 &901  $;$  900 	 &10 	 &1205  $|$  1582  $|$  11000 	 &   9.9-7 $|$    9.9-7 $|$    9.9-7	 &   -2.0-7 $|$    -1.8-7 $|$    -2.2-7 	 &4:18 $|$ 4:16 $|$ 36:54\\[2pt] 
 \hline 
spambase-medium-4 
	 &901  $;$  900 	 &10 	 &2560  $|$  2990  $|$  4045 	 &   9.7-7 $|$    9.8-7 $|$    9.9-7	 &   {\green{ -2.3-6}} $|$    {\green{ 2.5-6}} $|$    {\green{ 1.1-6}} 	 &9:06 $|$ 8:04 $|$ 13:37\\[2pt] 
 \hline 
spambase-medium-5 
	 &901  $;$  900 	 &10 	 &1414  $|$  1900  $|$  2901 	 &   9.9-7 $|$    9.9-7 $|$    9.0-7	 &   7.4-8 $|$    3.8-8 $|$    {\green{ -1.1-6}} 	 &5:06 $|$ 5:17 $|$ 9:58\\[2pt] 
 \hline 
spambase-medium-6 
	 &901  $;$  900 	 &10 	 &1607  $|$  2107  $|$  2698 	 &   9.9-7 $|$    9.9-7 $|$    9.9-7	 &   -1.0-8 $|$    3.7-8 $|$    {\green{ -1.3-6}} 	 &6:01 $|$ 6:16 $|$ 9:25\\[2pt] 
 \hline 
spambase-medium-7 
	 &901  $;$  900 	 &10 	 &1805  $|$  2508  $|$  2846 	 &   9.9-7 $|$    9.9-7 $|$    9.9-7	 &   -8.7-8 $|$    -4.5-8 $|$    {\green{ -1.4-6}} 	 &6:55 $|$ 7:36 $|$ 10:00\\[2pt] 
 \hline 
spambase-medium-8 
	 &901  $;$  900 	 &10 	 &1655  $|$  2309  $|$  2489 	 &   9.9-7 $|$    9.9-7 $|$    9.9-7	 &   -2.6-8 $|$    -6.7-8 $|$    4.6-7 	 &6:19 $|$ 6:54 $|$ 8:47\\[2pt] 
 \hline 
spambase-medium-9 
	 &901  $;$  900 	 &10 	 &1683  $|$  2330  $|$  2687 	 &   9.9-7 $|$    9.9-7 $|$    9.9-7	 &   2.6-8 $|$    -5.9-8 $|$    2.2-8 	 &6:23 $|$ 6:56 $|$ 9:38\\[2pt] 
 \hline 
spambase-medium-10 
	 &901  $;$  900 	 &10 	 &1641  $|$  2030  $|$  2617 	 &   9.9-7 $|$    9.9-7 $|$    9.8-7	 &   -6.5-7 $|$    -4.7-7 $|$    {\green{ 1.9-6}} 	 &6:11 $|$ 5:59 $|$ 9:22\\[2pt] 
 \hline 
spambase-medium-11 
	 &901  $;$  900 	 &10 	 &1608  $|$  1838  $|$  3210 	 &   9.9-7 $|$    9.9-7 $|$    9.9-7	 &   -5.0-7 $|$    5.4-7 $|$    9.0-7 	 &6:06 $|$ 5:20 $|$ 11:21\\[2pt] 
 \hline 
abalone-medium-2 
	 &401  $;$  400 	 &10 	 &500  $|$  682  $|$  1301 	 &   9.9-7 $|$    9.9-7 $|$    8.5-7	 &   -7.4-8 $|$    5.8-8 $|$    3.4-8 	 &16 $|$ 17 $|$ 40\\[2pt] 
 \hline 
abalone-medium-3 
	 &401  $;$  400 	 &10 	 &715  $|$  1011  $|$  1679 	 &   9.9-7 $|$    9.9-7 $|$    9.9-7	 &   -2.5-9 $|$    1.3-8 $|$    -1.1-8 	 &24 $|$ 28 $|$ 56\\[2pt] 
 \hline 
abalone-medium-4 
	 &401  $;$  400 	 &10 	 &372  $|$  626  $|$  684 	 &   9.9-7 $|$    9.9-7 $|$    9.9-7	 &   -5.3-8 $|$    3.6-9 $|$    6.3-9 	 &12 $|$ 16 $|$ 24\\[2pt] 
 \hline 
abalone-medium-5 
	 &401  $;$  400 	 &10 	 &524  $|$  779  $|$  942 	 &   9.9-7 $|$    9.9-7 $|$    9.9-7	 &   -3.8-8 $|$    -1.4-7 $|$    -9.6-8 	 &18 $|$ 21 $|$ 32\\[2pt] 
 \hline 
abalone-medium-6 
	 &401  $;$  400 	 &10 	 &536  $|$  946  $|$  1162 	 &   9.7-7 $|$    9.9-7 $|$    9.9-7	 &   -1.3-7 $|$    -2.3-7 $|$    -1.8-7 	 &22 $|$ 27 $|$ 38\\[2pt] 
 \hline 
abalone-medium-7 
	 &401  $;$  400 	 &10 	 &1046  $|$  1676  $|$  2013 	 &   9.9-7 $|$    9.9-7 $|$    9.9-7	 &   -8.9-8 $|$    -4.2-8 $|$    -3.3-8 	 &37 $|$ 47 $|$ 1:09\\[2pt] 
 \hline 
abalone-medium-8 
	 &401  $;$  400 	 &10 	 &745  $|$  1123  $|$  1641 	 &   9.6-7 $|$    9.7-7 $|$    9.9-7	 &   -3.9-8 $|$    -2.2-7 $|$    -9.1-8 	 &27 $|$ 32 $|$ 55\\[2pt] 
 \hline 
abalone-medium-9 
	 &401  $;$  400 	 &10 	 &1035  $|$  1504  $|$  1709 	 &   9.9-7 $|$    9.5-7 $|$    9.9-7	 &   -8.3-8 $|$    7.1-8 $|$    -1.2-8 	 &38 $|$ 43 $|$ 1:02\\[2pt] 
 \hline 
abalone-medium-10 
	 &401  $;$  400 	 &10 	 &1349  $|$  1803  $|$  1904 	 &   9.9-7 $|$    9.4-7 $|$    9.8-7	 &   -1.7-7 $|$    -2.0-7 $|$    -2.2-7 	 &49 $|$ 51 $|$ 1:07\\[2pt] 
 \hline 
abalone-medium-11 
	 &401  $;$  400 	 &10 	 &1066  $|$  1504  $|$  1704 	 &   9.9-7 $|$    9.7-7 $|$    9.5-7	 &   -1.1-7 $|$    -1.6-7 $|$    -1.6-7 	 &40 $|$ 45 $|$ 1:02\\[2pt] 
 \hline 
abalone-large-2 
	 &1001  $;$  1000 	 &10 	 &594  $|$  734  $|$  909 	 &   9.9-7 $|$    9.8-7 $|$    9.9-7	 &   4.6-7 $|$    4.5-7 $|$    1.3-7 	 &3:16 $|$ 2:35 $|$ 3:54\\[2pt] 
 \hline 
abalone-large-3 
	 &1001  $;$  1000 	 &10 	 &656  $|$  1014  $|$  1901 	 &   9.9-7 $|$    9.9-7 $|$    9.9-7	 &   -1.4-8 $|$    -7.2-8 $|$    -4.4-8 	 &3:03 $|$ 3:37 $|$ 8:20\\[2pt] 
 \hline 
abalone-large-4 
	 &1001  $;$  1000 	 &10 	 &505  $|$  749  $|$  995 	 &   9.9-7 $|$    9.9-7 $|$    9.8-7	 &   -1.3-9 $|$    -1.6-8 $|$    -6.6-8 	 &2:42 $|$ 2:39 $|$ 4:24\\[2pt] 
 \hline 
abalone-large-5 
	 &1001  $;$  1000 	 &10 	 &752  $|$  1187  $|$  1550 	 &   9.8-7 $|$    9.9-7 $|$    9.9-7	 &   -6.8-8 $|$    -1.8-7 $|$    -1.2-7 	 &4:11 $|$ 4:16 $|$ 6:53\\[2pt] 
 \hline 
abalone-large-6 
	 &1001  $;$  1000 	 &10 	 &886  $|$  1364  $|$  1670 	 &   9.9-7 $|$    9.9-7 $|$    9.9-7	 &   -9.5-8 $|$    -1.1-7 $|$    -1.2-7 	 &4:09 $|$ 4:56 $|$ 7:27\\[2pt] 
 \hline 
abalone-large-7 
	 &1001  $;$  1000 	 &10 	 &1206  $|$  1614  $|$  2251 	 &   9.9-7 $|$    9.9-7 $|$    9.9-7	 &   -1.1-7 $|$    1.8-8 $|$    -7.5-8 	 &5:40 $|$ 5:47 $|$ 9:59\\[2pt] 
 \hline 
abalone-large-8 
	 &1001  $;$  1000 	 &10 	 &1092  $|$  1721  $|$  2046 	 &   9.9-7 $|$    9.9-7 $|$    9.9-7	 &   -3.1-7 $|$    -1.8-7 $|$    -2.9-7 	 &5:08 $|$ 6:14 $|$ 9:07\\[2pt] 
 \hline 
abalone-large-9 
	 &1001  $;$  1000 	 &10 	 &1557  $|$  2407  $|$  2746 	 &   9.8-7 $|$    9.9-7 $|$    9.9-7	 &   -3.8-7 $|$    -3.5-7 $|$    -2.8-7 	 &8:30 $|$ 8:47 $|$ 12:15\\[2pt] 
 \hline 
abalone-large-10 
	 &1001  $;$  1000 	 &10 	 &1682  $|$  2488  $|$  2821 	 &   9.9-7 $|$    9.9-7 $|$    9.9-7	 &   -1.6-7 $|$    -2.6-7 $|$    -2.5-7 	 &8:00 $|$ 9:06 $|$ 12:39\\[2pt] 
 \hline 
abalone-large-11 
	 &1001  $;$  1000 	 &10 	 &1923  $|$  3005  $|$  3723 	 &   9.8-7 $|$    9.9-7 $|$    9.9-7	 &   1.3-7 $|$    3.6-8 $|$    -3.5-8 	 &9:17 $|$ 11:00 $|$ 16:39\\[2pt] 
 \hline 
segment-medium-2 
	 &701  $;$  700 	 &10 	 &1016  $|$  1541  $|$  1880 	 &   9.7-7 $|$    9.8-7 $|$    9.9-7	 &   {\green{ 1.3-6}} $|$    {\green{ -1.1-6}} $|$    2.5-7 	 &2:07 $|$ 2:13 $|$ 3:26\\[2pt] 
 \hline 
segment-medium-3 
	 &701  $;$  700 	 &10 	 &713  $|$  714  $|$  1801 	 &   9.4-7 $|$    9.5-7 $|$    9.2-7	 &   -4.0-7 $|$    -9.7-7 $|$    -8.7-7 	 &1:24 $|$ 1:03 $|$ 3:20\\[2pt] 
 \hline 
segment-medium-4 
	 &701  $;$  700 	 &10 	 &2282  $|$  2710  $|$  17881 	 &   9.9-7 $|$    9.9-7 $|$    9.9-7	 &   -7.1-8 $|$    -6.5-8 $|$    -6.5-8 	 &4:30 $|$ 4:25 $|$ 34:11\\[2pt] 
 \hline 
segment-medium-5 
	 &701  $;$  700 	 &10 	 &2322  $|$  3100  $|$  18701 	 &   9.9-7 $|$    9.9-7 $|$    9.9-7	 &   -1.2-7 $|$    -9.5-8 $|$    -7.3-8 	 &4:40 $|$ 5:02 $|$ 35:56\\[2pt] 
 \hline 
segment-medium-6 
	 &701  $;$  700 	 &10 	 &2966  $|$  3916  $|$  25000 	 &   9.9-7 $|$    9.9-7 $|$    {\green{ 1.4-6}}	 &   -1.7-7 $|$    -1.4-7 $|$    -1.3-7 	 &6:12 $|$ 6:29 $|$ 51:26\\[2pt] 
 \hline 
segment-medium-7 
	 &701  $;$  700 	 &10 	 &3185  $|$  4268  $|$  25000 	 &   9.9-7 $|$    9.9-7 $|$    {\green{ 1.6-6}}	 &   -1.7-7 $|$    -1.7-7 $|$    -1.6-7 	 &7:03 $|$ 7:34 $|$ 53:28\\[2pt] 
 \hline 
segment-medium-8 
	 &701  $;$  700 	 &10 	 &2998  $|$  4140  $|$  25000 	 &   9.9-7 $|$    9.9-7 $|$    {\green{ 1.1-6}}	 &   -1.6-7 $|$    -1.7-7 $|$    -6.7-8 	 &6:28 $|$ 7:09 $|$ 52:54\\[2pt] 
 \hline 
segment-medium-9 
	 &701  $;$  700 	 &10 	 &2123  $|$  2635  $|$  8801 	 &   9.9-7 $|$    9.9-7 $|$    9.9-7	 &   -1.9-7 $|$    -3.0-8 $|$    -4.3-8 	 &4:32 $|$ 4:25 $|$ 18:04\\[2pt] 
 \hline 
segment-medium-10 
	 &701  $;$  700 	 &10 	 &1695  $|$  2414  $|$  6101 	 &   9.9-7 $|$    9.9-7 $|$    9.8-7	 &   -2.4-7 $|$    -1.2-7 $|$    -2.2-8 	 &3:35 $|$ 4:07 $|$ 12:27\\[2pt] 
 \hline 
segment-medium-11 
	 &701  $;$  700 	 &10 	 &1454  $|$  2437  $|$  2101 	 &   9.4-7 $|$    9.7-7 $|$    8.6-7	 &   6.4-8 $|$    -6.3-7 $|$    -1.5-7 	 &3:01 $|$ 4:00 $|$ 4:13\\[2pt] 
 \hline 
segment-large-2 
	 &1001  $;$  1000 	 &10 	 &1348  $|$  1823  $|$  2038 	 &   9.6-7 $|$    9.9-7 $|$    9.9-7	 &   {\green{ -1.3-6}} $|$    {\green{ -1.3-6}} $|$    {\green{ -1.4-6}} 	 &6:30 $|$ 6:15 $|$ 8:40\\[2pt] 
 \hline 
segment-large-3 
	 &1001  $;$  1000 	 &10 	 &479  $|$  533  $|$  1601 	 &   9.9-7 $|$    9.9-7 $|$    8.7-7	 &   -4.0-7 $|$    {\green{ -1.0-6}} $|$    -4.4-7 	 &2:10 $|$ 1:53 $|$ 6:49\\[2pt] 
 \hline 
segment-large-4 
	 &1001  $;$  1000 	 &10 	 &2157  $|$  2802  $|$  20226 	 &   9.9-7 $|$    9.9-7 $|$    9.9-7	 &   -9.1-8 $|$    -9.5-8 $|$    -7.1-8 	 &9:57 $|$ 9:57 $|$   1:27:58\\[2pt] 
 \hline 
segment-large-5 
	 &1001  $;$  1000 	 &10 	 &2618  $|$  3404  $|$  25000 	 &   9.9-7 $|$    9.9-7 $|$    {\green{ 1.0-6}}	 &   -1.1-7 $|$    -9.3-8 $|$    -8.3-8 	 &12:13 $|$ 12:12 $|$   1:50:29\\[2pt] 
 \hline 
segment-large-6 
	 &1001  $;$  1000 	 &10 	 &3236  $|$  4143  $|$  25000 	 &   9.9-7 $|$    9.9-7 $|$    {\green{ 1.4-6}}	 &   -1.8-7 $|$    -1.8-7 $|$    -1.2-7 	 &15:28 $|$ 15:20 $|$   1:52:58\\[2pt] 
 \hline 
segment-large-7 
	 &1001  $;$  1000 	 &10 	 &3505  $|$  4318  $|$  25000 	 &   9.9-7 $|$    9.9-7 $|$    {\green{ 1.8-6}}	 &   -1.8-7 $|$    -1.7-7 $|$    -1.9-7 	 &17:07 $|$ 16:39 $|$   1:56:00\\[2pt] 
 \hline 
segment-large-8 
	 &1001  $;$  1000 	 &10 	 &3063  $|$  3749  $|$  25000 	 &   9.9-7 $|$    9.9-7 $|$    {\green{ 1.2-6}}	 &   -9.3-8 $|$    -7.8-8 $|$    -1.0-7 	 &14:55 $|$ 14:18 $|$   1:56:05\\[2pt] 
 \hline 
segment-large-9 
	 &1001  $;$  1000 	 &10 	 &2497  $|$  3248  $|$  15649 	 &   9.9-7 $|$    9.9-7 $|$    9.9-7	 &   -1.4-7 $|$    -1.2-7 $|$    -5.1-8 	 &12:05 $|$ 13:16 $|$   1:11:25\\[2pt] 
 \hline 
segment-large-10 
	 &1001  $;$  1000 	 &10 	 &1723  $|$  2226  $|$  4901 	 &   9.9-7 $|$    9.9-7 $|$    9.9-7	 &   7.4-9 $|$    1.4-8 $|$    -2.1-8 	 &8:00 $|$ 8:12 $|$ 21:45\\[2pt] 
 \hline 
segment-large-11 
	 &1001  $;$  1000 	 &10 	 &1571  $|$  2331  $|$  3417 	 &   9.9-7 $|$    9.7-7 $|$    9.9-7	 &   1.9-7 $|$    -5.1-7 $|$    -1.7-8 	 &7:20 $|$ 8:30 $|$ 15:23\\[2pt] 
 \hline 
housing-2 
	 &507  $;$  506 	 &10 	 &3183  $|$  5358  $|$  4689 	 &   9.4-7 $|$    9.7-7 $|$    9.7-7	 &   -1.9-7 $|$    1.8-7 $|$    2.0-7 	 &2:54 $|$ 3:22 $|$ 3:48\\[2pt] 
 \hline 
housing-3 
	 &507  $;$  506 	 &10 	 &845  $|$  1970  $|$  1714 	 &   9.9-7 $|$    9.9-7 $|$    9.9-7	 &   -1.5-7 $|$    1.2-7 $|$    -2.2-8 	 &48 $|$ 1:16 $|$ 1:24\\[2pt] 
 \hline 
housing-4 
	 &507  $;$  506 	 &10 	 &805  $|$  1742  $|$  2057 	 &   9.4-7 $|$    9.9-7 $|$    9.9-7	 &   -2.5-8 $|$    -4.8-8 $|$    -3.4-8 	 &45 $|$ 1:09 $|$ 1:45\\[2pt] 
 \hline 
housing-5 
	 &507  $;$  506 	 &10 	 &874  $|$  1262  $|$  1774 	 &   9.9-7 $|$    9.9-7 $|$    9.9-7	 &   2.4-7 $|$    -2.3-7 $|$    -2.6-7 	 &1:10 $|$ 1:14 $|$ 3:08\\[2pt] 
 \hline 
housing-6 
	 &507  $;$  506 	 &10 	 &586  $|$  826  $|$  1005 	 &   9.9-7 $|$    9.9-7 $|$    9.9-7	 &   -1.9-8 $|$    2.9-9 $|$    -8.6-8 	 &1:41 $|$ 1:26 $|$ 1:39\\[2pt] 
 \hline 
housing-7 
	 &507  $;$  506 	 &10 	 &583  $|$  906  $|$  1069 	 &   9.9-7 $|$    9.9-7 $|$    9.9-7	 &   -1.3-7 $|$    -2.7-7 $|$    -1.7-7 	 &32 $|$ 37 $|$ 56\\[2pt] 
 \hline 
housing-8 
	 &507  $;$  506 	 &10 	 &682  $|$  904  $|$  1074 	 &   9.9-7 $|$    9.3-7 $|$    9.9-7	 &   -1.1-7 $|$    -6.9-9 $|$    -6.6-8 	 &39 $|$ 38 $|$ 59\\[2pt] 
 \hline 
housing-9 
	 &507  $;$  506 	 &10 	 &765  $|$  1208  $|$  1590 	 &   8.5-7 $|$    9.9-7 $|$    9.8-7	 &   -1.5-7 $|$    -1.3-8 $|$    8.5-8 	 &44 $|$ 53 $|$ 1:26\\[2pt] 
 \hline 
housing-10 
	 &507  $;$  506 	 &10 	 &1027  $|$  1381  $|$  1541 	 &   9.9-7 $|$    9.9-7 $|$    9.9-7	 &   -6.4-8 $|$    -1.6-7 $|$    -1.0-7 	 &58 $|$ 1:02 $|$ 1:27\\[2pt] 
 \hline 
housing-11 
	 &507  $;$  506 	 &10 	 &867  $|$  1327  $|$  1359 	 &   9.9-7 $|$    9.9-7 $|$    9.9-7	 &   -1.0-7 $|$    -9.0-8 $|$    -9.2-8 	 &49 $|$ 1:01 $|$ 1:19\\[2pt] 
 \hline